\documentclass[11pt]{article}
\usepackage{dina4}
\usepackage{amssymb, amsmath,amsthm}
\usepackage{graphicx}
\usepackage{multirow}
\usepackage{float}
\usepackage{diagbox}
\usepackage{subfig}
\usepackage{color}
\usepackage[margin=1in]{geometry}
\usepackage{lscape}

\newcommand{\p}{\partial}

\newcommand{\Og}{\Omega}
\newcommand{\fl}[2]{\frac{#1}{#2}}
\newcommand{\dt}{\delta}

\newcommand{\nn}{\nonumber}
\newcommand{\ap}{\alpha}

\newcommand{\Dt}{\Delta}

\newcommand{\be}{\begin{equation}}
\newcommand{\ee}{\end{equation}}
\newcommand{\ba}{\begin{array}}
\newcommand{\ea}{\end{array}}
\newcommand{\bea}{\begin{eqnarray}}
\newcommand{\eea}{\end{eqnarray}}
\newcommand{\beas}{\begin{eqnarray*}}
\newcommand{\eeas}{\end{eqnarray*}}

\newcommand{\bx}{{\bf x} }

\title{A comparative study on nonlocal diffusion operators related to the fractional Laplacian}
\author{
Siwei Duo\thanks{Department of Mathematics and Statistics, Missouri University of Science and Technology, Rolla, MO 65409-0020 (Email: sddy9@mst.edu)},  \ \ 
Hong Wang\thanks{Department of Mathematics, University of South Carolina, Columbia, SC 29208 (Email: hwang@math.sc.edu) },  \ \ 
Yanzhi Zhang\thanks{Department of Mathematics and Statistics, Missouri University of Science and Technology, Rolla, MO 65409-0020 (Email: zhangyanz@mst.edu)}}

\begin{document}
\date{}
\maketitle

\begin{abstract} 
In this paper, we study four nonlocal diffusion operators, including the fractional Laplacian, spectral fractional Laplacian, regional fractional Laplacian,  and  peridynamic operator. 
These operators represent the infinitesimal generators of different stochastic processes, and especially their differences on a bounded domain are significant.
We provide extensive numerical experiments  to understand and compare their differences. 
We find that these four operators collapse to the classical Laplace operator as $\ap \to 2$.
The eigenvalues and eigenfunctions of these four operators are different, and the $k$-th (for $k \in {\mathbb N}$) eigenvalue of the spectral fractional Laplacian is always larger than those of the fractional Laplacian and regional fractional Laplacian. 
 For any $\ap \in (0, 2)$, the peridynamic operator can provide a good approximation to the fractional Laplacian, if the horizon size  $\dt$ is sufficiently large.
We  find that the solution of the peridynamic model converges to that of the fractional Laplacian model at a rate of ${\mathcal O}(\dt^{-\ap})$.
In contrast, although the regional fractional Laplacian can be used to approximate the fractional Laplacian as $\ap \to 2$, it generally provides inconsistent result from that of the fractional Laplacian if $\ap \ll 2$.
Moreover,  some conjectures are made from our numerical results, which could contribute to the mathematics analysis on these operators. 
\end{abstract}


\section{Introduction}
\setcounter{equation}{0}
\label{section1}

In the last couple of decades, nonlocal or fractional differential models have become a powerful tool for modeling challenging phenomena including anomalous transport, long-range interactions, or from local to nonlocal dynamics, which cannot be described properly by integer-order partial differential equations. 
So far, numerous fractional differential models have been proposed, among which  models with the fractional Laplacian have been well applied.
The fractional Laplacian operator $(-\Dt)^{\ap/2}$,  representing the infinitesimal generator of a symmetric $\alpha$-stable L\'{e}vy process,  has been  used to model anomalous diffusion or dispersion \cite{DuoJuZhang0016, Delia2013},   turbulent flows \cite{Carreras2001, Shlesinger1987},  systems of stochastic dynamics \cite{Bogdan2003, Chen2005}, finance \cite{Cont2004}, and so on. 
Nevertheless, the fractional Laplacian  is a nonlocal operator defined on the entire space. Moreover, an equation involving the fractional Laplacian has to be enclosed by a nonconventional, nonlocal boundary condition imposed on the complement of the physical domain where the governing equation is defined. 
The combination of these two facts introduces many significant challenges in the mathematical modeling, numerical simulations, and corresponding mathematical analysis. 
These issues have not been encountered in the context of integer-order partial differential equations. 

To avoid evaluation and analysis over the entire space,  one common approach in the literature is to ``{\it truncate}" and approximate the integral of the fractional Laplacian. 
Hence, some other nonlocal operators that are closely related to the fractional Laplacian have been proposed in recent years, including the regional fractional Laplacian, the spectral fractional Laplacian,  and the peridynamic operator. 
Similar to the fractional Laplacian, these operators are nonlocal, and on the entire space they are equivalent to the fractional Laplacian.  
On a bounded domain,  the spectral fractional Laplacian is defined via the spectra of the classical Laplace operator on the same domain \cite{Chen2005, Abatangelo0015}.
Thus, local boundary conditions are imposed to the equations involving the spectral fractional Laplacian. 
In contrast to the fractional Laplacian,  the integration domain in the regional fractional Laplacian is reduced from the entire space to a  finite domain where the governing equation is defined \cite{Guan2005, Guan2006, Chen2010}. 
Consequently, the regional fractional Laplacian model considerably simplifies the numerical computations of the original nonlocal problem with the fractional Laplacian.

Peridynamic model was originally proposed as a reformulation of the classical solid mechanics \cite{Silling2000}. 
The classical theory of solid mechanics assumes that all internal forces act through zero distance, and  the corresponding mathematical models are expressed in terms of partial differential equations, which cannot describe problems with spontaneous formation of discontinuities and other singularities. The peridynamic model leads to a nonlocal framework that does not explicitly involve the notion of deformation gradients and thus provides a more accurate description of problems with discontinuities and  singularities. 
The peridynamic operator has been recently used to approximate the fractional Laplacian so as to reduce the computational costs \cite{Delia2013, Guan2017}.

In this paper, we study and compare the properties of the fractional Laplacian, spectral fractional Laplacian, regional fractional Laplacian, and the peridynamic operator. 
We show that the spectral fractional Laplacian and regional fractional Laplacian are significantly different from the Dirichlet fractional Laplacian, although they can be used to approximate each other as the power $\ap\to 2$.
In contrast, the peridynamic operator can provide a consistent approximation to the fractional Laplacian for any $\ap \in (0, 2)$, but a large horizon size is required to obtain a good approximation,  especially when $\ap$ is small. 
The rest of this paper is organized as follows. 
In Section \ref{section2}, we introduce  the fractional Laplacian and its related nonlocal diffusion operators, including  the spectral fractional Laplacian,  regional fractional Laplacian, and  peridynamic operator. 
In Section \ref{section3}, we numerically compare these four operators by studying their nonlocal effects, eigenvalues and eigenfunctions, solutions to their corresponding Poisson equations, and dynamics of the nonlocal diffusion equations.
Finally, we draw some conclusions in Section \ref{section4}.

\section{Nonlocal diffusion operators}
\setcounter{equation}{0}
\label{section2}
In this section, we introduce the fractional Laplacian and its related  nonlocal diffusion operators, including  the spectral fractional Laplacian,  regional fractional Laplacian, and  peridynamic operator.
The peridynamic operator with a specially chosen kernel function can be used to approximate the fractional Laplacian \cite{Delia2013, Guan2017}. 
If a bounded domain is considered, the spectral fractional Laplacian and regional fractional Laplacian are significantly different from the fractional Laplacian, although they are freely interchanged with the fractional Laplacian in some  literature. 
In the following, we will introduce and compare the properties of these four operators from various aspects. 
Let $\Og \subset {\mathbb R}^n$ ($n = 1,2$, or $3$) denote an open bounded domain,  and $\Og^c = {\mathbb R}^n\backslash\Og$ represents the complement of $\Og$.

\subsection{Fractional Laplacian}
\label{section2-1}
On the entire space ${\mathbb R}^n$, the {\it fractional Laplacian} $(-\Dt)^{\ap/2}$ is defined via a pseudo-differential operator with symbol $|{\bf \xi}|^\alpha$  \cite{Landkof1972, Samko}:
\begin{eqnarray}
\label{pseudo}
(-\Delta)^{\alpha/2}u({\bx}) = \mathcal{F}^{-1}\big[|\xi|^\alpha \mathcal{F}[u]\big], \qquad \mbox{for} \ \ \ap > 0,
\end{eqnarray}
where $\mathcal{F}$ represents the Fourier transform, and $\mathcal{F}^{-1}$ is the inverse Fourier transform. 
From the probabilistic point of view, the fractional Laplacian $(-\Delta)^{\alpha/2}$ represents the infinitesimal generator of a symmetric $\alpha$-stable L\'{e}vy process \cite{Banuelos2004, Chen2008}. 
In a special case of $\ap = 2$, the definition in (\ref{pseudo}) reduces to the standard Laplace operator $-\Dt$.
The definition in (\ref{pseudo}) enables one to utilize the fast Fourier transform to efficiently solve problems involving the fractional Laplacian, however,  it is suitable only for problems defined either on the whole space ${\mathbb R}^n$ or on a bounded domain with periodic boundary conditions. 

In the literature,  an equivalent hypersingular integral definition of the  fractional Laplacian $(-\Dt)^{\ap/2}$ is given by \cite{Samko, Stein, Landkof1972, Du2012}:
\bea\label{riesz}
(-\Dt)^{\ap/2}u(\bx) &=& c_{n, \ap}\ {\rm P.V.} \int_{{\mathbb R}^n} \frac{u({\bf x})-u({\bf y})}{|{\bf y}-\bx|^{n+\alpha}}\,d{\bf y}, \qquad \mbox{for} \ \ \ap \in (0, 2),
\eea
where P.V. stands for the principal value, and $c_{n,\ap}$ is the normalization constant  given by
\bea\label{DefC1ap}
\displaystyle c_{n,\ap} = \fl{2^{\ap-1}\ap\,\Gamma\big(\displaystyle ({n+\ap})/{2}\big)}{\sqrt{\pi^n}\,\Gamma\big(\displaystyle 1-{\ap}/{2}\big)}
\eea
with $\Gamma(\cdot)$ denoting the Gamma function. 
In contrast to (\ref{pseudo}), the definition in (\ref{riesz}) can easily  incorporate with non-periodic bounded  domains.
Note that the integral representation in (\ref{riesz}) is defined for $0 < \ap < 2$, while the pseudo-differential definition in (\ref{pseudo}) is valid for all $\ap >0$. 
The equivalence of definitions (\ref{pseudo}) and (\ref{riesz}) for $\alpha \in (0,2)$ are studied in \cite[Proposition 3.3]{DiNezza2012}. 
More discussion can be found in \cite{DiNezza2012, Samko, Kwasnicki0015,Yolcu2015} and references therein. 

The fractional Laplacian on a bounded domain is of great interest, not only from the mathematical point of view, but also in practical applications. 
Recently, many studies  have been carried out on   {\it the Dirichlet fractional Laplacian} (also known as  the {\it restricted fractional Laplacian}), i.e., the fractional Laplacian on a bounded domain  $\Og$ with extended homogeneous Dirichlet boundary condition ($u(\bx) \equiv 0$ for $\bx\in\Og^c$).  
However, the current understanding of this topic still remains limited, and the main challenge is  from  the non-locality of the operator. In the following, we  will discuss some fundamental properties of the Dirichlet fractional Laplacian. 

Probabilistically, the Dirichlet fractional Laplacian  $(-\Dt)^{\ap/2}$  represents the infinitesimal generator of a symmetric $\alpha$-stable L\'{e}vy process that particles are killed upon leaving the domain $\Og$ \cite{Chen2008, Chen2010-1, Yolcu2015, Song2008}. 
One fundamental issue in the study of the Dirichlet fractional Laplacian is its eigenvalues and eigenfunctions. 
So far, their exact results  still remain unknown, and only some estimates and approximations can be found in the literature. 
In \cite{Chen2005}, it shows that  the $k$-th eigenvalue $\lambda_k$ (for $k\in{\mathbb N}$) of the Dirichlet fractional Laplacian on {a convex domain $\Og\subset{\mathbb R}^n$} is  bounded by:
\bea \label{Eig_n}
\fl{1}{2}\mu_k^{\ap/2} \leq \lambda_k \leq \mu_k^{\ap/2}, \qquad \mbox{for} \ \ \ap \in (0, 2), 
\eea
where $\mu_k$ represents the $k$-th eigenvalue of the  standard Dirichlet  Laplace operator $-\Dt$ on the same domain $\Og$. 
That is, the eigenvalue of the fractional Laplacian is always smaller than that of the standard Laplacian $-\Dt$. 
If a one-dimensional (i.e., $n = 1$) domain is considered, the estimates in (\ref{Eig_n}) can be improved, and sharper bounds can be found for two special cases, such as $k = 1$ and $\ap\in(0, 2)$ in \cite{Banuelos2004, Dyda2012}, and $\ap = 1$ and $k = 1, 2, 3$ in \cite{Banuelos2004}. 
More discussion on the eigenvalue bounds can be found in \cite{Kaleta2012, Yolcu2013, Yolcu2015, Frank0016} and references therein. 
Furthermore, in a one-dimensional interval $(-1, 1)$,  the asymptotic approximation of the eigenvalue $\lambda_k$ is given by \cite{Kwasnicki2012}:
\bea \label{Eig_Asy}
\lambda_k = \left(\fl{k\pi}{2}-\fl{(2-\ap)\pi}{8}\right)^\ap+ {\mathcal O}\left(\fl{2-\ap}{k\sqrt{\ap}}\right), \qquad \mbox{for} \ \ k\in{\mathbb N}.
\eea
It further shows that if $\ap \ge 1$, the eigenvalue $\lambda_k$ (for $k\in{\mathbb N}$) is simple, and the corresponding eigenfunction satisfies $\phi_k(-x) = (-1)^{k-1}\phi_k(x)$. 
Compared to the understanding of eigenvalues, the knowledge of eigenfunctions is even less.  
As shown in \cite{Servadei2014}, on a smooth bounded domain $\Og \subset {\mathbb R}^n$, the eigenfunctions $\phi_k$ (for $k\in{\mathbb N}$) are H\"{o}lder continuous up to the boundary. 
Recent numerical results on eigenvalues and eigenfunctions of the Dirichlet fractional Laplacian can be found in \cite{Duo2015}.

The fractional Poisson equation is one of the building blocks in the study of fractional PDEs. It takes the following form \cite{Dyda2012, Delia2013, Acosta0016, DuoWykZhang}:
\bea\label{fPoisson0}
(-\Dt)^{\ap/2}u(\bx) = f(\bx),  &\ & \mbox{for} \ \ \bx \in {\Og},\qquad\\
\label{fBC0}
u(\bx) = 0, && \mbox{for} \ \ \bx\in{\Og}^c.
\eea 
In (\ref{fBC0}), the extended homogeneous boundary conditions are imposed on the complement $\Og^c$, distinguishing from the classical Poisson problem where boundary conditions are given on $\p\Og$.
This difference can be explained from probabilistic interpretation of the standard  and  fractional Laplacian. 
The standard Laplace operator represents the infinitesimal generator of a Brownian motion with continuous sample paths; thus for a particle in domain $\Og$, it must leave the domain via the boundary points on $\p\Og$.  
By contrast, the fractional Laplacian is the infinitesimal generator of a symmetric $\alpha$-stable L\'{e}vy process with discontinuous sample paths; particles may ``jump" out of the domain  without touching any boundary points on $\p\Og$. 
Hence, the solution on $\Og$ can be determined by the values at $\p\Og$  in the context of classical Poisson equations but not in the context of fractional Poisson equations.

The nonlocal problem (\ref{fPoisson0})--(\ref{fBC0}) plays an important role in studying stationary behaviors of various problems. 
Its solutions can be expressed in terms of the Green's function, however, it is challenging to find the explicit form of the  Green's function. 
Some estimates of the Green's function can be found in \cite{Burcur0016, Dyda2012, Kwasnicki2012, Chen2010-1,  Chen2014}. 
{The regularity of the solution to (\ref{fPoisson0})--(\ref{fBC0}) is discussed in \cite{Acosta0016, RosOton2014, Serra2015-1,  RosOton2016-3, RosOton2016-2}.
As shown in \cite{RosOton2014}, assume that $\Og \subset {\mathbb R}^n$ is a bounded Lipschitz domain  which satisfies the exterior ball condition;  if the function $f \in L^{\infty}(\Omega)$,   the solution of the fractional Poisson equation (\ref{fPoisson0})--(\ref{fBC0}) satisfies $u\in C^{\alpha/2}(\mathbb{R}^n)$, and furthermore $\|u\|_{C^{\ap/2}({\mathbb R}^n)} \le c \|f\|_{L^\infty(\Og)}$
with $c$ a constant depending on $\Og$ and $\ap$. 
}

\subsection{Spectral fractional Laplacian}
\label{section2-2}
On a bounded domain $\Og \subset {\mathbb R}^n$, the {\it spectral fractional Laplacian} (also known as {\it the fractional power of the Dirichlet Laplacian},  or  the {\it ``Navier" fractional Laplacian}) is defined via the spectral decomposition of the standard Laplace operator  \cite{Servadei2014, Abatangelo0015, Musina2014},  i.e., 
\bea\label{fL}
(-\Dt_\Og)^{\ap/2} u(\bx) = \sum_{{k} \in {\mathbb N}} c_{k}\, \mu_{k}^{\ap/2} \phi_{k}(\bx), \qquad  \mbox{for} \ \  \ap > 0, 
\eea
where $\mu_{k}$ and $\phi_{k}$ are the ${k}$-th eigenvalue and normalized eigenfunction of the standard  {Dirichlet} Laplace operator $-\Dt$ on the domain $\Og$. 
From a probabilistic point of view, it represents the infinitesimal generator of a  subordinate killed Brownian motion, i.e., the process that first kills Brownian motion in a bounded domain $\Og$ and then subordinates it via a $\ap/2$-stable subordinator \cite{Song2003, Song2008}. 
Here, we include the domain $\Og$ in the notation $(-\Dt_\Og)^{\ap/2}$ to reflect this process and to distinguish it from the fractional Laplacian $(-\Dt)^{\ap/2}$ that was discussed in Section  \ref{section2-1}.  
Specially, if $\ap = 2$ the definition in (\ref{fL}) reduces to the standard Dirichlet Laplace operator $-\Dt$ on the domain $\Og$. 

The spectral fractional Laplacian is a nonlocal operator, and it is often used in the analysis of (partial) differential equations.  
The eigenvalues and eigenfunctions of the spectral fractional Laplacian are clearly suggested from its definition in (\ref{fL}), that is, the $k$-th eigenvalue of $(-\Dt_\Og)^{\ap/2}$ is  $\mu_k^{\ap/2}$, and the corresponding eigenfunction is $\phi_k(\bx)$. 
We remark that the spectral fractional Laplacian and the Dirichlet fractional Laplacian represent generators of different processes,  which is also reflected by their eigenvalues and eigenfunctions. 
{The eigenfunctions of the spectral fractional Laplacian are smooth up to the boundary as the boundary allows, while those of the Dirichlet fractional Laplacian are only H\"{o}lder continuous up to the boundary \cite{Servadei2014}. }
Additionally, it is easy to conclude from (\ref{Eig_n}) that the $k$-th  (for $k\in{\mathbb N}$) eigenvalue of the Dirichlet fractional Laplacian is always smaller than that of the spectral fractional Laplacian. 
 
The nonlocal Poisson problem  with the spectral fractional Laplacian reads \cite{DuoJuZhang0016}: 
\bea\label{fPoisson1}
(-\Dt_\Og)^{\ap/2} u(\bx) = f(\bx),  &\  & \mbox{for} \  \ \bx  \in {\Og},\qquad\\
\label{fBC1}
u(\bx) = 0, &&  \mbox{for} \  \  \bx\in\p{\Og}.
\eea
It can be viewed as an analogy to the fractional Poisson equation in (\ref{fPoisson0})--(\ref{fBC0}).
However, inheriting from the definition of the spectral fractional Laplacian,  the boundary condition in (\ref{fBC1}) is defined locally on $\partial\Og$. 
From the definition in (\ref{fL}),  one can formally obtain the solution of (\ref{fPoisson1})--(\ref{fBC1})  as: 
\bea \label{fPoisson_Solution}
u(\bx)  = \sum _{k\in \mathbb{N}} \mu_k ^{-\ap/2} \widehat{f}_k\,\phi _k ({\bf x}), \qquad \mbox{for} \  \ \bx  \in {\Og},
\eea
where the coefficient $\widehat{f}_k$ is computed by
\begin{eqnarray*}
\widehat{f}_k = \int _{\Omega} f(\bx)\phi_k(\bx)\,d\bx, \qquad k\in{\mathbb N}.
\end{eqnarray*}

\subsection{Regional fractional Laplacian}
\label{section2-3}

On a bounded domain $\Omega\subset{\mathbb R}^n$,  the {\it regional fractional Laplacian}  
(also known the {\it censored fractional Laplacian}) is defined as \cite{Bogdan2003, Guan2005, Guan2006, Guan2006-1}:
\begin{eqnarray}\label{Lap-L}
(-\Dt)_\Og^{\alpha/2} u(\bx) = c_{n,\alpha}\, {\rm P. V.}\int _\Og\frac{u(\bx)-u({\bf y})}{|{\bf y}-\bx|^{n+\alpha}}\,d{\bf y}, \qquad  \mbox{for} \ \ \ap\in(0, 2),
\end{eqnarray}
with the constant $c_{n,\ap}$ defined in (\ref{DefC1ap}). 
In contrast to the fractional Laplacian, the regional fractional Laplacian $(-\Dt)_\Og^{\alpha/2}$ represents the infinitesimal generator of a censored $\ap$-stable process that is obtained from a symmetric $\ap$-stable L\'evy  process by restricting its measure to $\Og$.  
If the domain $\Og = {\mathbb R}^n$, the regional fractional Laplacian collapses to the fractional Laplacian $(-\Dt)^{\ap/2}$. 
To distinguish it from the fractional Laplacian $(-\Dt)^{\ap/2}$, we include the subscript `$\Og$' in the operator $(-\Dt)_\Og^{\ap/2}$ to indicate the restriction of the $\ap$-stable L\'evy  process to the domain $\Og$. 

The regional fractional Laplacian is different from the Dirichlet fractional Laplacian, although they are freely interchanged in some literature.  
In fact, a symmetric $\alpha$-stable L\'evy process killed upon leaving the domain $\Og$ (represented by the Dirichlet fractional Laplacian) is a subprocess of the censored  $\alpha$-stable process  (represented by the regional fractional Laplacian) killing inside the domain $\Og$, i.e., the trajectories may be killed inside $\Og$ through Feynman--Kac transform \cite{Song2008}.  
Moreover, we will  illustrate their difference using a simple example. 
Consider a one-dimensional interval $\Og = (-l, l)$. 
Let $u$ be a smooth function satisfying $u(x) = 0$ for $x\in\Og^c$. 
Then the difference between the regional fractional Laplacian and the Dirichlet fractional Laplacian  can be computed as:
\begin{eqnarray}
\label{Error_L}
Q_1 u(x) &=& \Big((-\Dt)^{\alpha/2}-(-\Delta)_\Og^{\alpha/2}\Big)u(x)\nn\\
&=& c_{1,\alpha}\left(
 \int _{\mathbb{R}} \frac{u(x)-u(y)}{|x-y|^{1+\alpha}}\,dy 
- \int _{-l}^l \frac{u(x)-u(y)}{|x-y|^{1+\alpha}}\,dy \right)\nn\\
&=& c_{1,\alpha} \left(\int _{-\infty}^{-l} \frac{1}{|x-y|^{1+\alpha}}\,dy
+ \int _{l}^{\infty} \frac{1}{|x-y|^{1+\alpha}}\,dy \right) u(x)\qquad\quad\nn\\
&=& \fl{c_{1,\ap}}{\ap} 
\left(\frac{1}{(l+x)^\alpha}+\frac{1}{(l-x)^\alpha}\right)u(x), \qquad\ \mbox{for} \ \ x\in\Og. 
\end{eqnarray}
We find that in the limiting case of $\ap \to 2$,  the difference between the regional fractional Laplacian and the Dirichlet fractional Laplacian vanishes, i.e., $Q_1 \to 0$, due to the constant $c_{1, \ap}\to 0$. 
In other words, the regional fractional Laplacian can be used to approximate the Dirichlet fractional Laplacian as $\ap \to 2$. 
While in the limit of $\ap \to 0$,  the difference in (\ref{Error_L}) reduces to $Q_1u \to  u$, i.e., the Dirichlet fractional Laplacian can be  written as the summation of the regional fractional Laplacian and an identity operator. 
Otherwise, if $\ap\gg 0$ and $\ap \ll2$,  the difference $Q_1u \sim {\mathcal O}\big(1/(l - |x|)^\ap\big)$, which does not tend to zero for any fixed $l$, and as $x \rightarrow \pm l$,  there is $|Q_1u| \to\infty$.

In contrast to the fractional Laplacian, the current understanding of the regional fractional Laplacian still remains very limited.  
Recently, the interior regularity of the regional fractional Laplacian is discussed in \cite{Guan2006, Mou2015}. 
{ It shows that $(-\Delta)^{\ap/2}_{\Og}u \in C^{p}(\Og)$ (for $p \in{\mathbb N}$), if $u \in C^{p,\, s}(\Og)$ for $s\in(\ap, 1]$ or $u \in C^{p+1,\,s}(\Og)$ for $s \in (\ap - 1,  \min(\ap, 1)]$.}
So far, no results on the eigenvalues or eigenfunctions of the regional fractional Laplacian can be found in the literature. 
Here, we expect that our numerical results in Section \ref{section3-2} could  provide insights into the understanding of the spectrum of the regional fractional Laplacian in the future.

The counterpart of the fractional Poisson equation (\ref{fPoisson0})--(\ref{fBC0}), but with the regional fractional Laplacian reads: 
\bea\label{fPoisson2}
(-\Dt)_\Og^{\ap/2} u(\bx) = f(\bx),  &\  & \mbox{for} \  \ \bx  \in {\Og},\qquad\\
\label{fBC2}
u(\bx) = 0, &&  \mbox{for} \  \  \bx\in\p{\Og}.
\eea
The above nonlocal problem has been used to approximate the fractional Poisson problem, in order to reduce computational complexity in solving (\ref{fPoisson0})--(\ref{fBC0}).  
However, as shown in (\ref{Error_L}),  the regional fractional model  does not necessarily provide a consistent approximation to the fractional Poisson equation if $\ap \ll 2$; see more numerical comparison and discussion in Section \ref{section3-3}.

\subsection{Peridynamic operator}
\label{section2-4}

The peridynamic models were originally proposed as a reformulation of the classical solid mechanics in \cite{Silling2000}. 
In contrast to the classical models, it  properly accounts for the near-field nonlocal interactions so as to effectively model elasticity problems with discontinuity and other singularities. 
The general form of this nonlocal operator has the following form: 
\bea\label{peri}
{\mathcal L}u(\bx) = \int_{B(\bx, \delta)} K(\bx, {\bf y})\big(u(\bx) -u({\bf y})\big) \,d{\bf y},
\eea
where $B(\bx, \delta)$,  denoting a ball with its center at point $\bx$ and radius $\dt$, represents the interaction region of point $\bx$. 
The kernel function $K({\bx}, {\bf y}) = K(|\bx - {\bf y}|)$ describes the interaction strength between points ${\bx}$ and ${\bf y}$. 
The  constant $\delta > 0$ denotes the size of material horizon, and it is often chosen to be a small number in practical applications. 
 
Recently, the operator (\ref{peri}) with specially chosen kernel function is used to approximate the fractional Laplacian \cite{Delia2013, Guan2017}.  
We will refer it as the {\it peridynamic operator} and denote it as
\begin{eqnarray}\label{Lap-R}
(-\Dt)_\delta^{\alpha/2} u(\bx) = c_{n,\alpha} \int_{B(\bx, \delta)} \frac{u({\bf x})-u({\bf y})}{|{\bf x}-{\bf y}|^{n+\alpha}}\,d{\bf y},
\end{eqnarray}
i.e.,  the kernel function in this case is taken as:
\beas
K_\dt(\bx, {\bf y}) = \left\{\begin{array}{ll} \displaystyle \fl{c_{n, \ap}}{|\bx - {\bf y}|^{n + \ap}},  \qquad & \mbox{if} \ \ {\bf y} \in B(\bx, \dt),\\
0, &\mbox{otherwise}.
\end{array}\right.
\eeas
In other words,  $K_\delta(\bx, {\bf y})$ in the peridynamic operator represents a hard-threshold of the kernel function $K(\bx, {\bf y}) = {c_{n, \ap}}/{|\bx - {\bf y}|^{n + \ap}}$ of the fractional Laplacian, which can be viewed as a truncation of  $K(\bx, {\bf y})$ in the fractional Laplacian. 
In the limiting case of $\dt\to\infty$, the peridynamic operator (\ref{Lap-R}) coincides with the fractional Laplacian (\ref{riesz}), and thus it is often used to approximate the fractional Laplacian by choosing a sufficiently large $\dt$  \cite{Delia2013, Guan2017}. 
On the other hand, note that the kernel function $K(\bx, {\bf y})$ has an algebraic decay of order $n+\ap$, which presents a heavy tail that accounts for considerable far field interactions. 
Hence, the cutoff of  the kernel function $K(\bx, {\bf y})$ outside of the horizon $B(\bx, \dt)$ may have a significant impact on its approximation  to the fractional Laplacian  as we shall show next. 

Similarly, we choose a smooth function $u$ satisfying $u(x) = 0$ for $x\in\Og^c$ with $\Og = (-l, l)$ to illustrate the difference between the peridynamic operator and the fractional Laplacian.   
Here, we assume that the horizon size $\delta$ in (\ref{Lap-R}) is large enough, such that $\delta > \max\{l-x, l+x\}$ for any point  $x \in (-l, l)$. 
Then, we can compute their difference as:
\begin{eqnarray}
\label{Error_R}
Q_2 u(x) &=& \Big((-\Delta)^{\alpha/2} -(-\Dt)_\delta^{\alpha/2}\Big)u(x)\nn\\
&=& c_{1,\alpha}\bigg(\int _{\mathbb{R}} \frac{u(x)-u(y)}{|x-y|^{1+\alpha}}\,dy
       - \int _{x-\delta}^{x+\delta} \frac{u(x)-u(y)}{|x-y|^{1+\alpha}}\,dy \bigg)\nn\\
&=& c_{1,\alpha}\left(\int _{-\infty}^{x-\delta} \frac{1}{|x-y|^{1+\alpha}}\,dy
       + \int _{x+\delta}^{\infty} \frac{1}{|x-y|^{1+\alpha}}\,dy \right)u(x)\qquad\qquad \nn\\
&=& \fl{c_{1,\ap}}{\ap}\fl{2}{\delta^\ap}\, u(x) , \qquad  \mbox{for} \ \ x\in\Og. 
\end{eqnarray}
It shows that the difference of these two operators is of order ${\mathcal O}(1/\delta^\ap)$ when $u(x)$ is uniformly bounded on $\Og$, hence their difference vanishes as $\dt \to \infty$. 
On the other hand,  the convergence of the peridynamic operator to the fractional Laplacian as $\dt \to \infty$ depends on the power $\ap$, and it may degenerate rapidly for small $\ap$. 
Additionally, in the limiting case of $\ap \to 2$, the difference $Q_2 u \to 0$,  because the coefficient $c_{n, \ap} \to 0$. 

A peridynamic model for describing the steady-state displacement of a finite microelastic bar can be formulated as follows \cite{GunLeh, Delia2013, MenDu}:
\begin{eqnarray}\
\label{fPoisson2}
\displaystyle(-\Dt)_\delta^{\alpha/2} u(\bx) = f(\bx), &\ &\bx \in\Og,\qquad\\
\label{fBC2}
\displaystyle u(\bx) = 0,    && \bx \in\Og_\dt,
\end{eqnarray}
where $\Og_\dt$  defines the boundary region with width of $\dt$. 
Note that the peridynamic model (\ref{fPoisson2})--(\ref{fBC2}) is enclosed with a nonlocal Dirichlet boundary condition on a nonconventional finite ``volume" boundary region of size $\delta$, i.e. $\Og_\dt$.  
Hence, the boundary condition in (\ref{fBC2}) is also referred to as a {\it finite volume constraint} in the literature \cite{Du2012, Delia2013}. 
The peridynamic operator $(-\Dt)_\dt^{\ap/2}$ is often used to approximate the fractional Laplacian $(-\Dt)^{\ap/2}$, while the  nonlocal problem  (\ref{fPoisson2})--(\ref{fBC2}) is often used  to approximate the fractional Poisson equation in (\ref{fPoisson0})--(\ref{fBC0}). 
The well-posedness of (\ref{fPoisson2})--(\ref{fBC2}) can be found in the literature \cite{Du2012}.

The peridynamic operator in (\ref{Lap-R}) can be viewed as an infinitesimal generator of a symmetric $\ap$-stable L{\'e}vy process by restricting its measure to $B(\bx, \dt)$.
In contrast to the regional fractional Laplacian operator, the interaction region of point ${\bf x}$ in  the peridynamic operator is symmetric with respect to itself. 
Hence, the peridynamic operator is expected to provide a symmetric approximation for a homogeneous elastic material.  \\
 
In summary, the fractional Laplacian (\ref{riesz}), spectral fractional Laplacian (\ref{fL}), regional fractional Laplacian (\ref{Lap-L}), and the peridynamic operator (\ref{Lap-R}) are all nonlocal operators in which every point $\bx$ interacts with other points ${\bf y}$ over certain long distance. 
For a point $\bx \in \Og$,  the fractional Laplacian $(-\Dt)^{\ap/2}$ accounts for the interactions between $\bx$ and ${\bf y}$ for all ${\bf y}  \in{\mathbb R}^n\backslash\{\bx\}$. 
By contrast, the interaction region of $\bx$ in the regional fractional Laplacian $(-\Dt)_\Og^{\ap/2}$ is truncated to $\Og\backslash\{\bx\}$, i.e., the same domain of $\bx$, while  the interaction region of the peridynamic operator $(-\Dt)_\delta^{\ap/2}$ reduces to $B(\bx, \dt)\backslash\{\bx\}$. 
We will further compare them in Section \ref{section3}.

\section{Numerical comparisons}
\setcounter{equation}{0}
\label{section3}

In this section, we further compare these four nonlocal operators by studying their nonlocal effects, eigenvalues and eigenfunctions, and the solution behavior of the corresponding nonlocal problems. 
In our simulations, the spectral fractional Laplacian is discretized by using the finite difference method combined with matrix transfer techniques introduced in \cite{DuoJuZhang0016},  while the other three operators are discretized by the finite difference method based on weighted trapezoidal rules proposed in \cite{DuoWykZhang}. 
Our numerical results provide insights not only to further understand these operators but also to improve the analytical results in the literature. 

In the following, we  will consider the one-dimensional cases. 
For notational simplicity, we will also  use ${\mathcal L}_h$ to represent the fractional Laplacian, ${\mathcal L}_s$ for the spectral fractional Laplacian, ${\mathcal L}_r$ for the regional fractional Laplacian, and ${\mathcal L}_p$ for the peridynamic operator.

\subsection{Nonlocal effects of operators}
\label{section3-1}

We compare the nonlocal effects of these four operators by acting them on functions with compact support on the domain  $\Og = (-1, 1)$.
\vskip 10pt
\noindent {\bf Example 1.\ }  Consider the  function 
\bea\label{usin}
u(x) = \left\{\begin{array}{ll}
\displaystyle \sin\Big(\fl{\pi(1+x)}{2}\Big),  \quad  \ \ &\text{if }\ \, x\in \Og,\\
0, &\text{otherwise}, 
\end{array}\right.\qquad x\in{\mathbb R},
\eea
which is continuous on the whole space ${\mathbb R}$. 
It is easy to obtain that
\beas \label{usin_exact}
(-\Dt_\Og)^{\ap/2} u(x) = \left(\fl{\pi}{2}\right)^\ap \sin\Big(\fl{\pi(1+x)}{2}\Big), \qquad \mbox{for\ \  $x\in{\Og}$},
\eeas
that is, the function from the spectral fractional Laplacian can be found exactly. 
While we will numerically compute the functions from the other three operators using the finite difference method proposed in \cite{DuoWykZhang}. 

\begin{figure}[htb!]
\centerline{
\includegraphics[height=5.160cm,width=6.960cm]{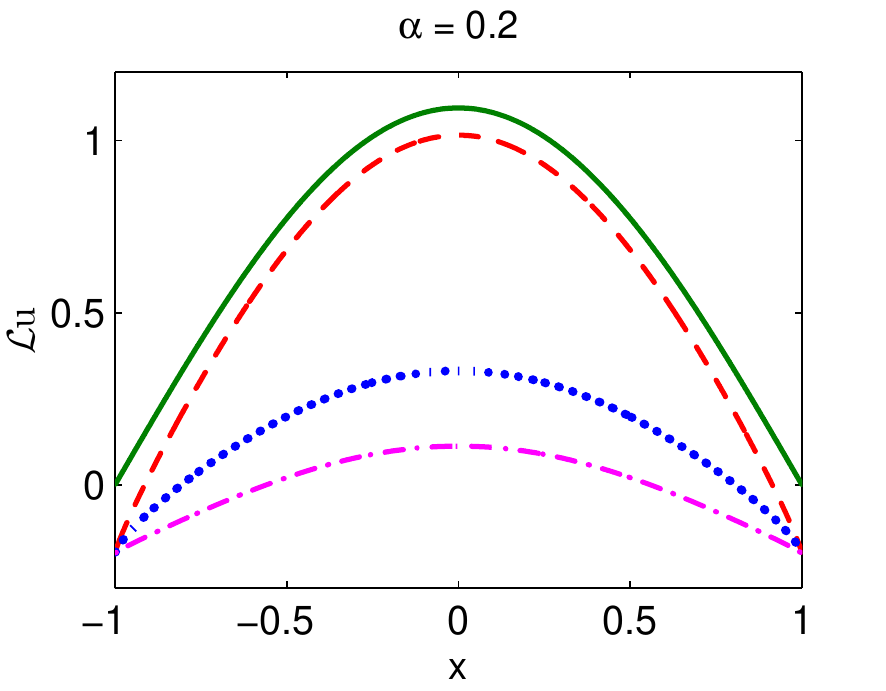}\quad
\includegraphics[height=5.160cm,width=6.960cm]{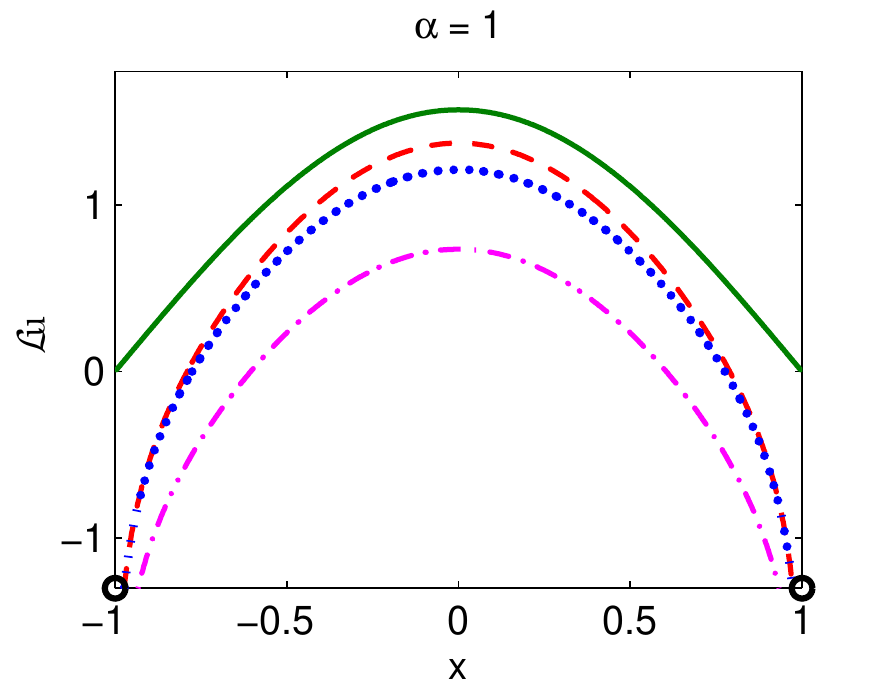}}
\centerline{
\includegraphics[height=5.160cm,width=6.960cm]{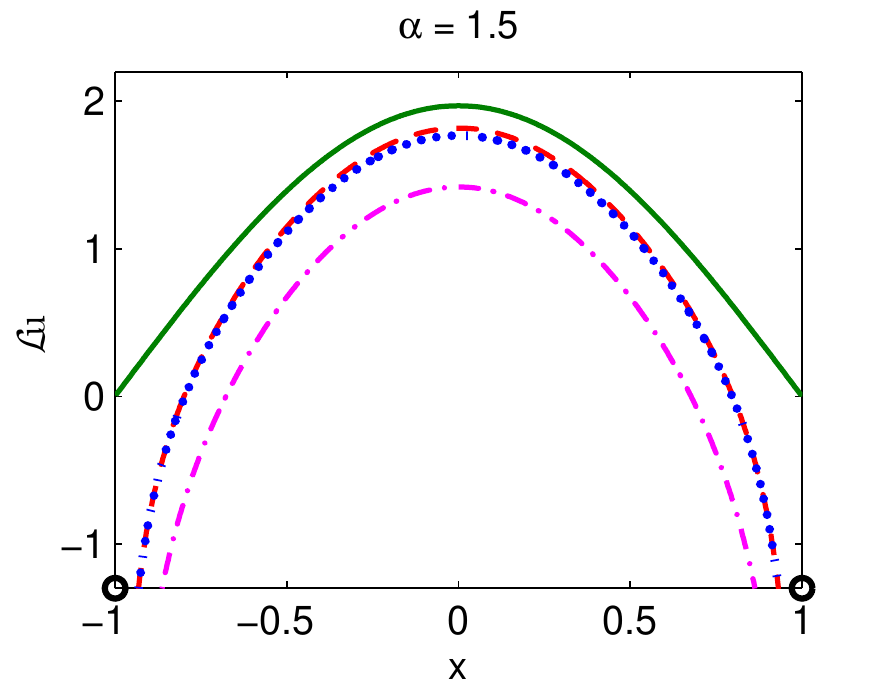}\quad
\includegraphics[height=5.160cm,width=6.960cm]{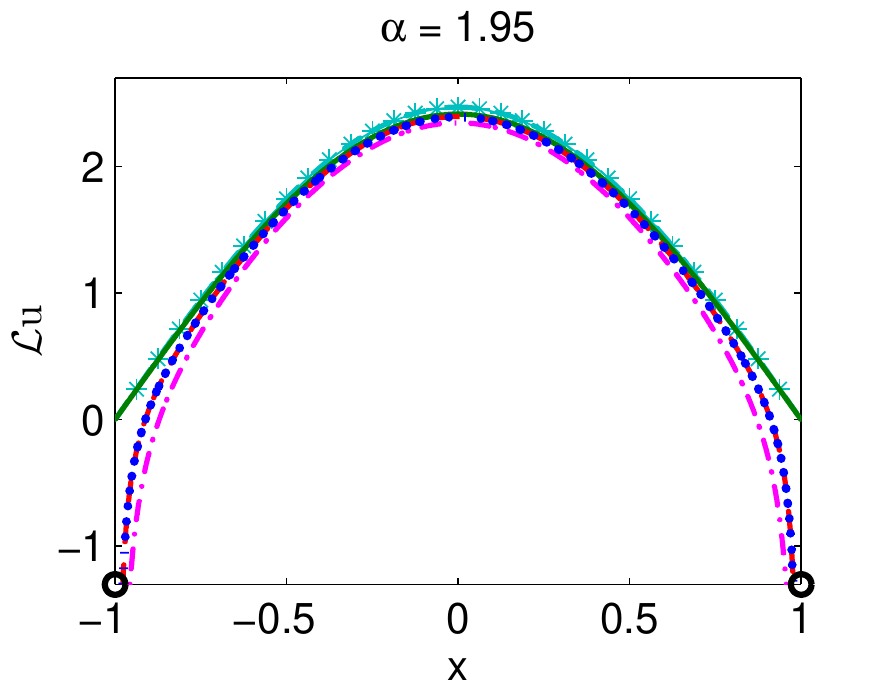}}
\caption{Comparison of the function ${\mathcal L}u$ with $u$ defined in (\ref{usin}), where the operator ${\mathcal L}$ represents ${\mathcal L}_s$ (solid line),  ${\mathcal L}_h$ (dashed line), ${\mathcal L}_r$ (dash-dot line), or ${\mathcal L}_p$ with $\dt = 4$ (dotted line). For easy comparison, the result for ${\mathcal L} = -\p_{xx}$ (line with symbols ``*") is included in the plot of $\ap = 1.95$.  
For $\ap = 1$, $1.5$ or $1.95$, the plots in $y$-direction are partially presented.}
\label{Fig1}
\end{figure}

In Figure \ref{Fig1}, we compare the functions ${\mathcal L}_iu$ for $i = s,\, h,\, r$, or $p$. 
The results clearly suggest the difference between these four operators, especially the function ${\mathcal L}_su$ from the spectral fractional Laplacian is significantly different from those of the other three operators. 
It shows that for any $\ap \in (0, 2)$,  the function ${\mathcal L}_su$ is proportional to the function $u$ on $(-1, 1)$. 
In contrast, the properties of ${\mathcal L}_iu$ (for $i = h, r$,  or $p$) significantly depend on the parameter $\ap$.
For  $\ap \in (0, 1)$, the functions ${\mathcal L}_iu$ (for $i = h, r$,  or $p$) exist on the closed domain $\overline{\Og}$, but they are very different between operators. The smaller the parameter $\ap$, the larger the differences.  
For $\ap \in [1, 2)$, the functions ${\mathcal L}_iu$ do not exist at the boundary points, i.e., $x= \pm 1$.
{As $\ap \to 2$,  the functions ${\mathcal L}_iu$ (for $i = h, r$,  or $p$)  converge to $-u_{xx}$ for $x\in (-1,1)$.} 

Additionally, Figure \ref{Fig1} shows that both the regional fractional Laplacian and peridynamic operator can be used to approximate the fractional Laplacian, if $\ap$ is close to 2 (see Fig. \ref{Fig1} for $\ap = 1.95$). 
\begin{figure}[htb!]
\centerline{
\includegraphics[height=5.560cm,width=7.860cm]{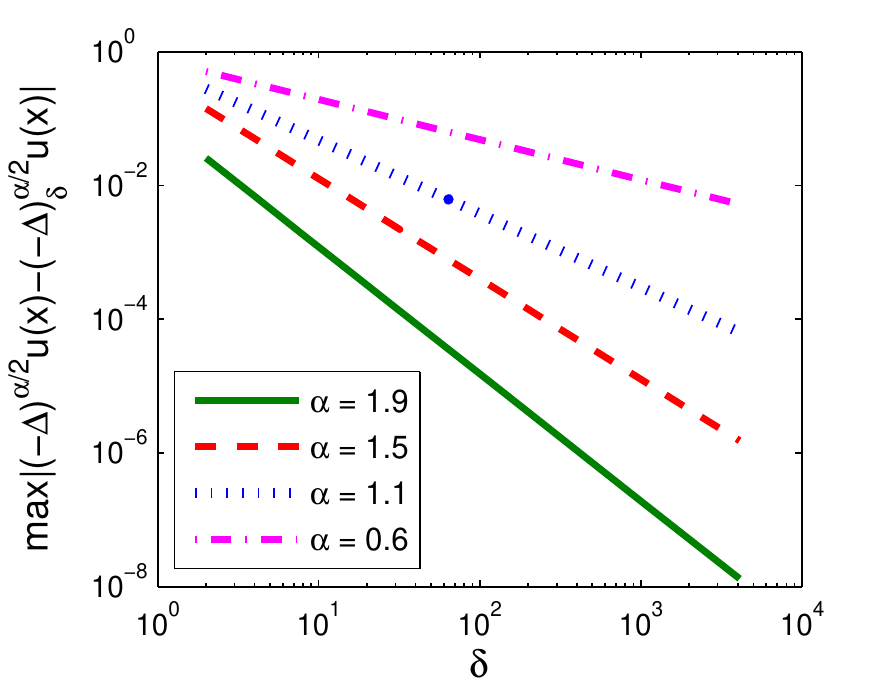}
} 
\caption{Difference between the peridynamic operator and the fractional Laplacian versus the parameter $\ap$, where $u(x)$  is defined in (\ref{usin}).}
\label{Fig1-1}
\end{figure}
For small $\ap$, the results from the regional fractional Laplacian is inconsistent with that from the fractional Laplacian. 
However, the peridynamic operator can still provide a good approximation to the fractional Laplacian by enlarging the horizon size $\dt$. 
Figure \ref{Fig1-1} presents the differences between the functions ${\mathcal L}_pu$ and ${\mathcal L}_hu$ for various $\ap$ and $\dt$. 
It shows that for a fixed horizon size $\dt$, the difference between these two operators dramatically decreases as $\ap$ increases. 
On the other hand, Figure \ref{Fig1-1} implies that for small $\ap$ the convergence of  the function ${\mathcal L}_pu$ to 
${\mathcal L}_hu$ could be very slow. 
For instance,  for $\ap = 0.6$,  the difference in Fig. \ref{Fig1-1} is around 0.005 for a horizon size $\dt = 4000$. 
In fact, the nonlocal interactions decay slowly for small $\ap$, and thus a large horizon size $\dt$ is needed for the peridynamic operator to better approximate the fractional Laplacian. 

\vskip 15pt
\noindent {\bf Example 2.\ }  Consider the  function 
\begin{eqnarray}\label{ufun}
u(x) = \left\{\begin{array}{ll}
\displaystyle (1-x^2)^{q+\fl{\ap}{2}},  \quad\  &\text{for}\ \  x\in \Og,\\
\displaystyle 0, &\text{otherwise}, 
\end{array}\right. \qquad x\in{\mathbb R},
\end{eqnarray}
for $q \in {\mathbb N}$.
For the fractional Laplacian, the analytical value can be found as: 
\beas
(-\Dt)^{\ap/2}u(x) = \fl{2^\ap\Gamma(\fl{\ap+1}{2})\Gamma(\fl{\ap}{2}+q+1)}{\sqrt{\pi}\Gamma(q+1)} {\ }_2F_1\bigg(\fl{\ap +1}{2},\,-q;\,\fl{1}{2};\,x^2 \bigg),  \qquad \text{for}\ \ \ x\in \Omega,
\eeas
where $_2F_1$ denotes the Gauss hypergeometric function.   
Moreover, we can obtain the exact values of $(-\Dt)_\Og^{\ap/2}u$ and $(-\Dt)_\dt^{\ap/2}u$ by using their relation to the fractional Laplacian in (\ref{Error_L}) and (\ref{Error_R}), respectively. 
For the spectral fractional Laplacian,  we numerically  computed $(-\Dt_\Og)^{\ap/2}u$ by the finite difference method proposed in \cite{DuoJuZhang0016}. 

\begin{figure}[htb!]
\centerline{
\includegraphics[height=5.160cm,width=6.960cm]{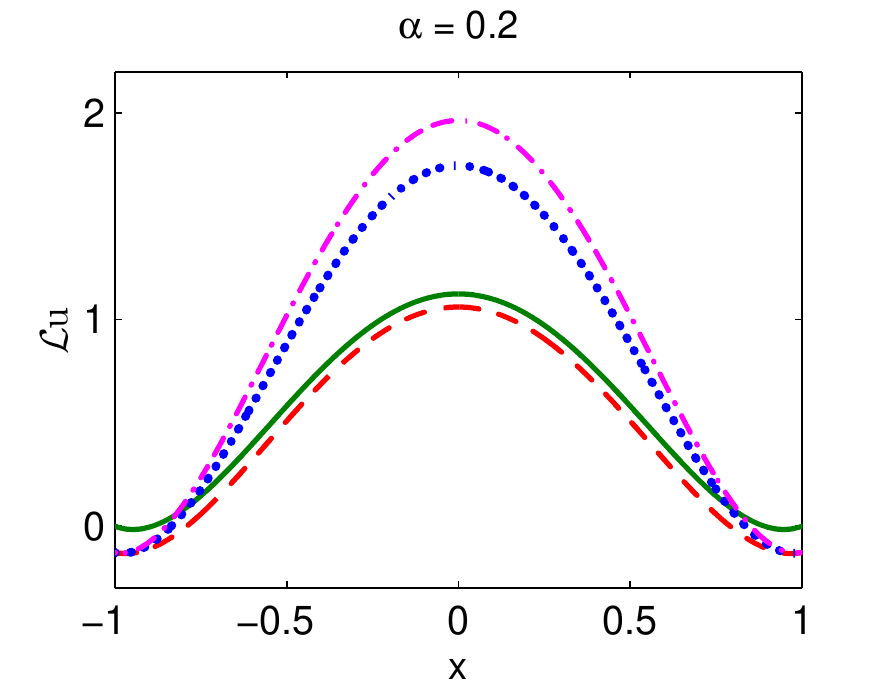}\quad
\includegraphics[height=5.160cm,width=6.960cm]{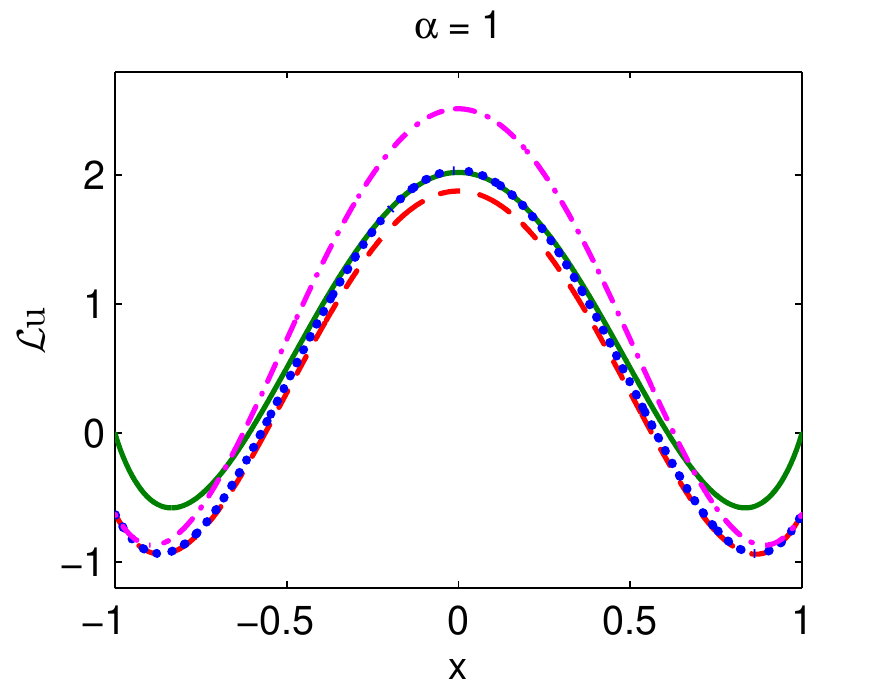}}
\centerline{
\includegraphics[height=5.160cm,width=6.960cm]{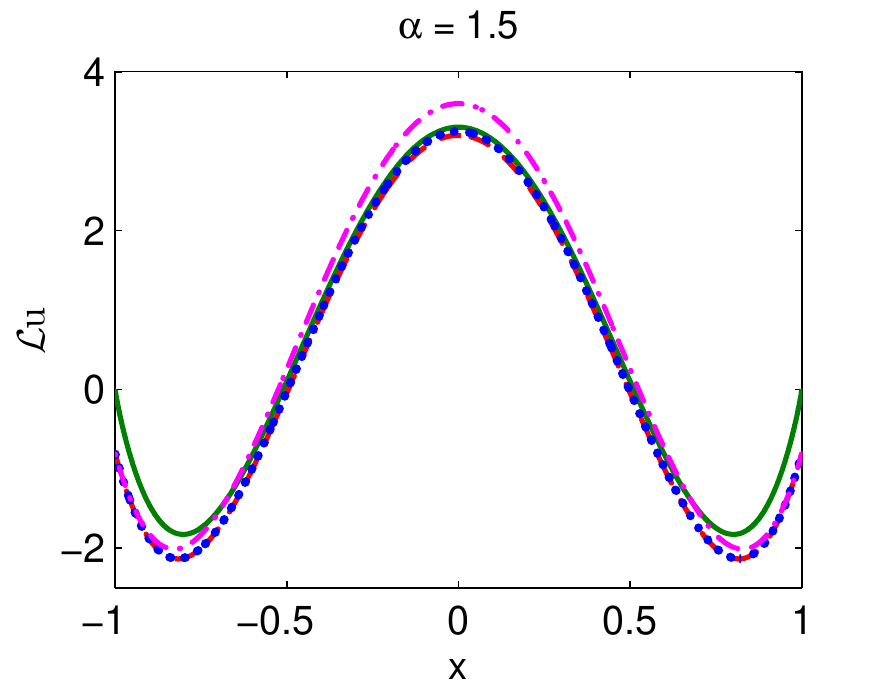}\quad
\includegraphics[height=5.160cm,width=6.960cm]{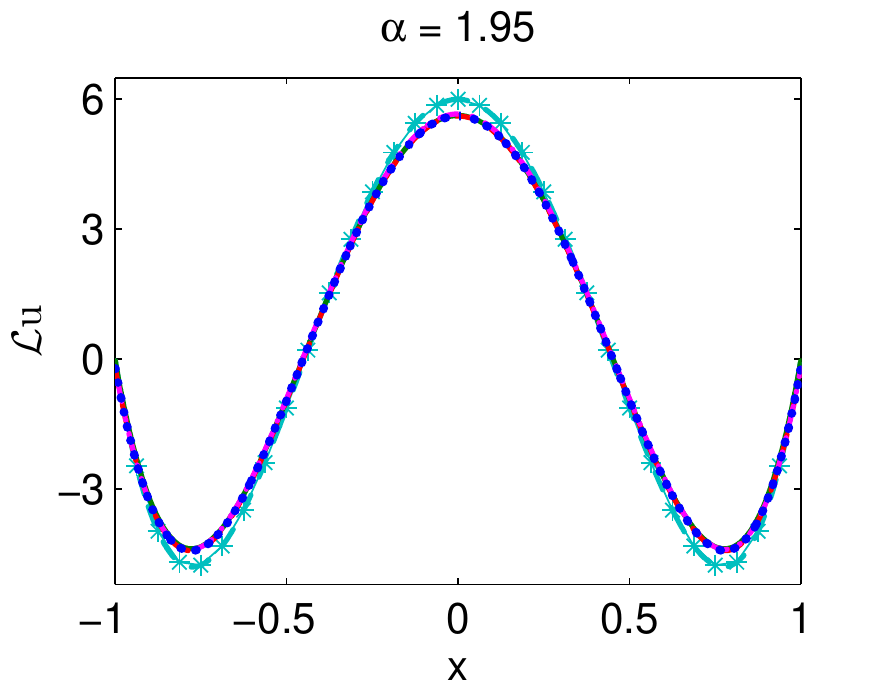}}
\caption{Comparison of the function ${\mathcal L}u$ with $u$ defined in (\ref{ufun}) with $q = 2$.  The operator ${\mathcal L}$ represents ${\mathcal L}_s$ (solid line),  ${\mathcal L}_h$ (dashed line), ${\mathcal L}_r$ (dash-dot line),  or ${\mathcal L}_p$ with $\dt = 4$ (dotted line). For easy comparison, the result for ${\mathcal L} = -\p_{xx}$ (line with symbols ``*") is included in the plot of $\ap = 1.95$.  }
\label{Fig3}
\end{figure}

Figure \ref{Fig3} displays the functions ${\mathcal L}_iu$ (for $i = s, h, r$, or $p$) for various $\ap$, where $u$ is defined in  (\ref{ufun}) with $q = 2$.  
It shows that the functions ${\mathcal L}_iu$ exist on the closed domain $[-1, 1]$ for any $\ap\in(0, 2)$, but their values are very different, especially for small $\ap$. 
For the spectral fractional Laplacian,  the values of ${\mathcal L}_su$ are always zero at boundary points. 
{For the regional fractional Laplacian, the function $u$ in (\ref{ufun}) with $q = 2$ satisfies the conditions that $u \in C^2([-1,1])$ and $u'(\pm 1)=0$, which  guarantee the existence of the function ${\mathcal L}_ru$ for any $\alpha \in (0 ,2)$  \cite{Guan2006}. }
Since the function $u(\pm 1) = 0$ and the relations in (\ref{Error_L}) and (\ref{Error_R}),  the values of ${\mathcal L}_i u$ (for $i = h, r$ and $p$) are the same at boundary points, but they are nonzero. 
Figure \ref{Fig3} also shows that both the regional fractional Laplacian and the peridynamic operator with relatively small $\dt$ could provide a good approximation to the fractional Laplacian, if $\ap$ is large (see Fig. \ref{Fig3} for $\ap = 1.95$).  
While $\ap$ is small, although the peridynamic operator can be still used to approximate the fractional Laplacian with a large $\dt$,  the regional fractional Laplacian is inconsistent with the fractional Laplacian.
 
\begin{figure}[htb!]
\centerline{
\includegraphics[height=5.160cm,width=6.960cm]{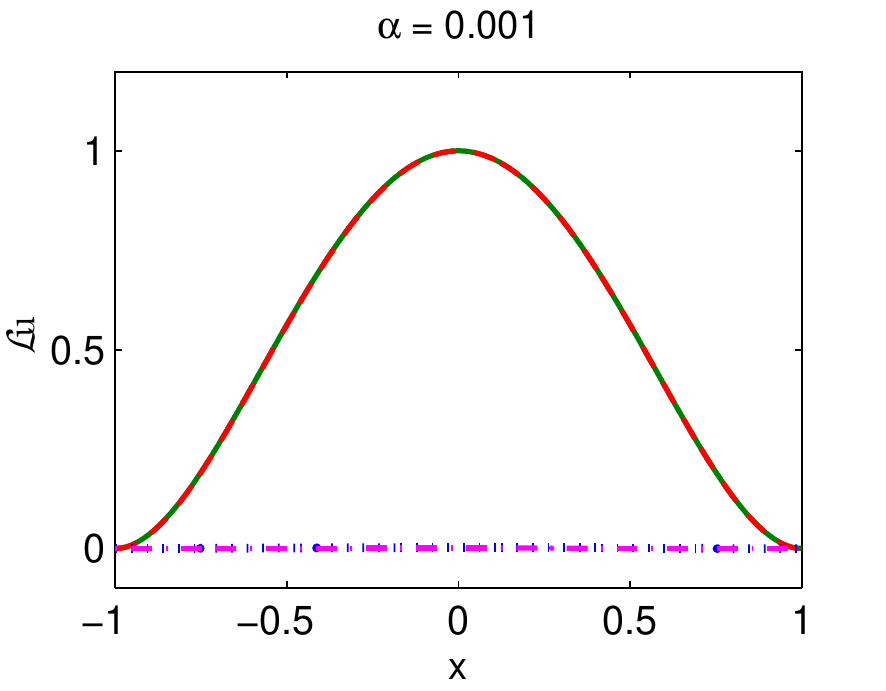}} 
\caption{Comparison of the function ${\mathcal L}u$ with $u$ defined in (\ref{ufun}) with $q = 2$.  The operator ${\mathcal L}$ represents ${\mathcal L}_s$ (solid line),  ${\mathcal L}_h$ (dashed line), ${\mathcal L}_r$ (dash-dot line), or ${\mathcal L}_p$ with $\dt = 4$ (dotted line).}\label{Fig3-1}
\end{figure}

\begin{figure}[htb!]
\centerline{
\includegraphics[height=5.160cm,width=6.960cm]{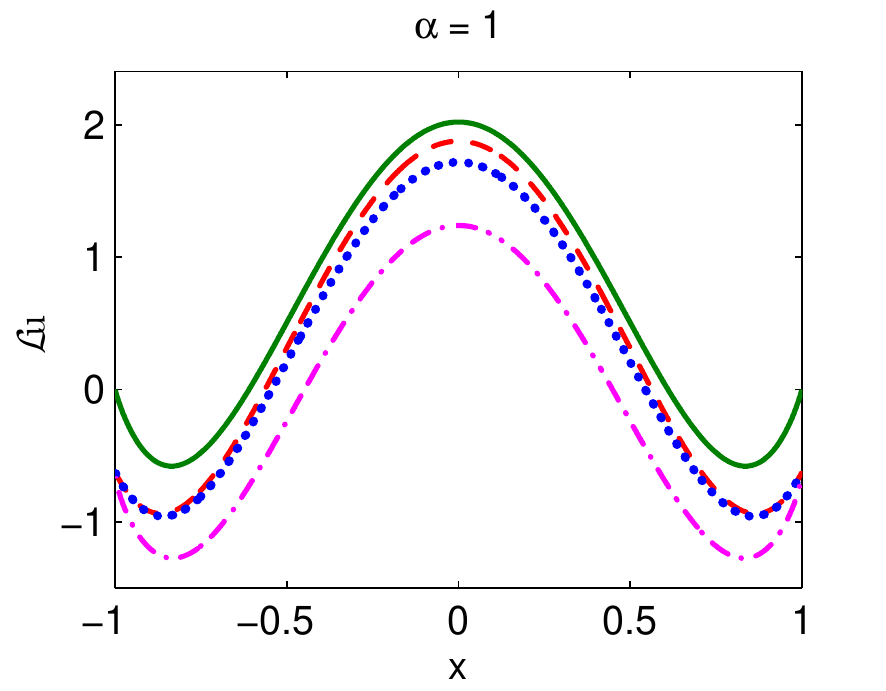}\quad
\includegraphics[height=5.160cm,width=6.960cm]{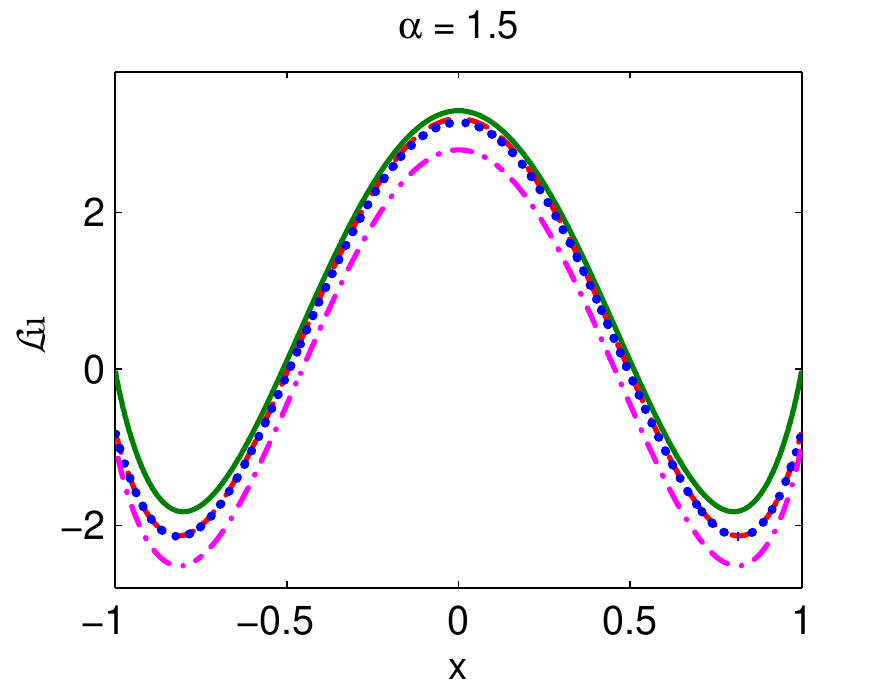}}
\caption{Comparison of the function ${\mathcal L}u$ with $u$ defined in (\ref{ufun}) with $q = 1$.  The operator ${\mathcal L}$ represents ${\mathcal L}_s$ (solid line),  ${\mathcal L}_h$ (dashed line), ${\mathcal L}_r$ (dash-dot line), or ${\mathcal L}_p$ with $\dt = 4$ (dotted line).} \label{Fig66}
\end{figure}

Figure \ref{Fig3} additionally shows that as $\ap \to 2$, the differences between the four operators become insignificant (see Fig. \ref{Fig3} for $\ap = 1.95$), and  the functions ${\mathcal L}_iu$ (for $i = h,\, s,\, r,\,$ or $p$) converge to $-\p_{xx}u$, that is, the four operators converge to the standard Dirichlet Laplace operator $-\Dt$.  
To understand their properties as $\ap  \to 0$, Figure \ref{Fig3-1} shows the functions ${\mathcal L}_iu$ (for $i = s, h, r$, or $p$) for a small value  of $\ap = 0.001$. 
The definition of the spectral fractional Laplacian in (\ref{fL}) implies that the function ${\mathcal L}_su$ converges to $u$, as $\ap \to 0$.
Fig. \ref{Fig3-1}  shows  that  the function ${\mathcal L}_hu$ from the fractional Laplacian  converges to $u$ as $\ap \to 0$, confirming the analytical results in \cite{DiNezza2012}. 
By contrast, the functions ${\mathcal L}_ru$ from the regional fractional Laplacian and ${\mathcal L}_pu$ from the peridynamic operator converge to a zero function. 
This can be easily obtained from their relation to the fractional Laplacian in (\ref{Error_L}) and (\ref{Error_R}), respectively. 
Moreover,  our numerical results seem to suggest that for $\ap\in[1, 2)$, if the function $u\in C^{1, \ap/2}(\overline{\Og})$ and $u'(\pm 1) = 0$, then the value ${\mathcal L}_ru$ from the regional fractional Laplacian exists; see Figure \ref{Fig66}. 
Hence, we conjecture that the regularity results in \cite{Guan2006} might be able to improve to $u\in C^{1, \ap/2}(\overline{\Og})$ at least for one-dimensional case. More analysis needs to be carried out for further understanding of this issue. 

\subsection{Eigenvalues and eigenfunctions}
\label{section3-2}

In this section, we compare the four nonlocal operators by studying their eigenvalues and eigenfunctions  on a one-dimensional bounded domain $\Og = (-l, \, l)$.
 
Denote $\lambda_k^i$ and  $\phi_k^i$ as the $k$-th (for $k\in{\mathbb N}$) eigenvalue and eigenfunction of the nonlocal operator ${\mathcal L}_i$  on $\Og$ with the corresponding homogeneous Dirichlet boundary conditions, where $i = h,\, s,\, r$, or $p$.
It is well known that the eigenvalues and eigenfunctions of the spectral fractional Laplacian ${\mathcal L}_s$ can be found analytically, i.e.,  
\beas\label{eig_s}
\lambda^{s}_{k} = \mu_k^{\ap/2}  = \left(\fl{k\pi}{2l}\right)^{\ap}, \qquad \phi_k^s(x) = \sqrt{\fl{1}{l}}\sin\bigg(\fl{k\pi}{2}\bigg(1+\fl{x}{l}\bigg)\bigg),\quad \ \ x\in(-l, l), 
\eeas
for $k\in{\mathbb N}$. 
For  the other operators, so far no analytical results can be found in the literature, and thus we will compute their eigenvalues and eigenfunctions numerically. 

In Table \ref{Tab1}, we present the eigenvalues of the Dirichlet fractional Laplacian ${\mathcal L}_h$, spectral fractional Laplacian ${\mathcal L}_s$, and regional fractional Laplacian ${\mathcal L}_r$,  on the domain $\Omega = (-1,1)$.  
We leave the peridynamic operator ${\mathcal L}_p$ out of our comparison here,  since its spectrum depends on the horizon size $\dt$.
From Table \ref{Tab1} and our extensive numerical studies, we find that 
\beas
\lambda_k^r <   \lambda_k^h < \lambda_k^s, \qquad \mbox{for} \ \ \ap\in(0, 2) \ \ \mbox{and} \ \  k\in{\mathbb N}, 
\eeas
that is, the eigenvalues of the regional fractional Laplacian are much smaller than those of the Dirichlet fractional Laplacian and the spectral fractional Laplacian. 
However, as $\ap \to 2$ the eigenvalue $\lambda_k^i$ of these three operators converges to  $\mu_k = k^2\pi^2/4$ -- the $k$th eigenvalue of the standard Dirichlet Laplace operator $-\Dt$ on $(-1, 1)$. 

\begin{figure}[htb!]
\centerline{
\includegraphics[height=5.160cm,width=7.560cm]{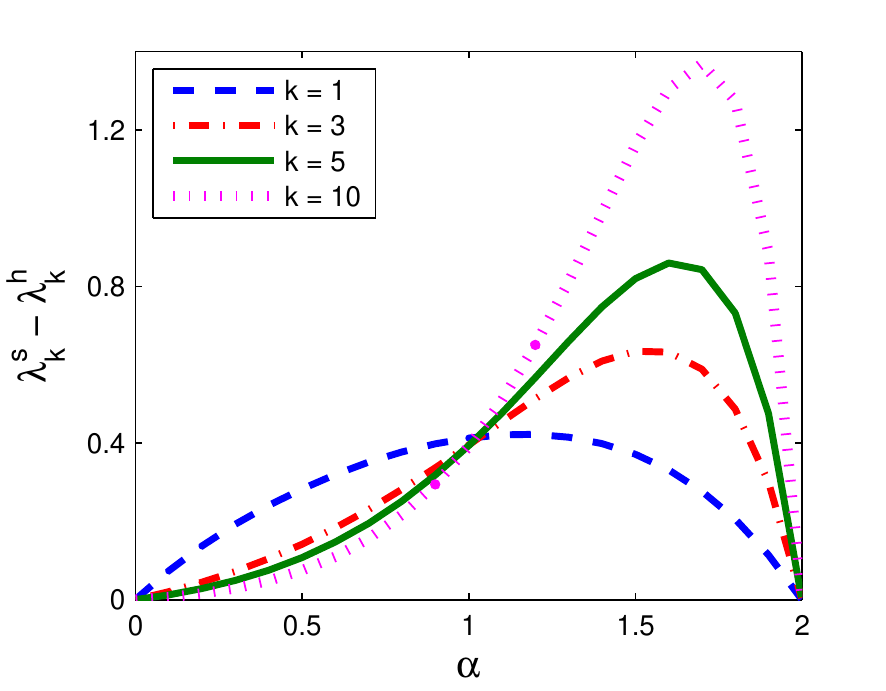}
\includegraphics[height=5.160cm,width=7.560cm]{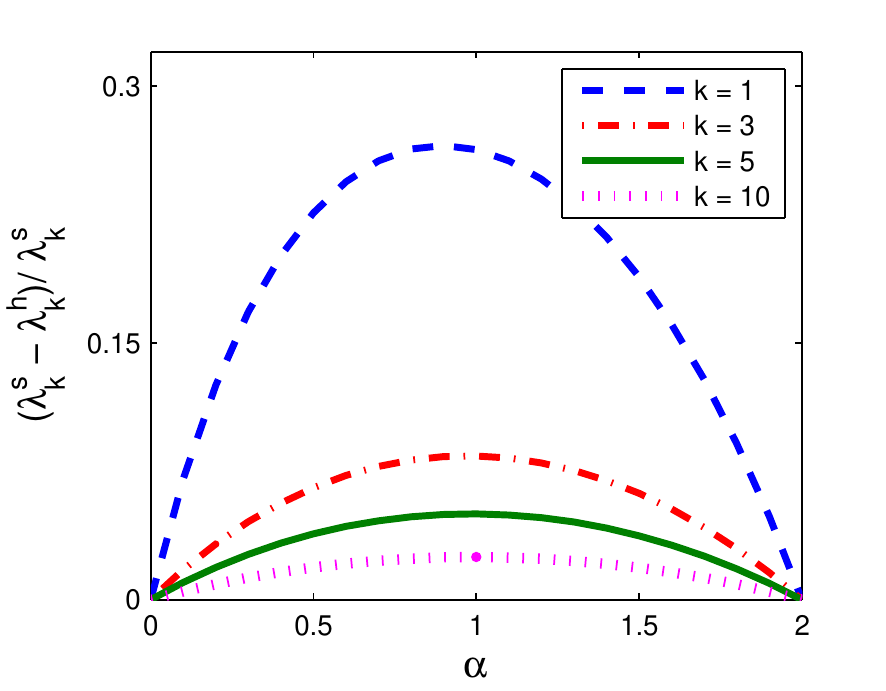}}
\caption{The absolute (left panel) and relative (right panel) differences in the eigenvalues of the fractional Laplacian and the spectral fractional Laplacian.}\label{fig01}
\end{figure}

In \cite{Servadei2014}, it is proved that the first eigenvalue of the Dirichlet fractional Laplacian is strictly smaller than that of the spectral fractional Laplacian, i.e., $\lambda_1^h < \lambda_1^s$,  for $\ap\in(0, 2)$.
Our numerical results  in Table \ref{Tab1}  confirm this conclusion and additionally suggest that the eigenvalue $\lambda_k^h$ is strictly smaller than $\lambda_k^s$,  for any $k \in{\mathbb N}$. 
Furthermore, we present the difference between the eigenvalues $\lambda_k^s$ and $\lambda_k^h$ for various $\ap$ and $k$ in Figure \ref{fig01}. 
It shows that the difference between the eigenvalues $\lambda_k^s$ and $\lambda_k^h$ depends on  both parameters $\ap$ and $k$.
For a given $k\in{\mathbb N}$, there exists a critical value $\ap_{k, {\rm cr}}$ where the gap between  $\lambda_k^s$ and $\lambda_k^h$ is maximized. 
The value of $\ap_{k, {\rm cr}}$ increases as $k\in{\mathbb N}$ increases (see Fig. \ref{fig01} left). 
On the other hand, as $k\to\infty$ the relative difference between the eigenvalues $\lambda_k^s$ and $\lambda_k^h$ decreases quickly (see Fig. \ref{fig01} right).

\begin{figure}[htb!]
\centerline{
\includegraphics[height=4.560cm,width=6.960cm]{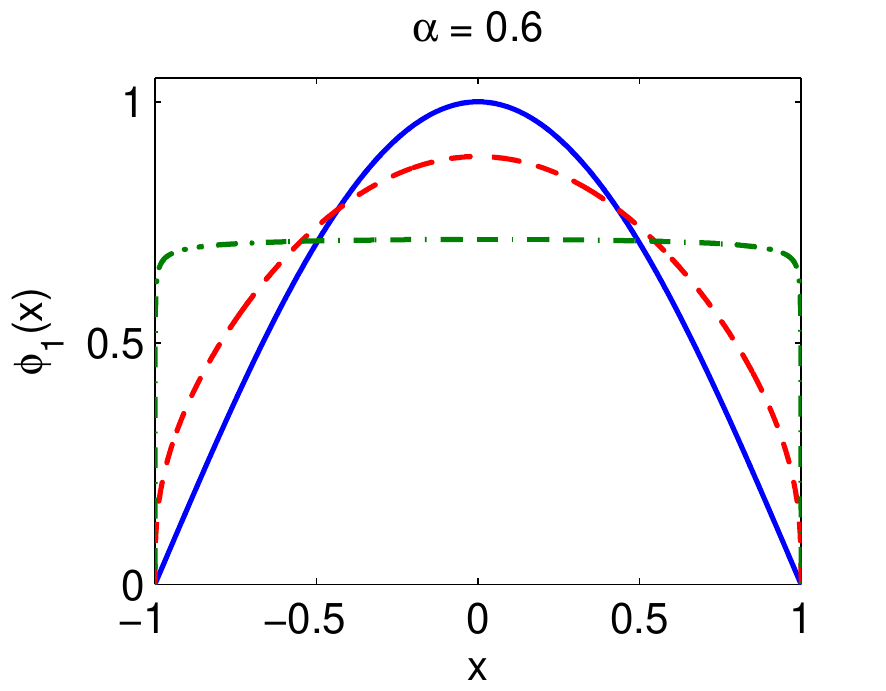}
\includegraphics[height=4.560cm,width=6.960cm]{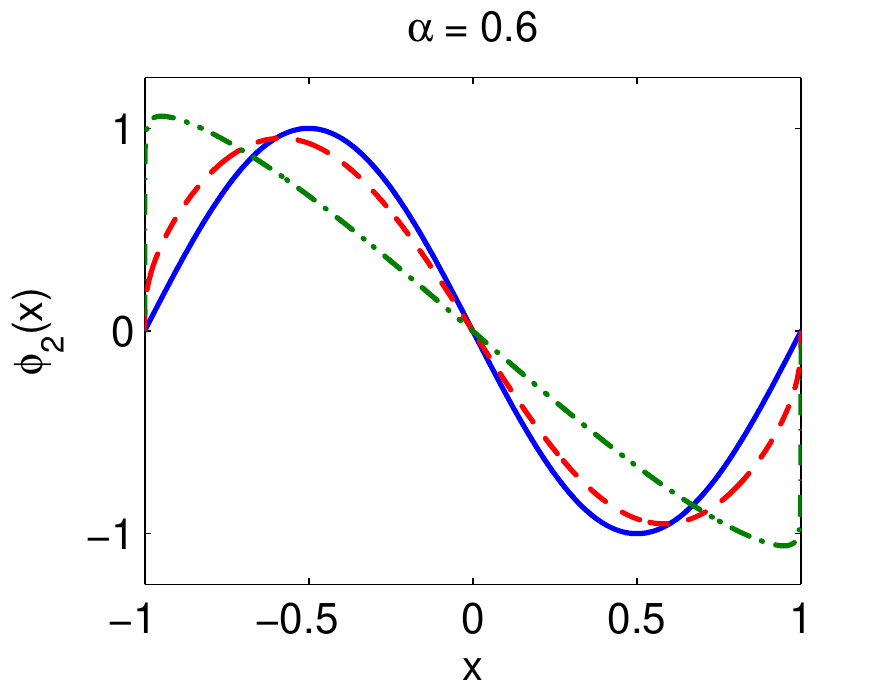}
}
\centerline{
\includegraphics[height=4.560cm,width=6.960cm]{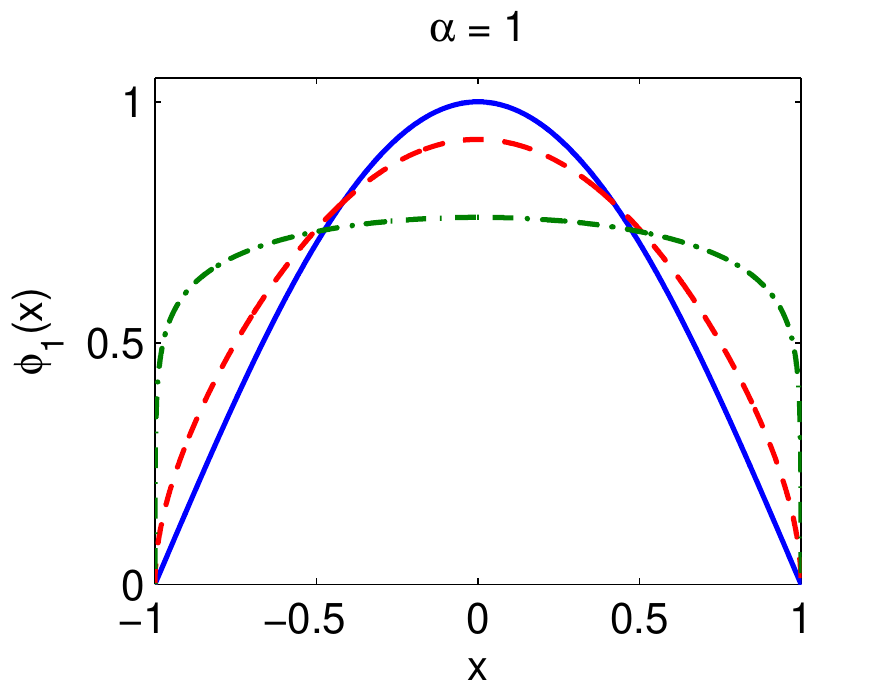}
\includegraphics[height=4.560cm,width=6.960cm]{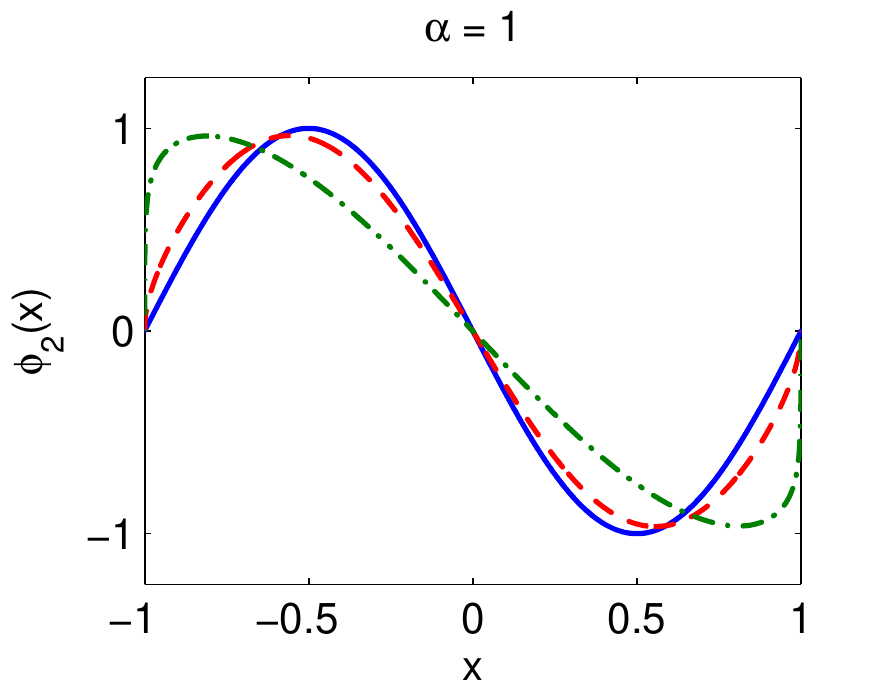}
}
\centerline{
\includegraphics[height=4.560cm,width=6.960cm]{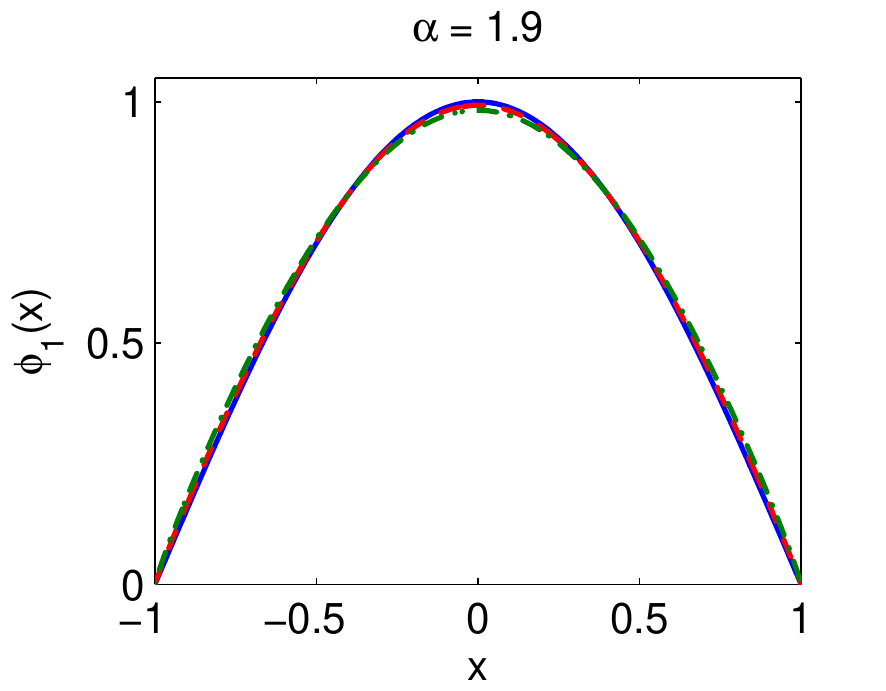}
\includegraphics[height=4.560cm,width=6.960cm]{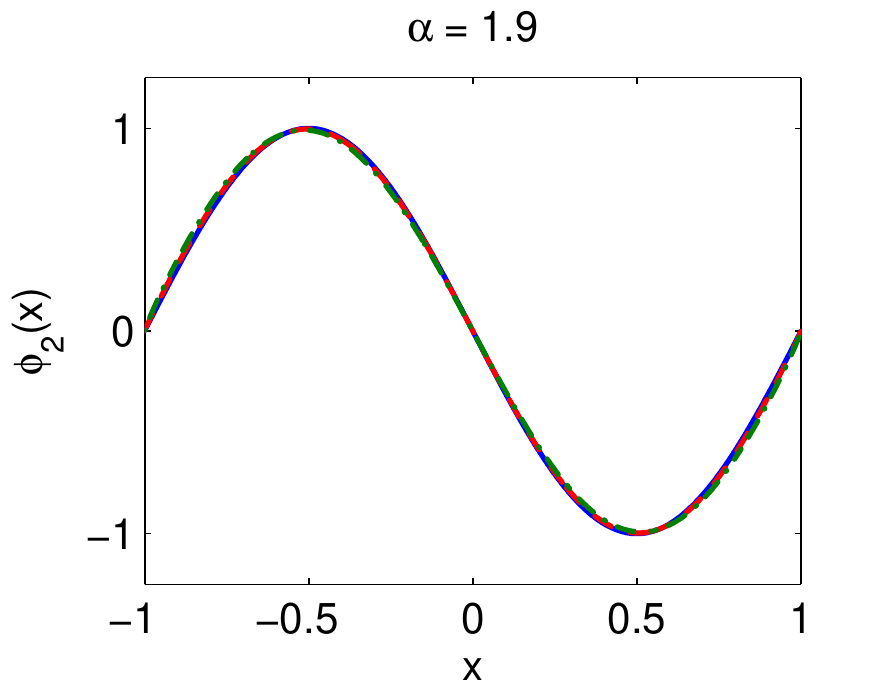}
}
\caption{The first (left panel) and second (right panel) eigenfunctions of the spectral fractional Laplacian ${\mathcal L}_s$ (solid line),  the fractional Laplacian ${\mathcal L}_h$ (dashed line),  and the regional fractional Laplacian ${\mathcal L}_r$ (dash-dot line). Note that the eigenfunctions of the spectral fractional Laplacian ${\mathcal L}_s$ are independent of $\ap > 0$. }\label{Fig6}
\end{figure}

In Figure \ref{Fig6}, we compare the first and second eigenfunctions of the fractional Laplacian, the spectral fractional Laplacian, and the regional fractional Laplacian. 
For any $\ap \in(0, 2)$, the eigenfunctions for these three operators are all symmetric (for odd $k$) or antisymmetric (for even $k$) with respect to the center of the domain $\Og$.  
Especially, the eigenfunctions of the spectral fractional Laplacian are independent of the parameter $\ap$,  which are also the eigenfunctions of the standard Dirichlet Laplace operator $-\Dt$.  
In contrast, the eigenfunctions of the other two operators significantly depend on $\alpha$, and as $\ap \to 2$, they converge to $\sin(k\pi(1+x)/2)$ -- the eigenfunctions of the standard Dirichlet Laplace operator $-\Dt$. 
Our numerical observations in Figure \ref{Fig6} justify the regularity results in \cite[Theorem 1]{Servadei2014}, that is, the eigenfunctions of the Dirichlet fractional Laplacian is no better than H\"{o}lder continuous up to the boundary, while the eigenfunctions of the spectral fractional Laplacian are smooth up to the boundary as the boundary allows.

From our extensive studies, we find that the eigenvalues of the Dirichlet fractional Laplacian ${\mathcal L}_h$, the spectral fractional Laplacian ${\mathcal L}_s$, and the regional fractional Laplacian ${\mathcal L}_r$ reduce as the domain size increases.  
In particular, if the domain size increases by a ratio of $\kappa$,  the $k$-th (for $k\in{\mathbb N}$) eigenvalues decreases by a ratio of $\kappa^\ap$.   
 Additionally, we explore the eigenvalues of the peridynamic operators for different $\dt$. 
It shows that the $k$-th (for $k \in {\mathbb N}$) eigenvalue increases as $\delta$ decreases. 

\subsection{Solutions to nonlocal  problems }
\label{section3-3}

In this section, we will further compare the four operators by studying the solutions of their corresponding nonlocal Poisson and diffusion problems. 
In the following, we choose the one-dimensional domain $\Omega = (-1,1)$ and consider homogeneous Dirichlet boundary conditions, i.e., 
\bea\label{bcc1}
u(x) = 0, \qquad \mbox{for} \ \ x\in\Gamma,
\eea
where $\Gamma = \Og^c$ for the fractional Laplacian, $\Gamma = \Og_\dt$ for the peridynamic operator,  and  $\Gamma = \p\Og$ for the spectral fractional Laplacian and the regional fractional Laplacian with $\ap\in(1, 2)$.

\vskip 10pt
\noindent {\bf Example 3. \ }  We consider a time-independent problem of the following form:
\bea\label{prob}
{\mathcal L}_i u(x) = 1,\qquad \mbox{for} \ \  x\in\Og, 
\eea
subject to the homogeneous Dirichlet boundary conditions as discussed in (\ref{bcc1}), where $i = h,\,s,\, r$,\, or $p$. 
This problem is often used as a benchmark to test numerical methods for the fractional Laplacian and spectral fractional Laplacian \cite{DuoJuZhang0016, DuoWykZhang,Delia2013}. 
In particular, the nonlocal problem (\ref{prob}) with the fractional Laplacian ${\mathcal L}_h$ has diverse applications in various areas  \cite{Dyda2012}. 
Its solution can be found exactly as: 
\begin{eqnarray*}
\label{u_example2}
u(x) = \frac{1}{\Gamma(\alpha +1)}(1-x^2)^{\alpha/2}, \quad\ \ \text{for}\ \ x\in\Og,
\end{eqnarray*}
which represents the probability density function of the first exit time of the symmetric $\alpha$-stable L\'{e}vy process from  the domain $\Og$.
While if the spectral fractional Laplacian ${\mathcal L}_s$ is considered,  the exact solution of (\ref{prob}) is expressed as
\beas
u(x) = \sum_{k = 0}^\infty \fl{2\big(1-(-1)^k\big)}{k\pi}\Big(\fl{k\pi}{2}\Big)^{-\ap} \sin\Big(\fl{k\pi}{2}(1+x)\Big), \quad\ \ \text{for}\ \ x\in\Og.
\eeas
We will compute the solutions of (\ref{prob}) with the regional fractional Laplacian ${\mathcal L}_r$ or  peridynamic operator ${\mathcal L}_p$ numerically. 

\begin{figure}[htb!]
\centerline{\includegraphics[height=5.1060cm,width=6.960cm]{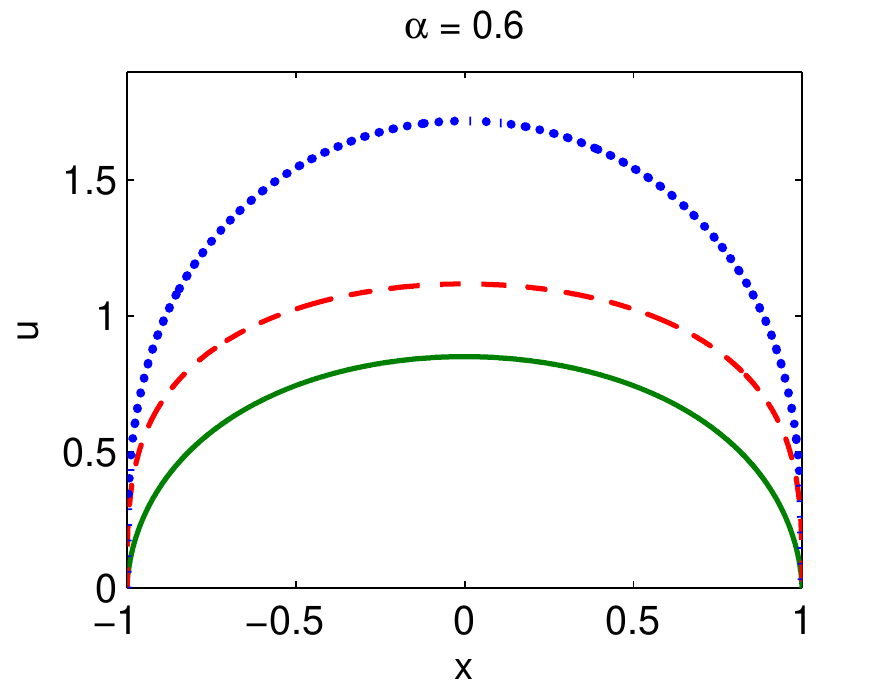}
\includegraphics[height=5.1060cm,width=6.960cm]{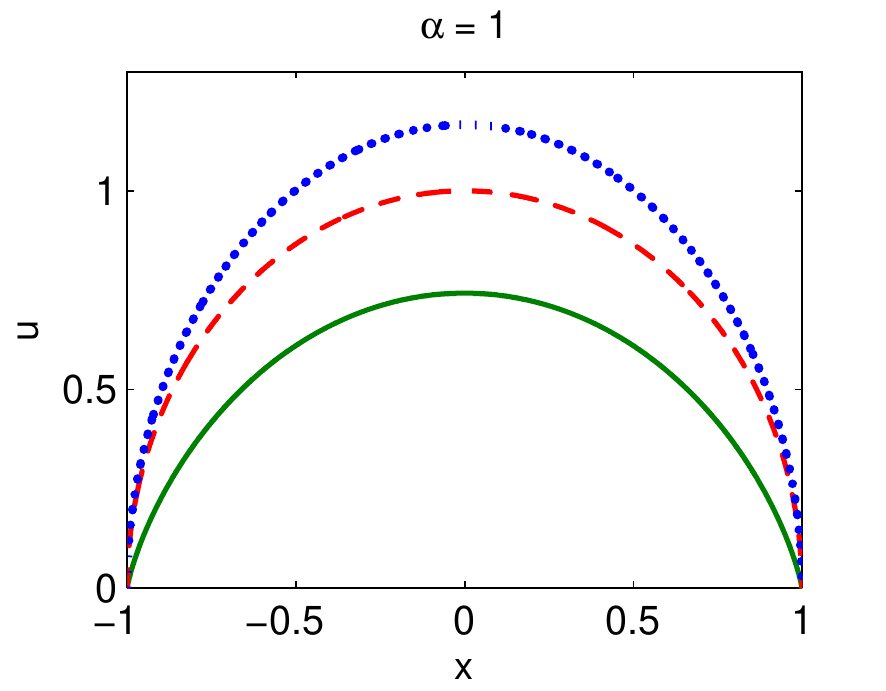}}
\centerline{
\includegraphics[height=5.1060cm,width=6.960cm]{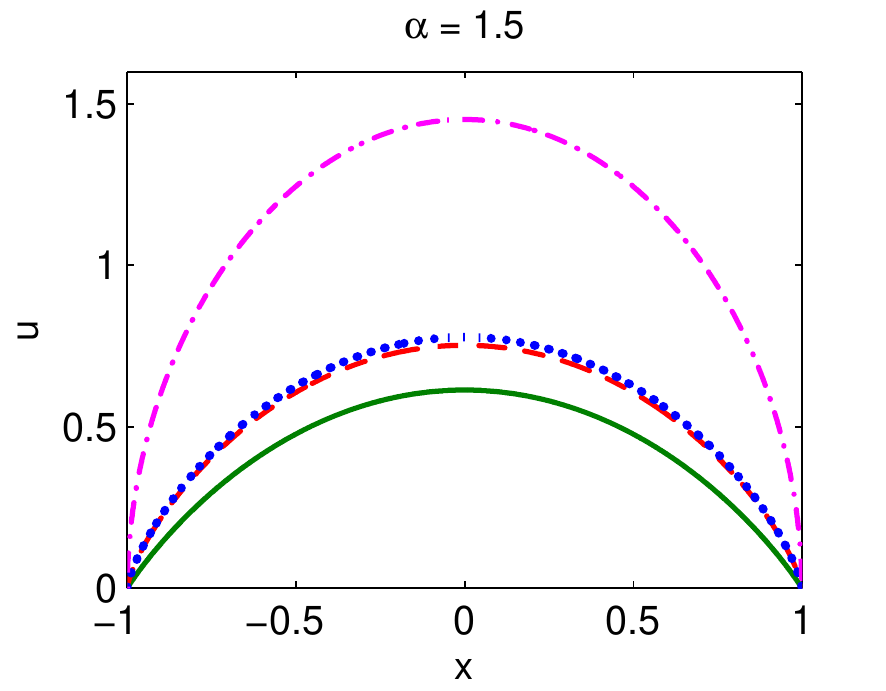}
\includegraphics[height=5.1060cm,width=6.960cm]{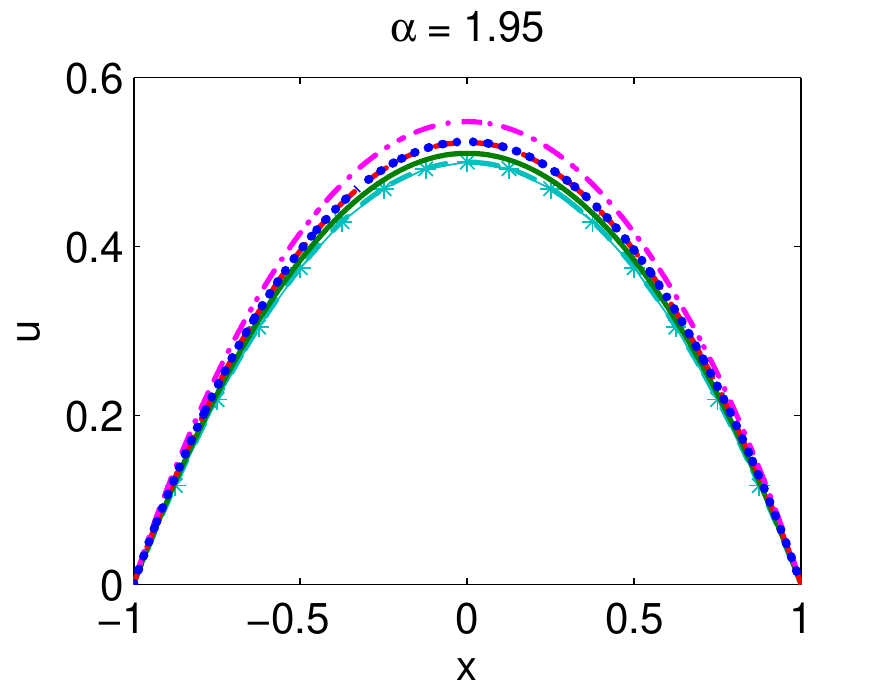}}
\caption{Comparison of the solution to (\ref{prob}) with ${\mathcal L}_s$ (solid line),  ${\mathcal L}_h$ (dashed line), ${\mathcal L}_r$ (dash-dot line),  or ${\mathcal L}_p$ with $\dt = 4$ (dotted line).  For easy comparison, the result for ${\mathcal L} = -\p_{xx}$ (line with symbols ``*") is included in the plot of $\ap = 1.95$. }
\label{Fig7}
\end{figure}

Figure \ref{Fig7} illustrates the solutions of the nonlocal problem (\ref{prob}) with different operators ${\mathcal L}_i$. 
In the cases of the regional fractional Laplacian ${\mathcal L}_r$,  since the solution does not exist for $\ap \leq 1$,  we only present those for $\ap > 1$. 
Generally,   the solutions from the four operators are significantly different, but as $\ap \to 2$  they all  converge to the function $u(x) = \fl{1}{2}(1-x^2)$ -- the solution to the classical  Poisson equation: 
\bea\label{poisson}
-\p_{xx} u(x) = 1, \quad\mbox{if} \ \ x\in\Og; \qquad u = 0, \quad\mbox{if} \ \ x\in\p\Og.
\eea
The solution from the regional fractional Laplacian ${\mathcal L}_r$ is very sensitive to the parameter $\ap$, and moreover it is inconsistent with that from the fractional Laplacian ${\mathcal L}_h$. 
The solution of the peridynamic model could serve a good approximation to that of the fractional Poisson problem (\ref{fPoisson0})--(\ref{fBC0}), if a proper horizon size $\dt$ is chosen.  
The smaller the parameter $\ap$ is, the larger the horizon size $\dt$ is needed, consistent with our observations in Section \ref{section3-1}. 

\begin{figure}[h!]
\centerline{\includegraphics[height=5.1060cm,width=6.960cm]{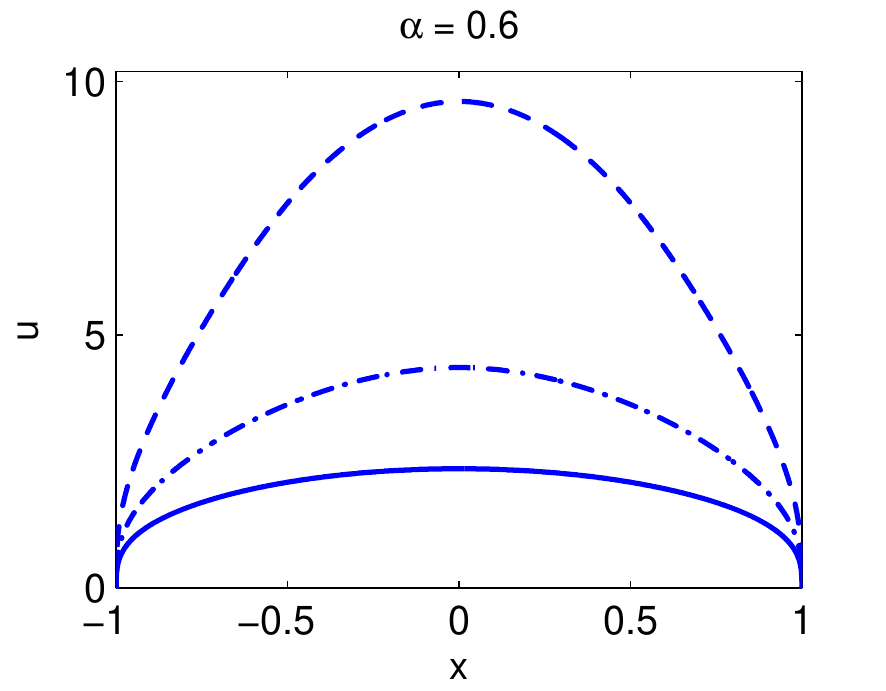}\quad
\includegraphics[height=5.1060cm,width=6.960cm]{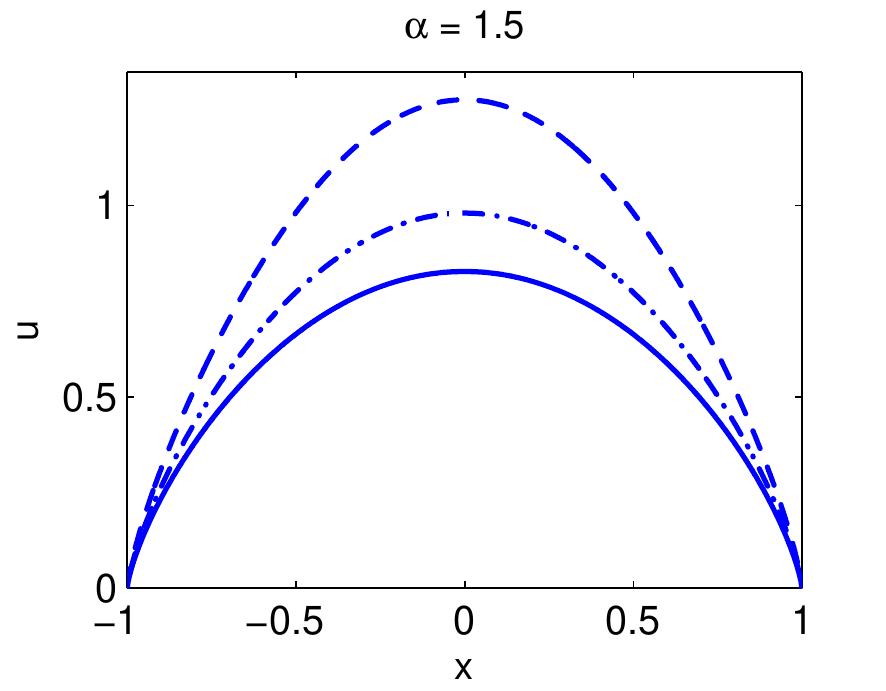}}
\caption{Effects of the horizon size $\dt$ on the solution of the nonlocal problem (\ref{prob}) with the peridynamic operator ${\mathcal L}_p$, where $\dt = 2$ (solid line), $1$ (dash-dot line), or $0.5$ (dashed line). }\label{fig00}
\end{figure}

In Figure \ref{fig00}, we additionally study the solution of (\ref{prob}) with the peridynamic operator ${\mathcal L}_p$ for various horizon size $\dt$.
It shows that the horizon size $\dt$ plays an important role in the solution of peridynamic models, especially when $\ap$ is small (see Fig. \ref{fig00} left).  
For the same $\ap$, the larger the horizon size $\dt$, the smaller the solution $u$.  
{ Our extensive studies show that as $\dt \to 0$, the solution of the peridynamic models converges to the function $u(x) = C(1-x^2)$ with $C$ a constant depending on $\ap$, which can be viewed as a rescaled solution of the classical Poisson problem in (\ref{poisson}).}
\vskip 10pt

\noindent {\bf Example 4.\ }  Here, we study the following nonlocal diffusion problem:
\begin{eqnarray}\label{fDiffusion}
\p_tu(x, t) = \mathcal{L}_iu, \qquad \mbox{for} \ \  x\in\Og, 
\eea
subject to the homogeneous Dirichlet boundary conditions as discussed in (\ref{bcc1}), where $i = h, \,s,\, r$, or $p$. 
At time $t = 0$, the initial condition is taken as a step function, i.e., 
\bea \label{initial5}
u(x,0) = \left\{\begin{array}{rcl}
1, & & \mbox{if} \  \ x\in (-0.2,0.2),\\
0, &  & \mbox{otherwise}, 
\end{array}\right.\quad \  \ x\in \mathbb{R}.
\eea
If the spectra of the nonlocal operator ${\mathcal L}_i$ are known, one can formally express the solution of (\ref{fDiffusion})  as: 
\begin{eqnarray}\label{fDiffusion_Solution}
u(x,t) = \sum _{k=1}^{\infty} c_k e^{-\lambda_k^i t} \phi_k^i(x),\qquad \mbox{for} \ \ x\in\Og, \quad t \ge 0,
\end{eqnarray}
where $\lambda _k^i$ and $\phi_k^i$ represent the $k$-th eigenvalue and eigenfunction of the  operator $\mathcal{L}_i$ on the domain $\Og$, and the coefficient $c_k$ is calculated by 
\begin{eqnarray}\label{ck}
c_k = \int _{\Omega} u(x, 0)\phi_k(x)dx.
\end{eqnarray}
Hence, the solution of the diffusion equation (\ref{fDiffusion}) with the spectral fractional Laplacian ${\mathcal L}_s$ can be obtained analytically as 
\beas
u(x,t) = \fl{2}{\pi}\sum _{k=1}^{\infty}\fl{\cos(2k\pi/5)-\cos(3k\pi/5)}{k}e^{-(\fl{k\pi}{2})^{\ap}t}\sin\Big(\fl{k\pi}{2}(1+x)\Big),
\eeas
for $x\in\Og$ and  $t \ge 0$.
For other cases with ${\mathcal L}_i$ ($i = h,r,p$) in (\ref{fDiffusion}), we will compute their numerical solutions by using the finite difference method proposed in \cite{DuoWykZhang}. 

\begin{figure}[htb!]
\centerline{\includegraphics[height=3.960cm,width=5.260cm]{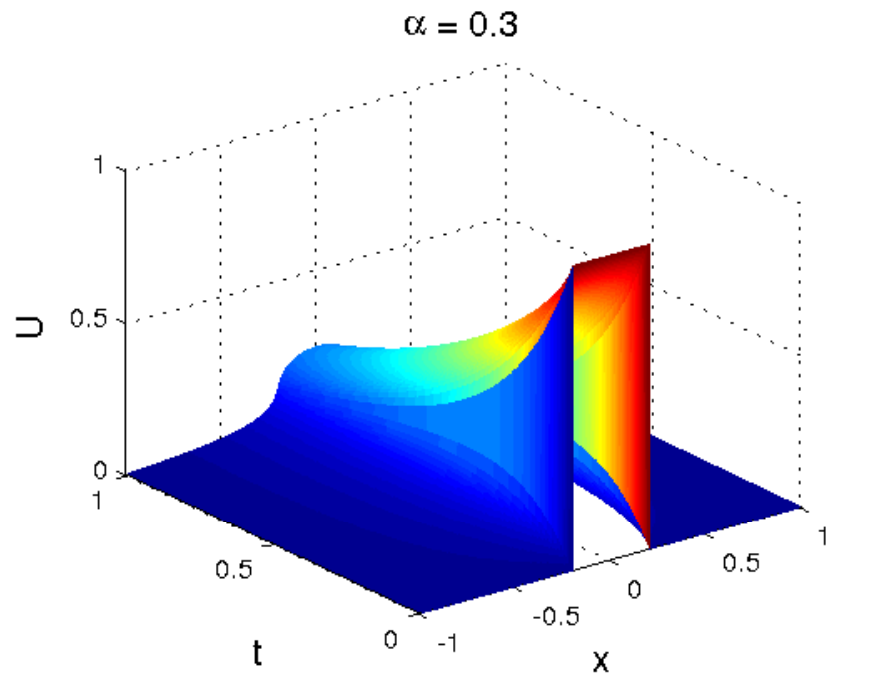}\hspace{-2mm}
\includegraphics[height=3.960cm,width=5.260cm]{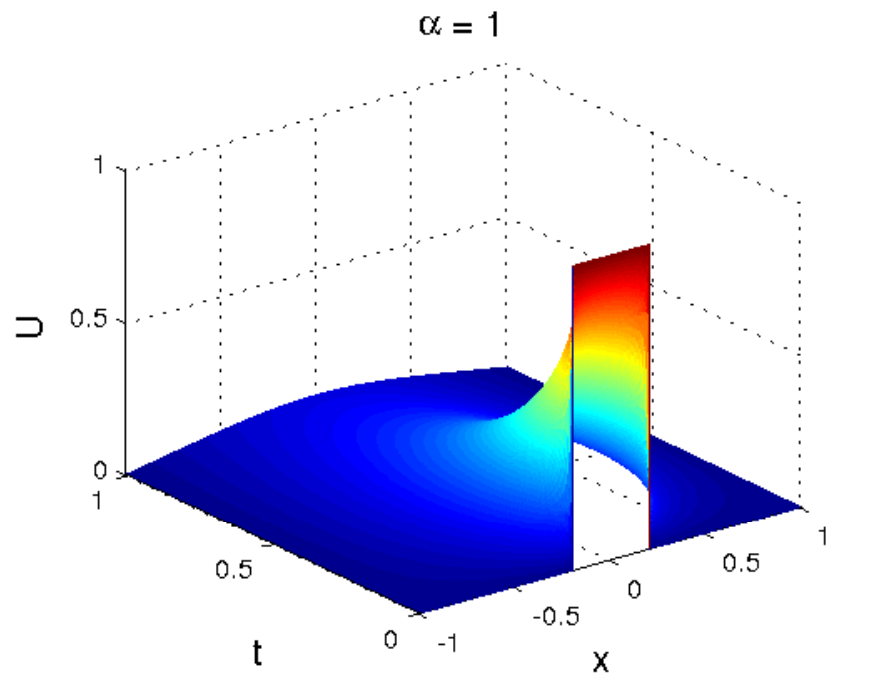} \hspace{-2mm}
\includegraphics[height=3.960cm,width=5.260cm]{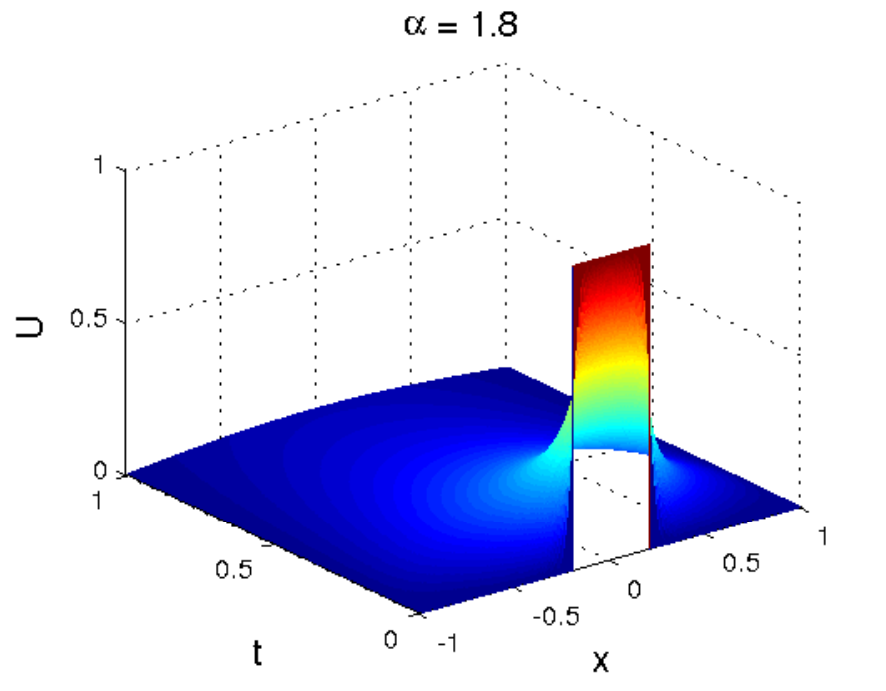}}
\centerline{\includegraphics[height=3.960cm,width=5.260cm]{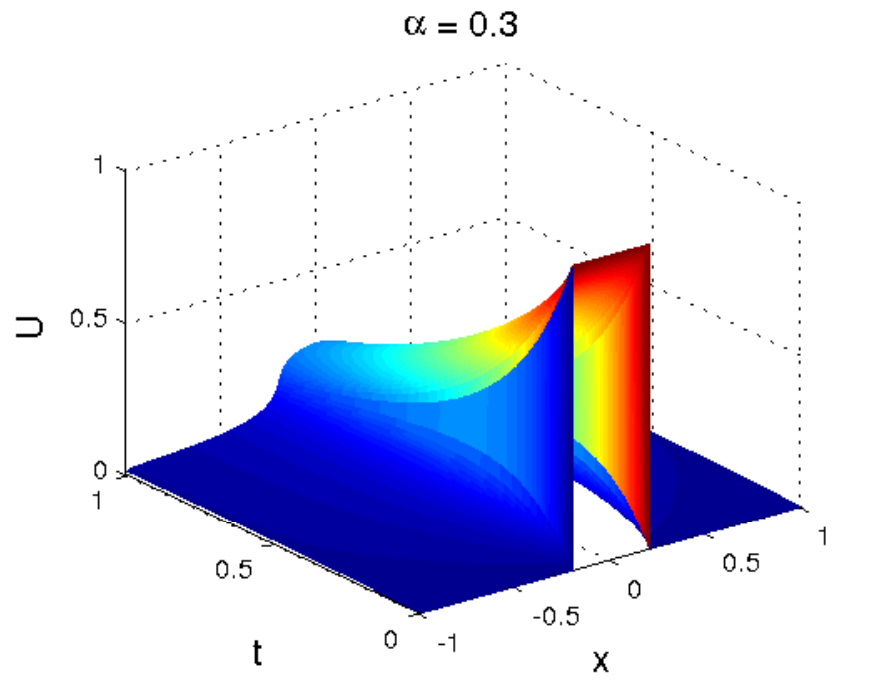}\hspace{-2mm}
\includegraphics[height=3.960cm,width=5.260cm]{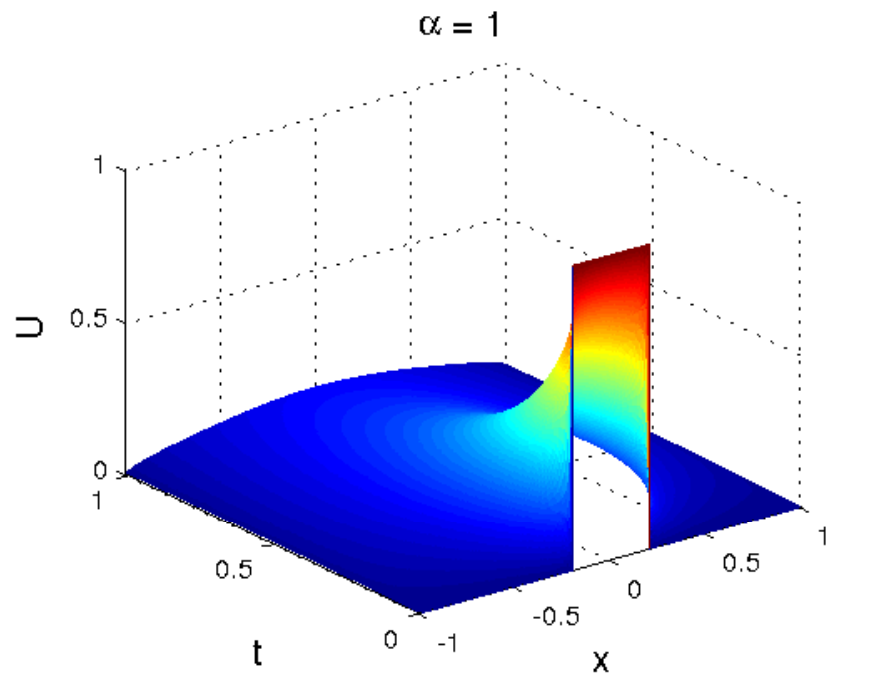} \hspace{-2mm}
\includegraphics[height=3.960cm,width=5.260cm]{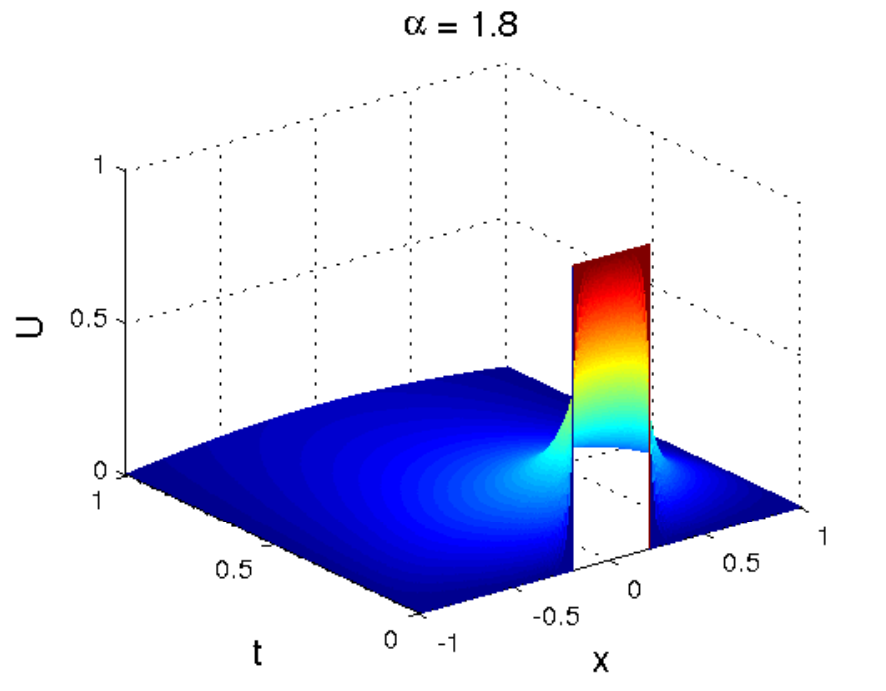}}
\centerline{\includegraphics[height=3.960cm,width=5.260cm]{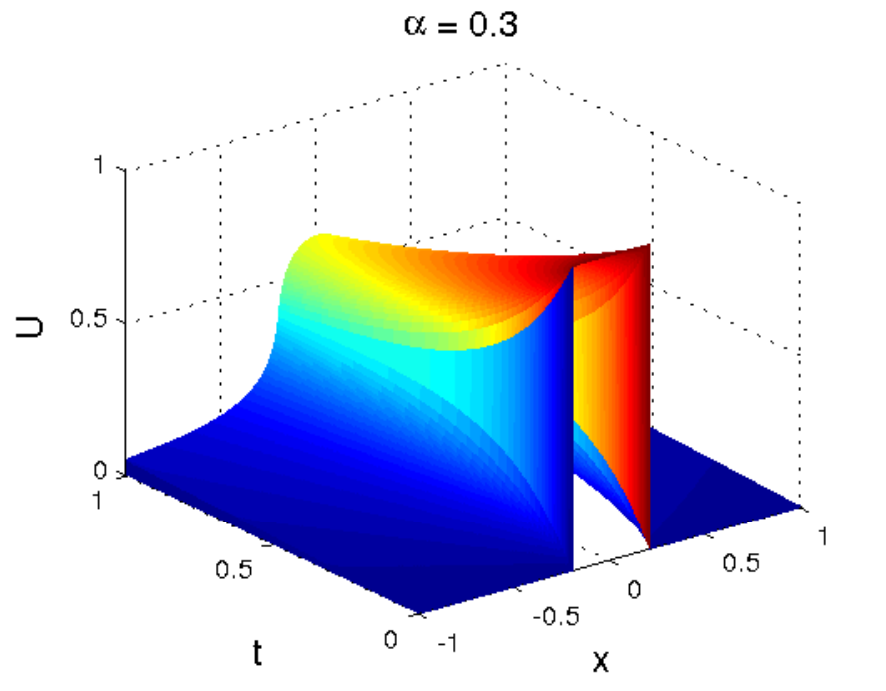}\hspace{-2mm}
\includegraphics[height=3.960cm,width=5.260cm]{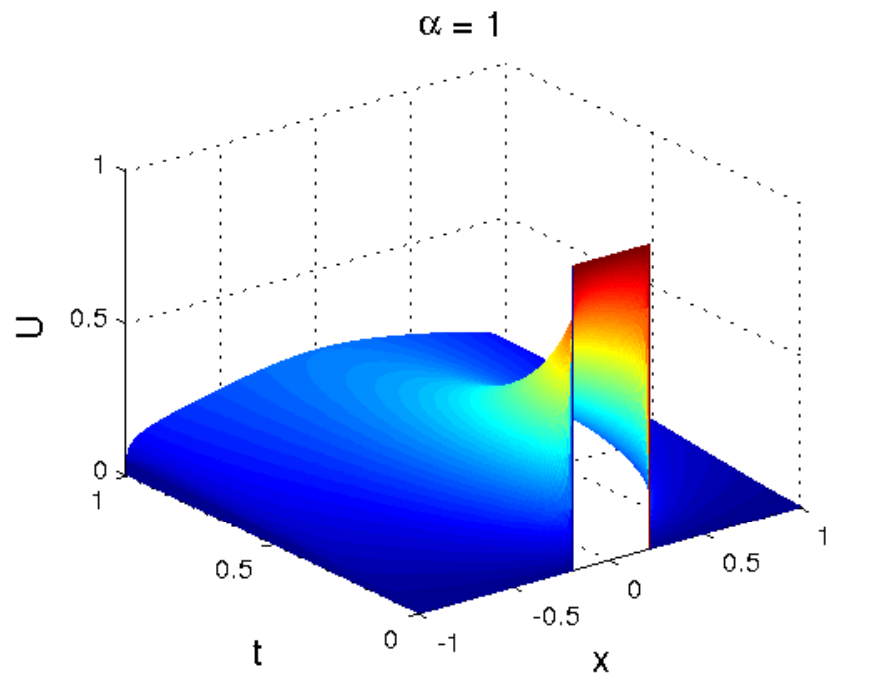} \hspace{-2mm}
\includegraphics[height=3.960cm,width=5.260cm]{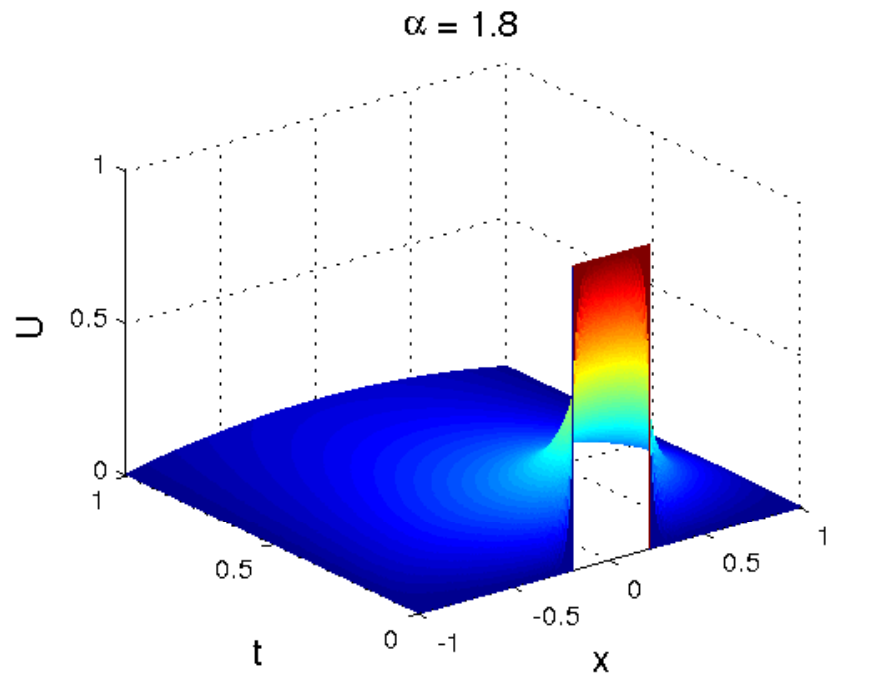}}
\caption{Time evolution of the solution $u(x, t)$ to the nonlocal diffusion equation (\ref{fDiffusion}) with ${\mathcal L}_s$ (upper row), ${\mathcal L}_h$ (middle row), and ${\mathcal L}_r$ (lower row).}\label{fig10}
\end{figure}
Figure \ref{fig10}  shows the time evolution of the solution $u(x, t)$ for  the cases of the spectral fractional Laplacian ${\mathcal L}_s$, fractional Laplacian ${\mathcal L}_h$, and regional fractional Laplacian ${\mathcal L}_r$, while the results for the peridynamic operator ${\mathcal L}_p$ are displayed in Figure \ref{fig11} for various horizon size $\dt$. 
\begin{figure}[htb!]
\centerline{
\includegraphics[height=3.960cm,width=5.260cm]{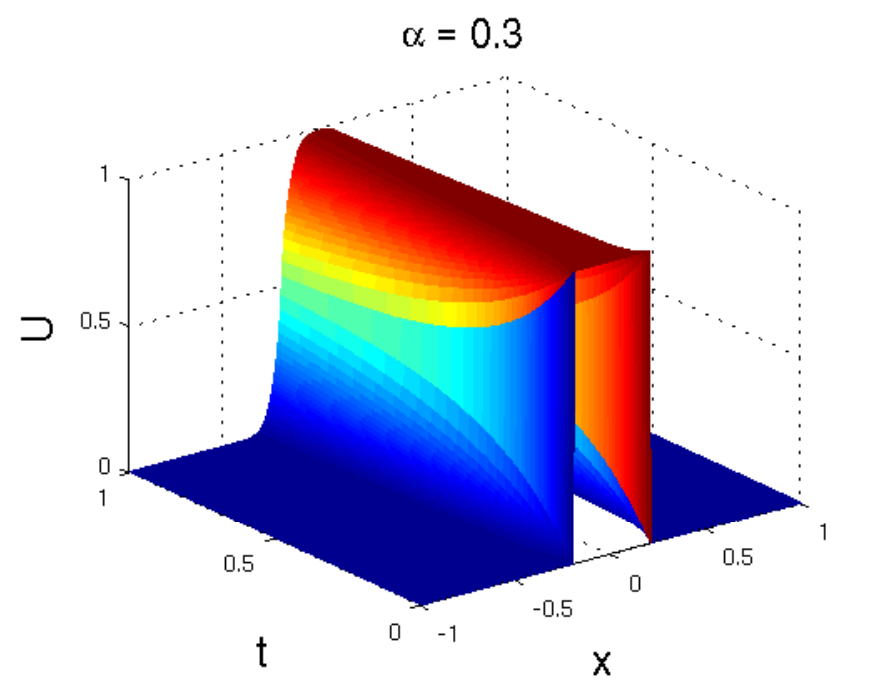}\hspace{-2mm}
\includegraphics[height=3.960cm,width=5.260cm]{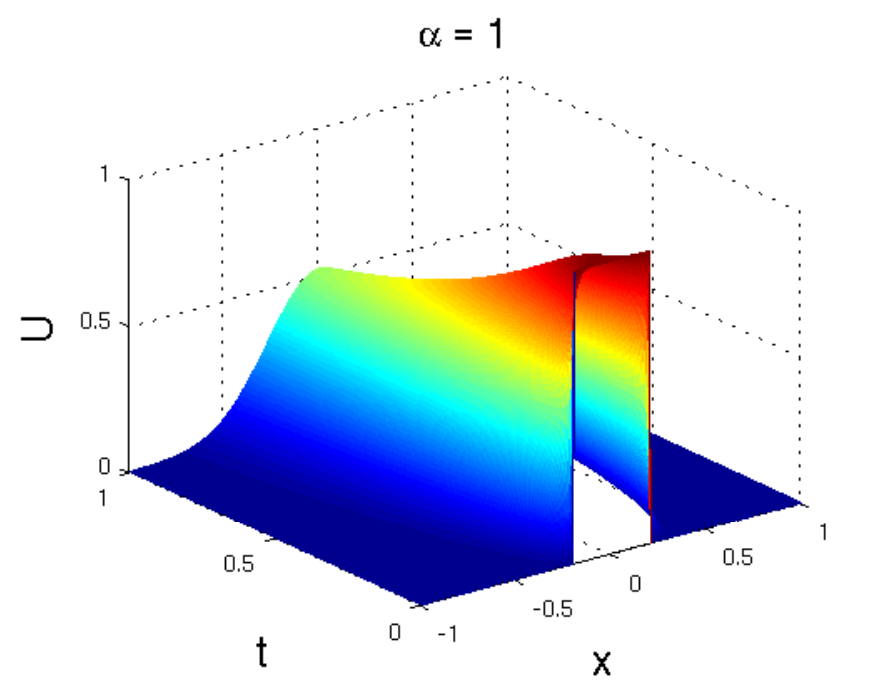} \hspace{-2mm}
\includegraphics[height=3.960cm,width=5.260cm]{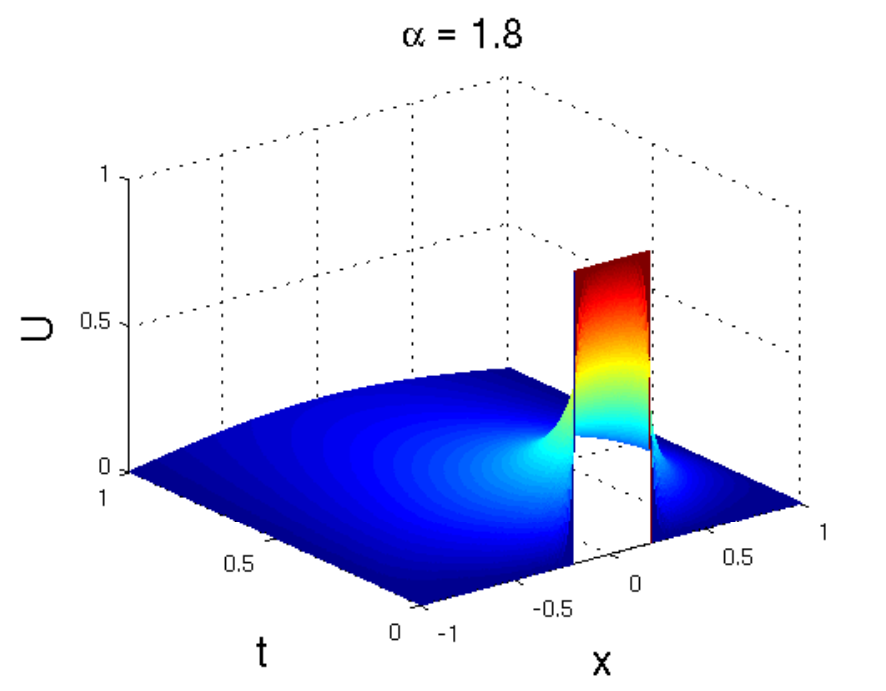}}
\centerline{
\includegraphics[height=3.960cm,width=5.260cm]{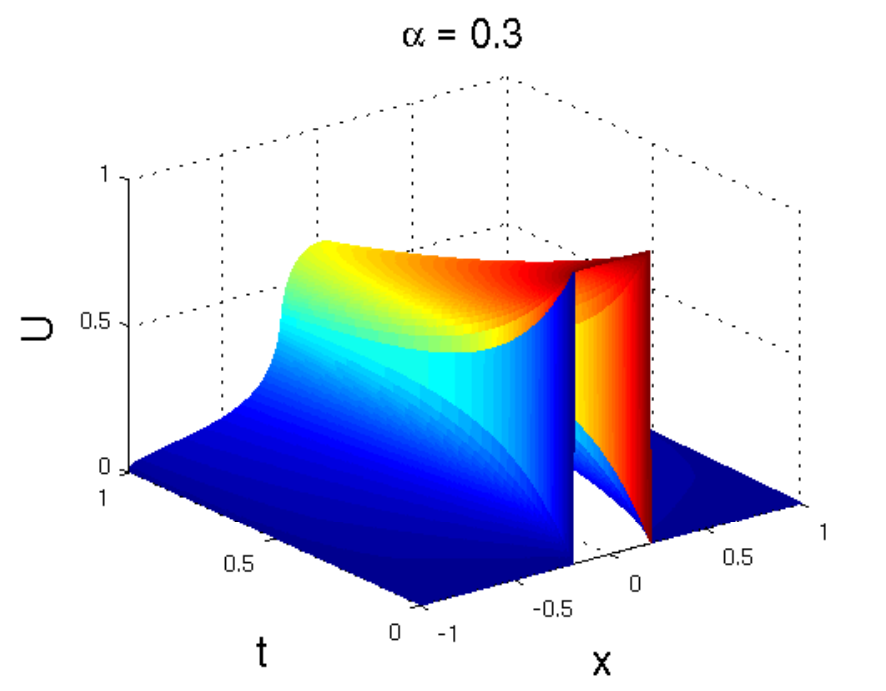}\hspace{-2mm}
\includegraphics[height=3.960cm,width=5.260cm]{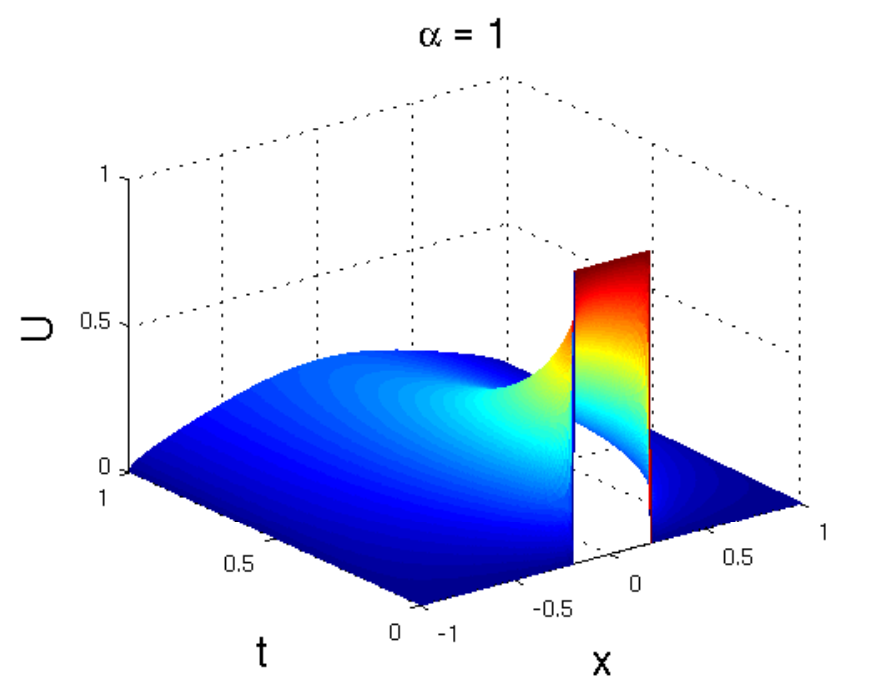} \hspace{-2mm}
\includegraphics[height=3.960cm,width=5.260cm]{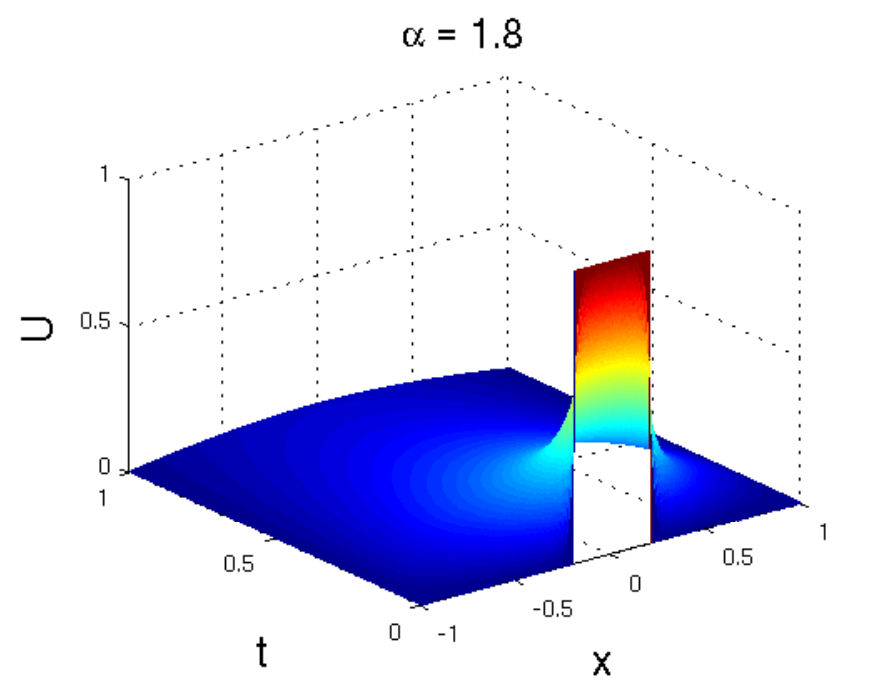}}
\caption{Time evolution of the solution $u(x, t)$ to the nonlocal diffusion equation (\ref{fDiffusion}) 
with the peridynamic operator ${\mathcal L}_p$,  where the horizon size $\dt = 0.1$ (top) or $\dt = 1$ (bottom).}\label{fig11}
\end{figure}
The similar phenomena to the classical diffusion equation are observed  -- the solution of nonlocal diffusion equation (\ref{fDiffusion}) decays over time, and the discontinuous initial condition smooths out quickly during the dynamics. 
For the same operator ${\mathcal L}_i$, the larger the parameter $\alpha$ is, the faster the solution diffuses. 
Moreover,  Fig. \ref{fig11} shows that for the peridynamic case, the diffusion speed of the solution depends not only on the parameter $\ap$, but also on the horizon size $\dt$. 
For fixed $\ap$, the larger the horizon size, the faster the diffusion. 

\begin{figure}[htb!]
\centerline{
\includegraphics[height=4.060cm,width=5.960cm]{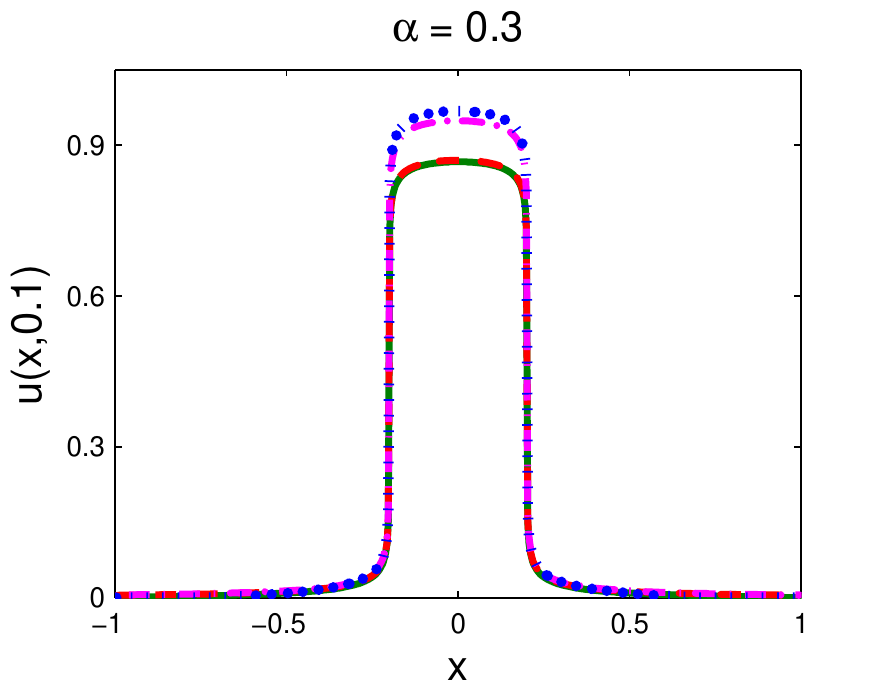}\hspace{-5mm}
\includegraphics[height=4.060cm,width=5.96cm]{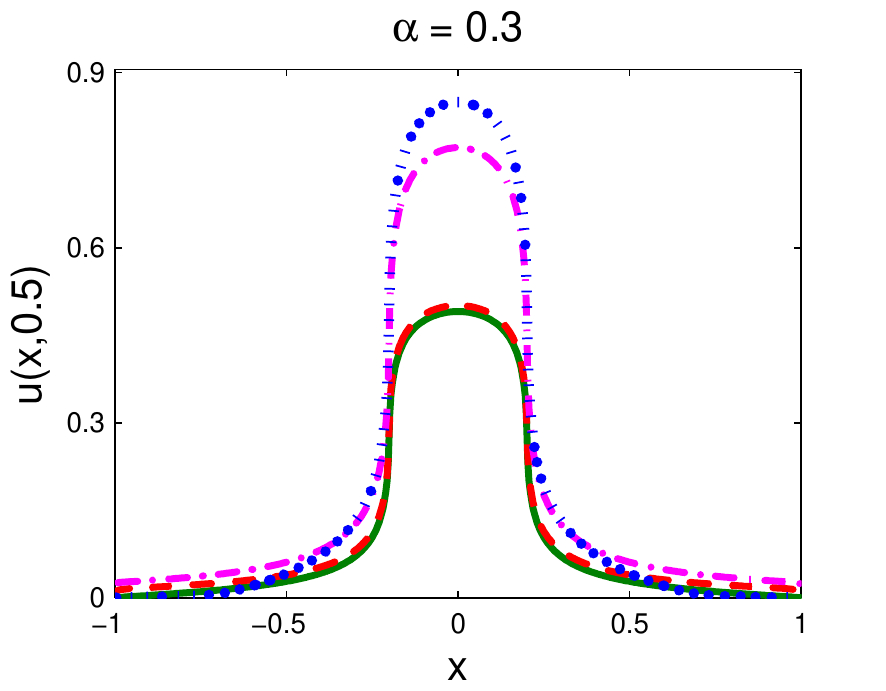}\hspace{-5mm}
\includegraphics[height=4.060cm,width=5.960cm]{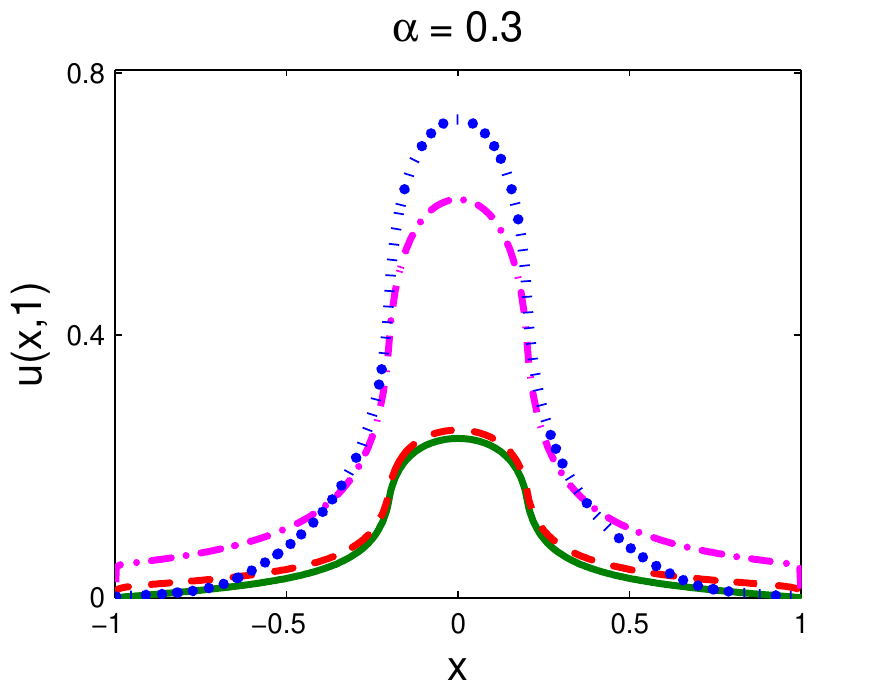}}
\centerline{
\includegraphics[height=4.060cm,width=5.960cm]{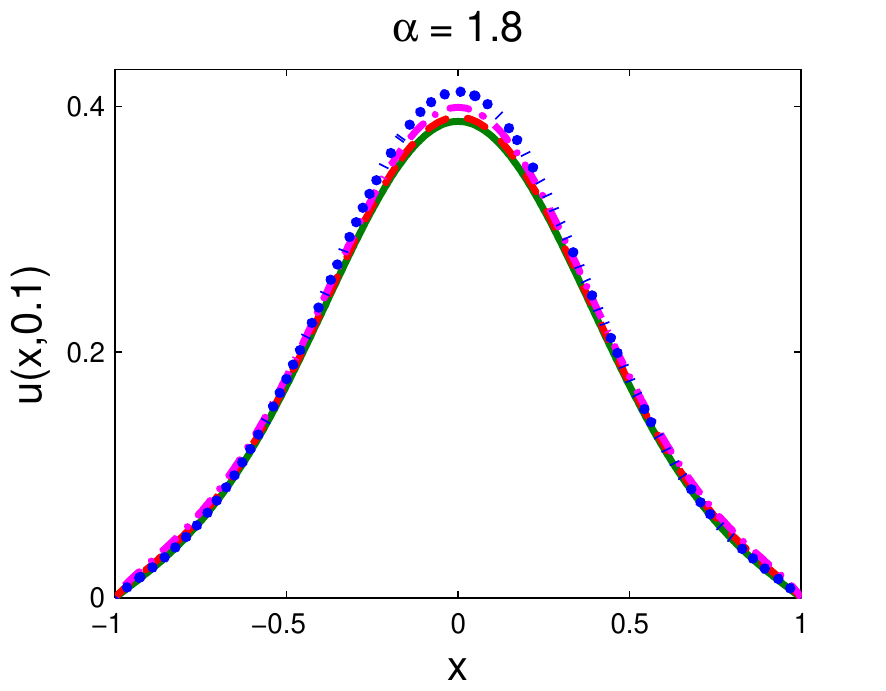}\hspace{-5mm}
\includegraphics[height=4.060cm,width=5.960cm]{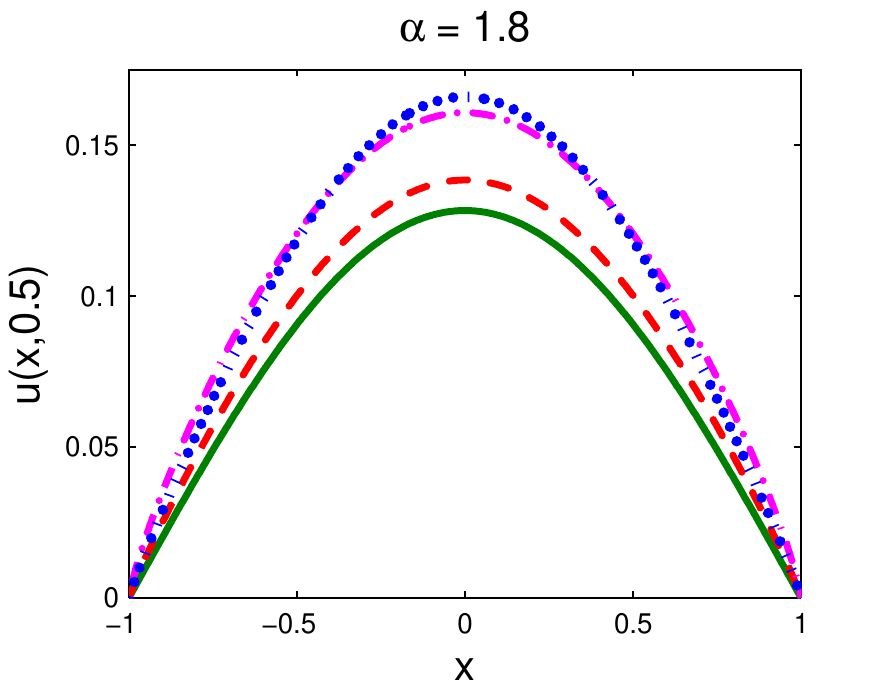} \hspace{-5mm}
\includegraphics[height=4.060cm,width=5.960cm]{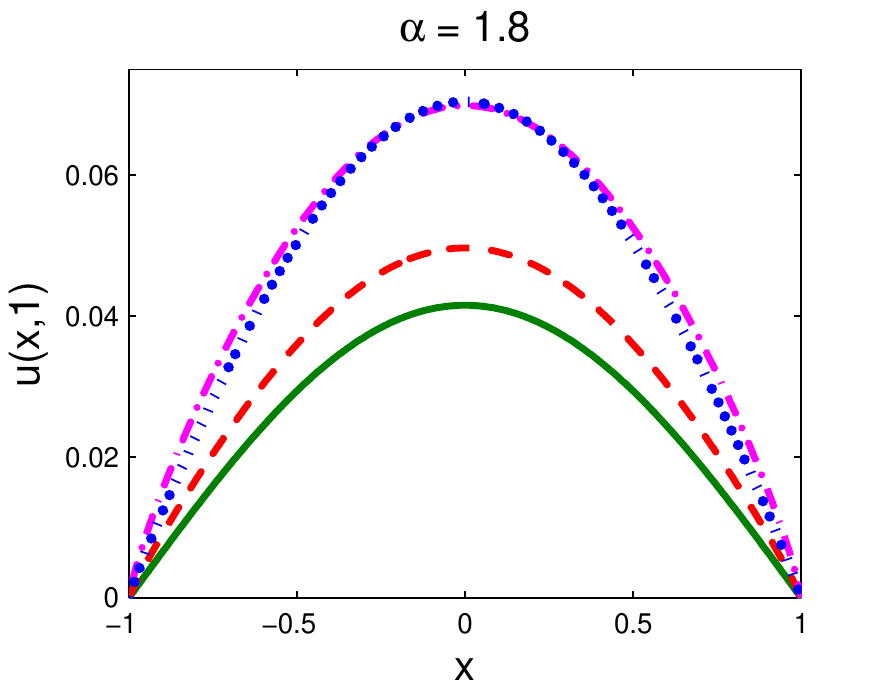}}
\centerline{
\includegraphics[height=4.060cm,width=5.960cm]{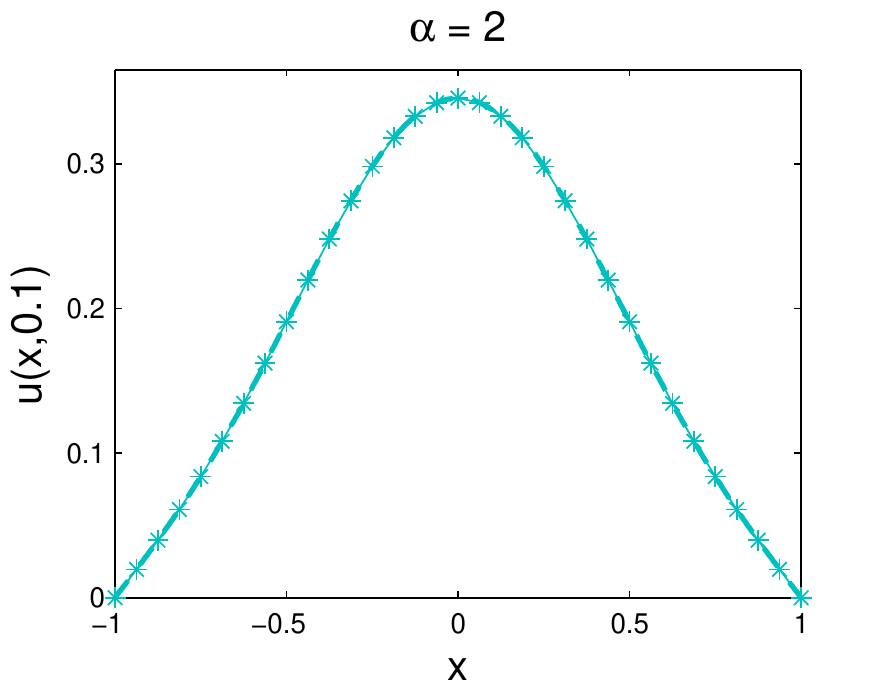}\hspace{-5mm}
\includegraphics[height=4.060cm,width=5.960cm]{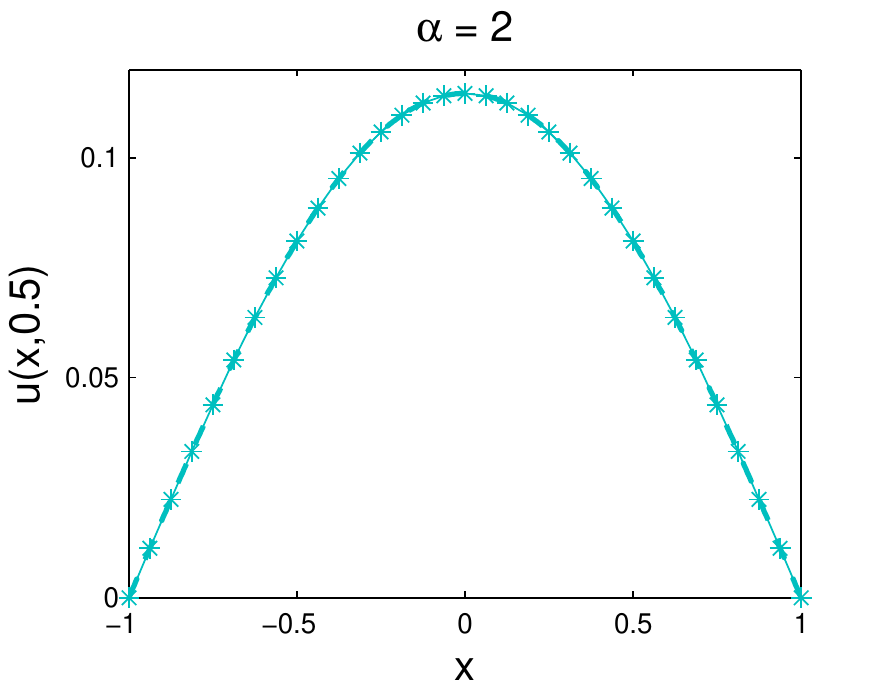} \hspace{-5mm}
\includegraphics[height=4.060cm,width=5.960cm]{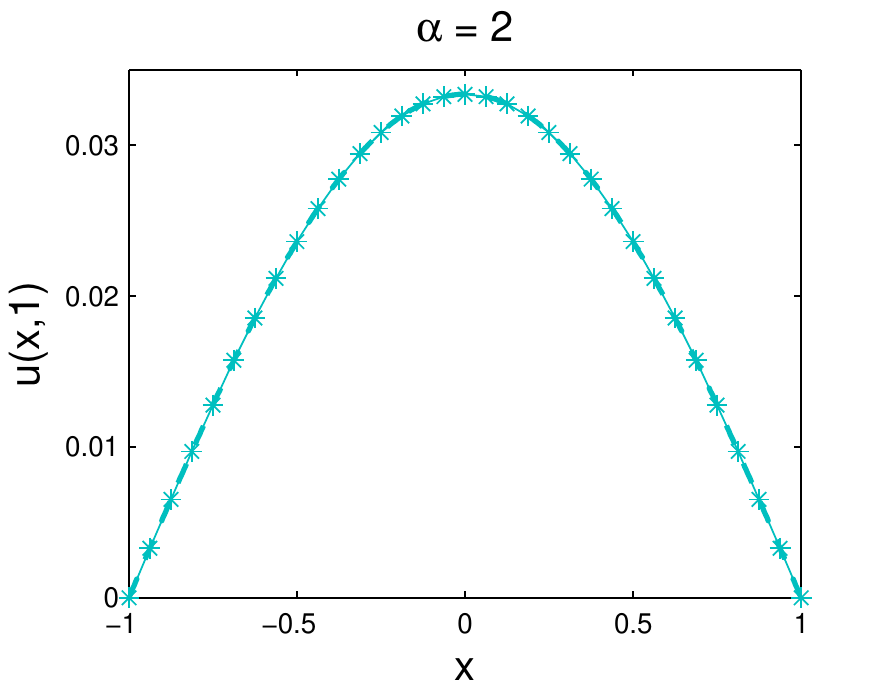}}
\caption{Solutions of the nonlocal diffusion equation (\ref{fDiffusion}) at time $t = 0.1, 0.5, 1$, where the operator is chosen as ${\mathcal L}_s$ (solid line),  ${\mathcal L}_h$ (dashed line), ${\mathcal L}_p$ (dotted line), or ${\mathcal L}_r$ (dash-dot line). 
For easy comparison, we include the solution of the classical diffusion equation (i.e., ${\mathcal L}_i = -\p_{xx}$ in (\ref{fDiffusion})) in the last row.}\label{fig12}
\end{figure}

In Figure \ref{fig12}, we further compare the solutions at various time $t$, where the horizon size $\dt = 0.5$ is used in the peridynamic operator. 
It shows that the solution quickly diffuses to the boundary if  $\ap$ is large (see Fig. \ref{fig12} with $\ap = 1.8$).  
For the same parameter $\ap$, the solution resulting from the spectral fractional Laplacian diffuses much faster than those from other operators, which can be explained by the solution in (\ref{fDiffusion_Solution}). 
The solution (\ref{fDiffusion_Solution}) implies that the larger the eigenvalues $\lambda_k^i$, the faster the solution decays over time. 
On the other hand, Table \ref{Tab1} suggests that the eigenvalues of different operators satisfy $\lambda_k^s > \lambda_k^h > \lambda_k^r$,  for $\ap\in(0, 2)$ and $k\in{\mathbb N}$. 
Hence, it is easy to conclude that the solutions from the spectral fractional Laplacian diffuse the fastest. 
As mentioned previously,  the solution of the peridynamic model depends on the horizon size $\dt$. 

\vskip 10pt
\noindent {\bf Example 5.\ } We consider a nonlocal diffusion-reaction equation of the following form: 
 \begin{eqnarray}\label{fDiffusion1}
\p_tu(x, t) = \mathcal{L}_iu + u, \qquad \mbox{for} \ \  x\in\Og, 
\eea
subject to the homogeneous Dirichlet boundary conditions as discussed in (\ref{bcc1}), where $i = h, \,s, \, r$, or $p$. 
The initial condition is taken as
\bea \label{initial6}
u(x,0) = \exp\big[-(4x)^2\big], \qquad x\in{\mathbb R}.
\eea
which decays quickly  to zero as $|x| \to\infty$.
Similarly, the  solution of (\ref{fDiffusion1}) can be formulated as: 
\begin{eqnarray}\label{fDiffusion_Solution1}
u(x,t) = \sum _{k=1}^{\infty} c_k e^{(1-\lambda_k^i) t} \phi_k^i(x),\qquad \mbox{for} \ \ x\in\Og, \quad t \ge 0,
\end{eqnarray}
with $\lambda _k^i$ and $\phi_k^i$ denoting the $k$-th eigenvalue and eigenfunction of the  operator $\mathcal{L}_i$ on the domain $\Og$, and the coefficient $c_k$ defined in (\ref{ck}). 

Figure \ref{fig13}  shows the time evolution of the solution $u(x, t)$ to the reaction-diffusion equation 
\begin{figure}[htb!]
\centerline{
\includegraphics[height=3.960cm,width=5.260cm]{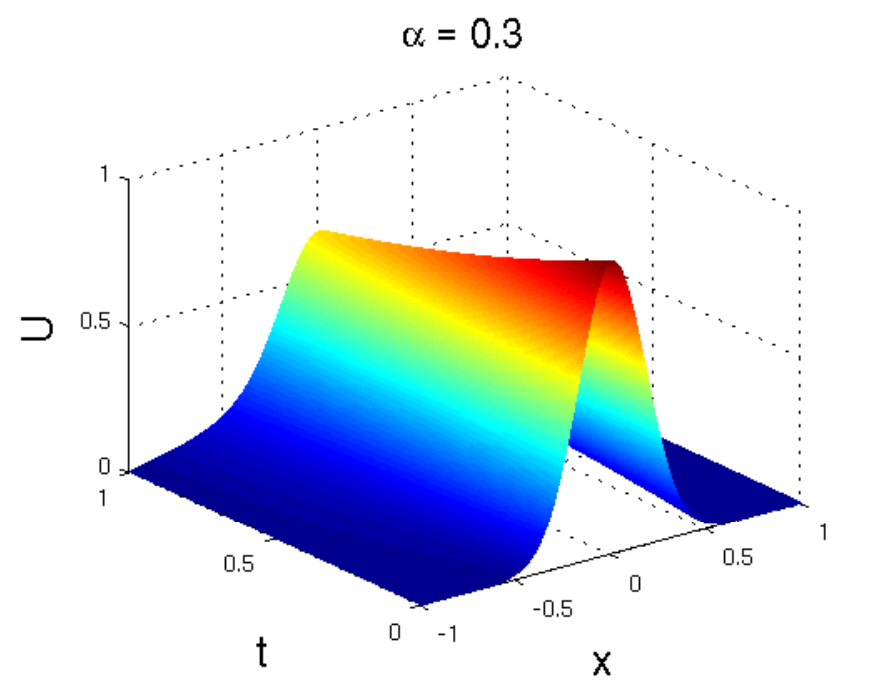}\hspace{-2mm}
\includegraphics[height=3.960cm,width=5.260cm]{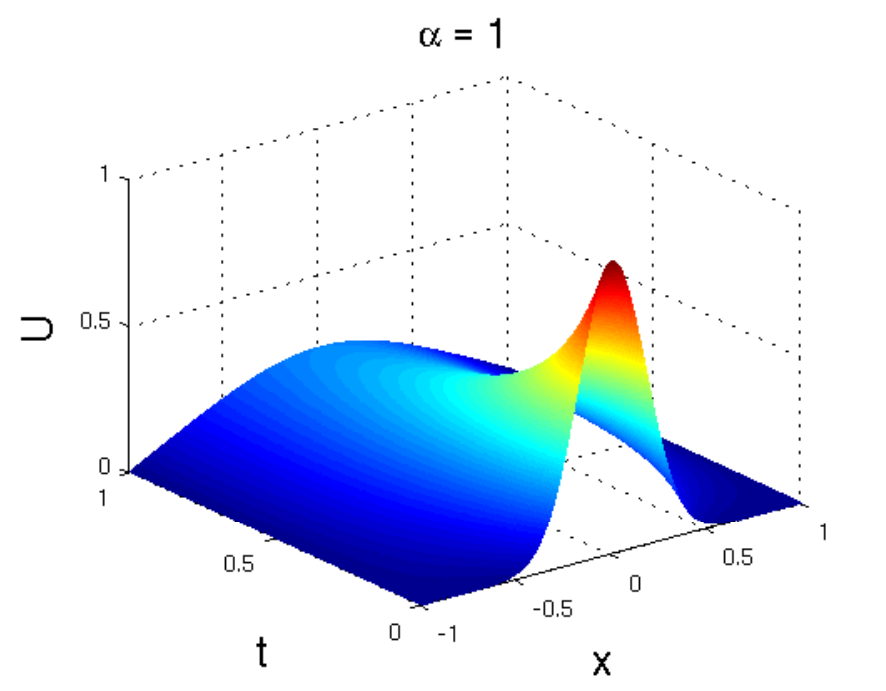}\hspace{-2mm}
\includegraphics[height=3.960cm,width=5.260cm]{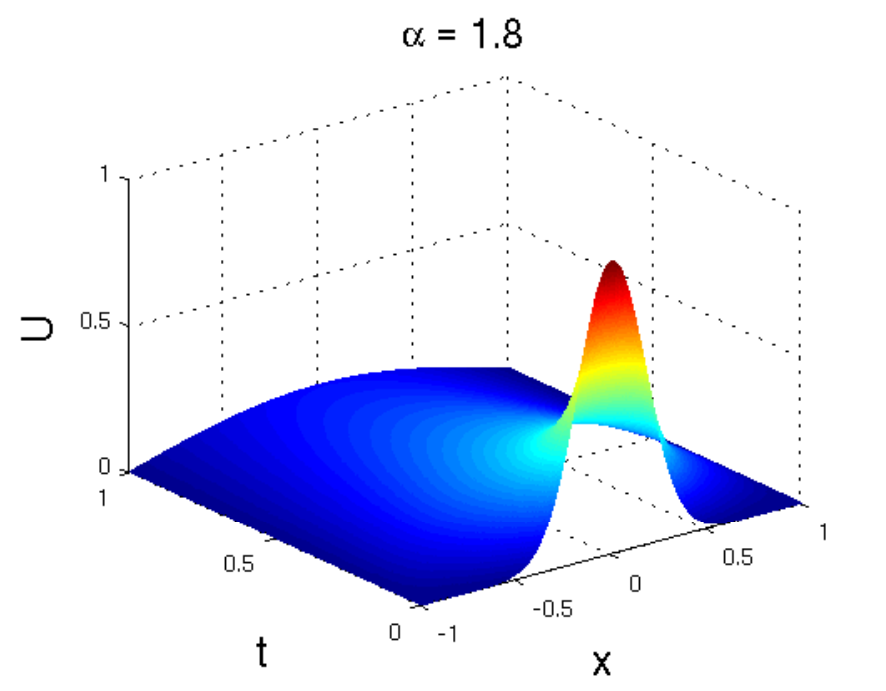}}
\centerline{
\includegraphics[height=3.960cm,width=5.260cm]{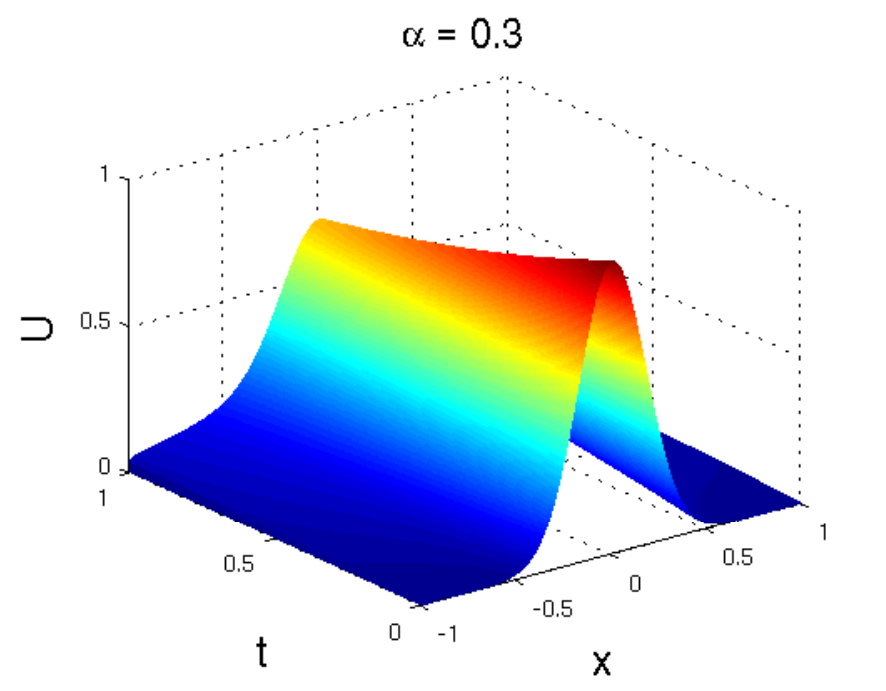}\hspace{-2mm}
\includegraphics[height=3.960cm,width=5.260cm]{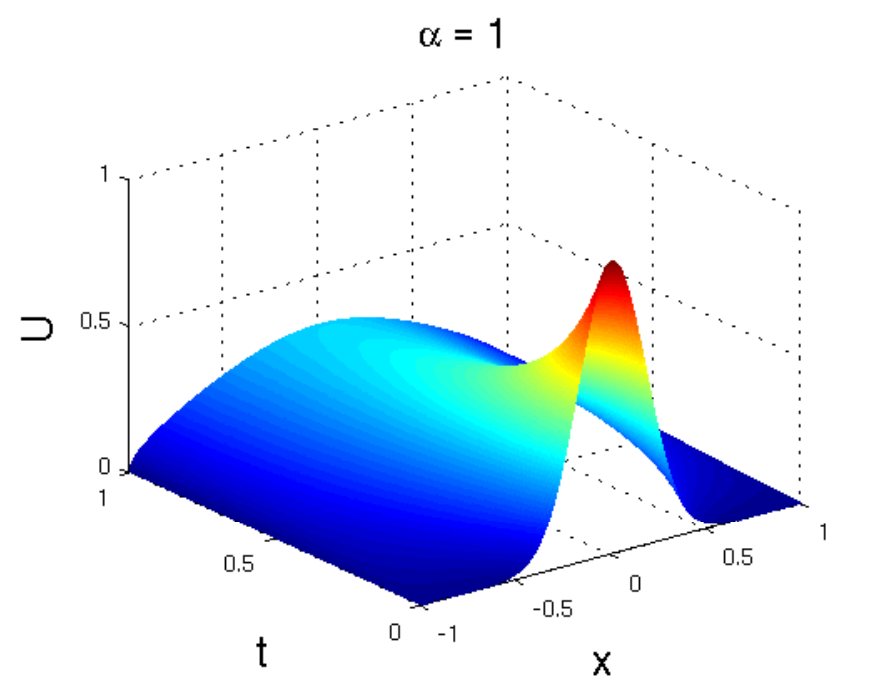}\hspace{-2mm}
\includegraphics[height=3.960cm,width=5.260cm]{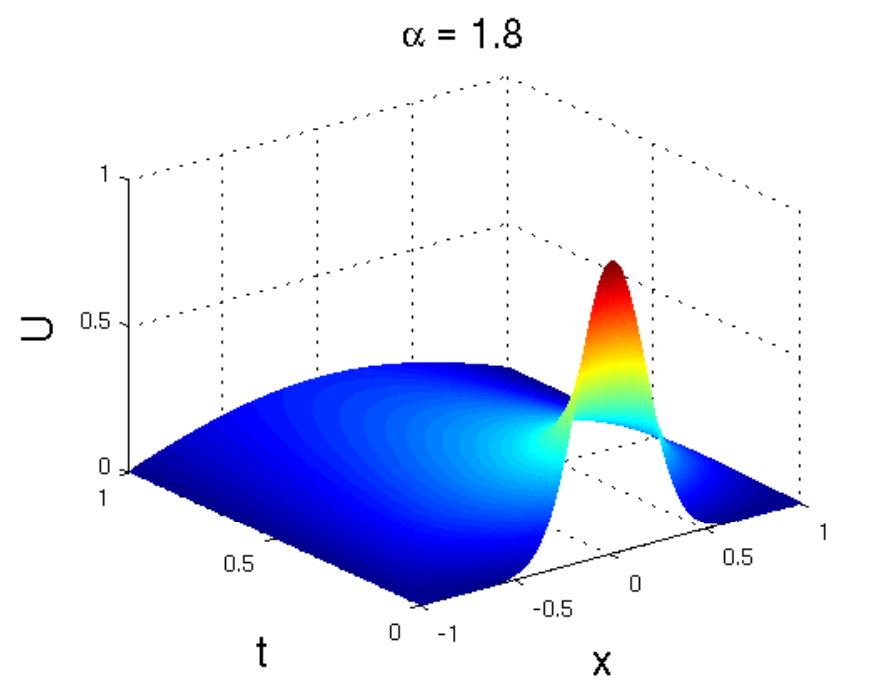}}
\centerline{
\includegraphics[height=3.960cm,width=5.260cm]{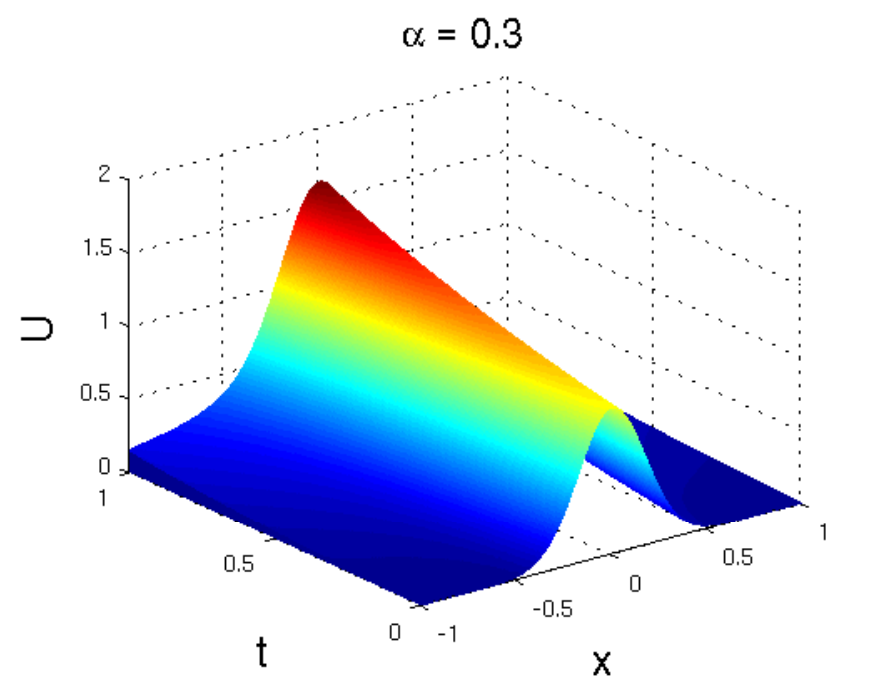}\hspace{-2mm}
\includegraphics[height=3.960cm,width=5.260cm]{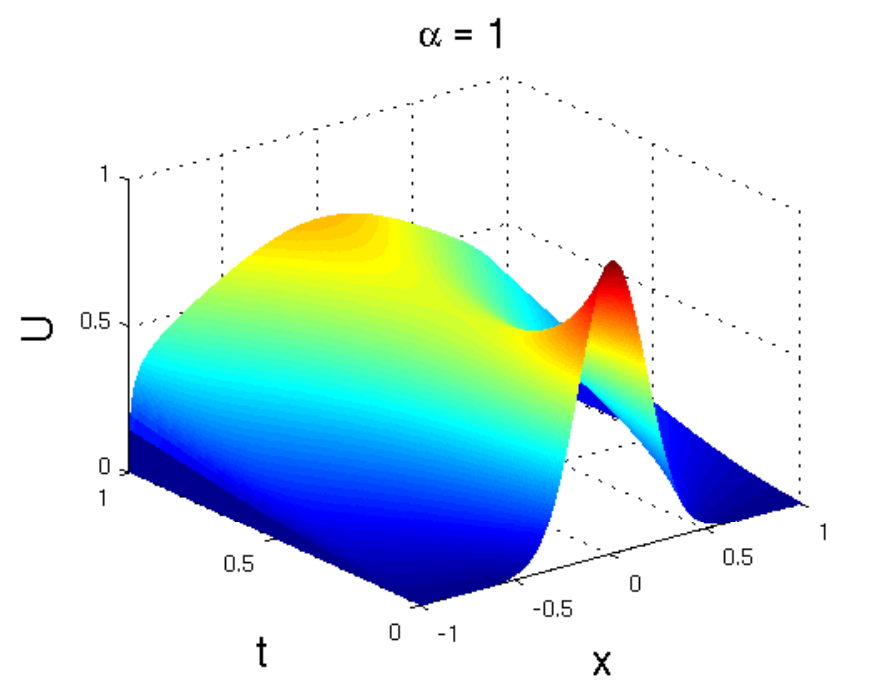}\hspace{-2mm}
\includegraphics[height=3.960cm,width=5.260cm]{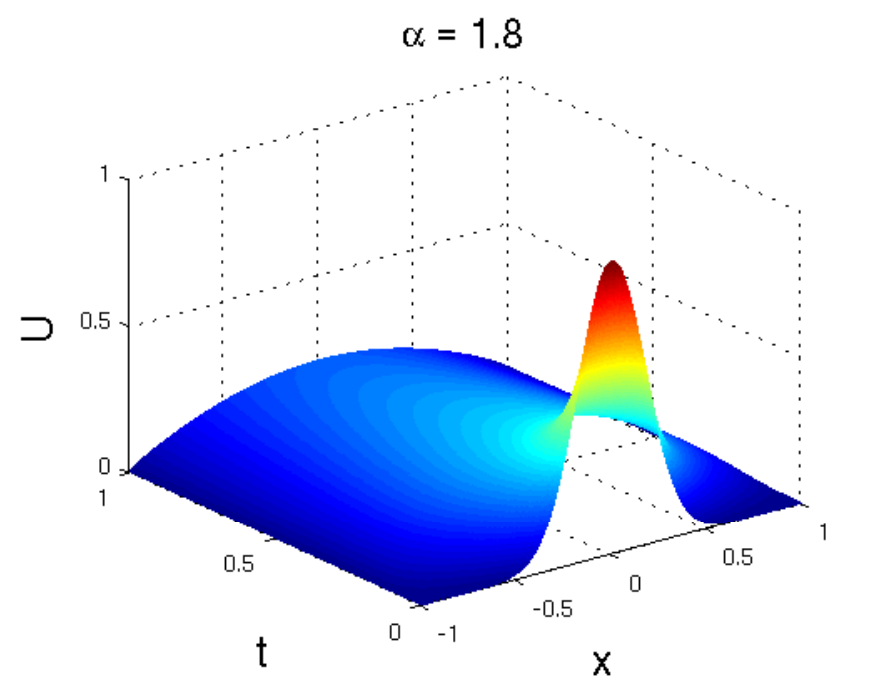}}
\caption{Time evolution of the solution $u(x, t)$ to the nonlocal diffusion-reaction equation (\ref{fDiffusion1}) with ${\mathcal L}_s$ (row one), ${\mathcal L}_h$ (row two), and ${\mathcal L}_r$ (row three).}\label{fig13}
\end{figure}
 (\ref{fDiffusion1}) with the spectral fractional Laplacian ${\mathcal L}_s$, fractional Laplacian ${\mathcal L}_h$, or regional fractional Laplacian ${\mathcal L}_r$, while the results for the peridynamic operator ${\mathcal L}_p$ are displayed in Figure \ref{fig14}. 
\begin{figure}[htb!]
\centerline{
\includegraphics[height=3.960cm,width=5.260cm]{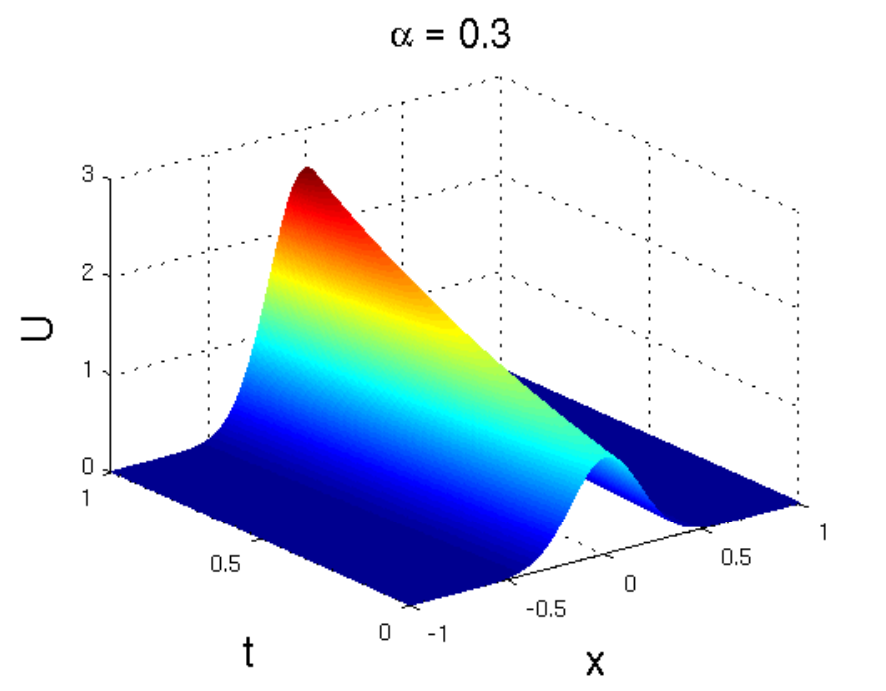}\hspace{-2mm}
\includegraphics[height=3.960cm,width=5.260cm]{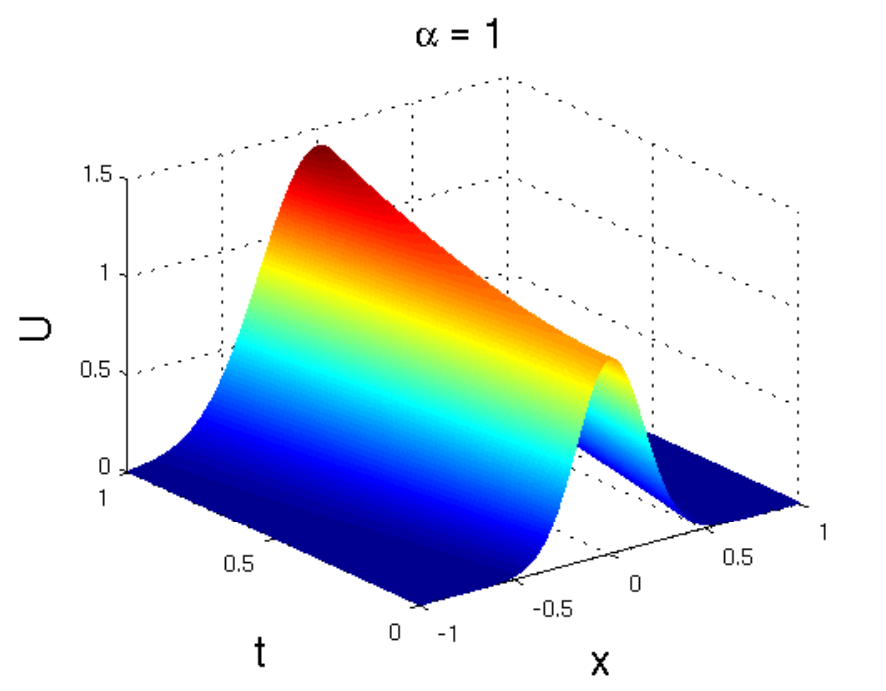} \hspace{-2mm}
\includegraphics[height=3.960cm,width=5.260cm]{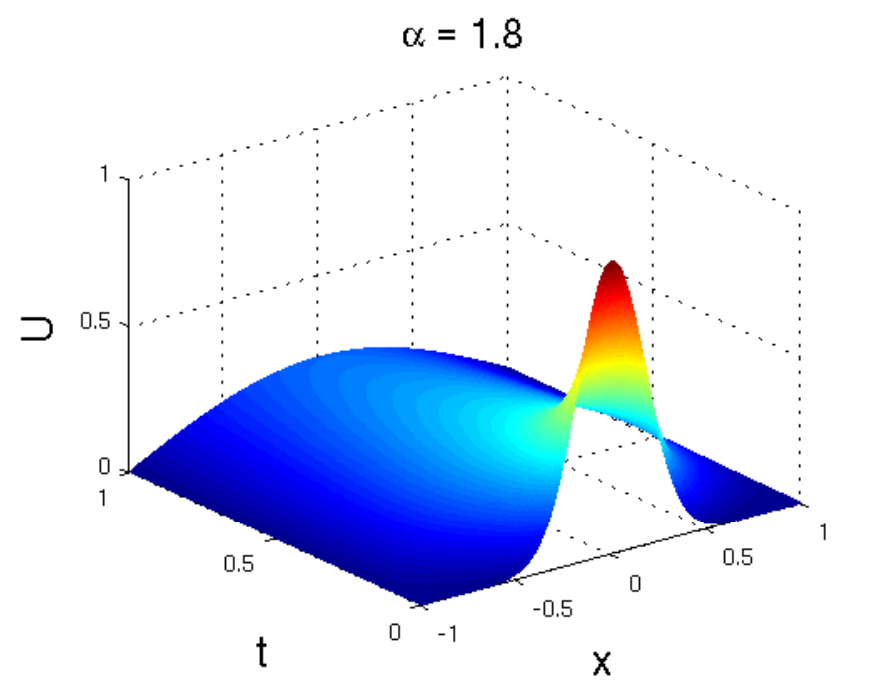}}
\centerline{
\includegraphics[height=3.960cm,width=5.260cm]{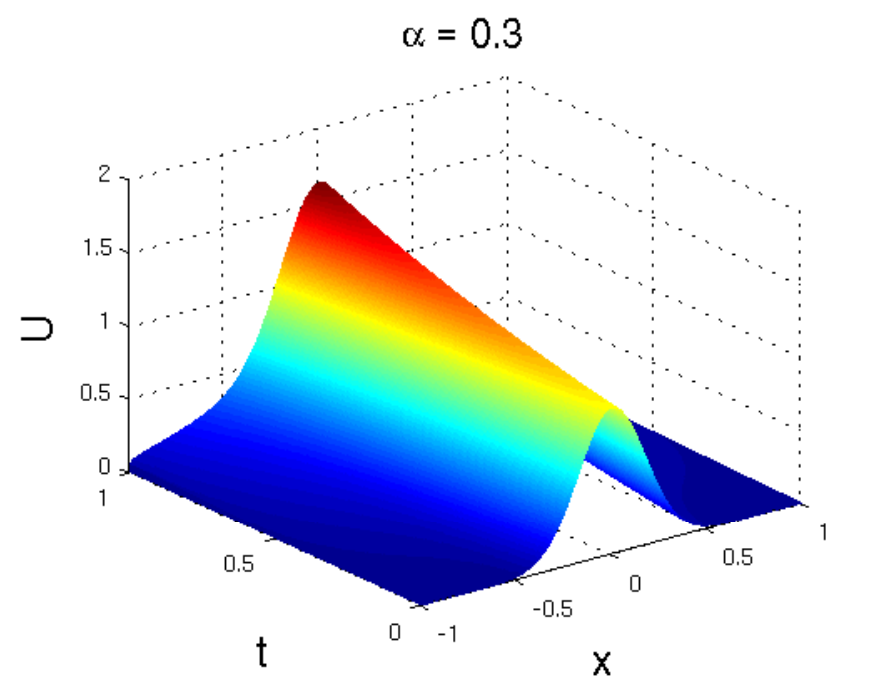}\hspace{-2mm}
\includegraphics[height=3.960cm,width=5.260cm]{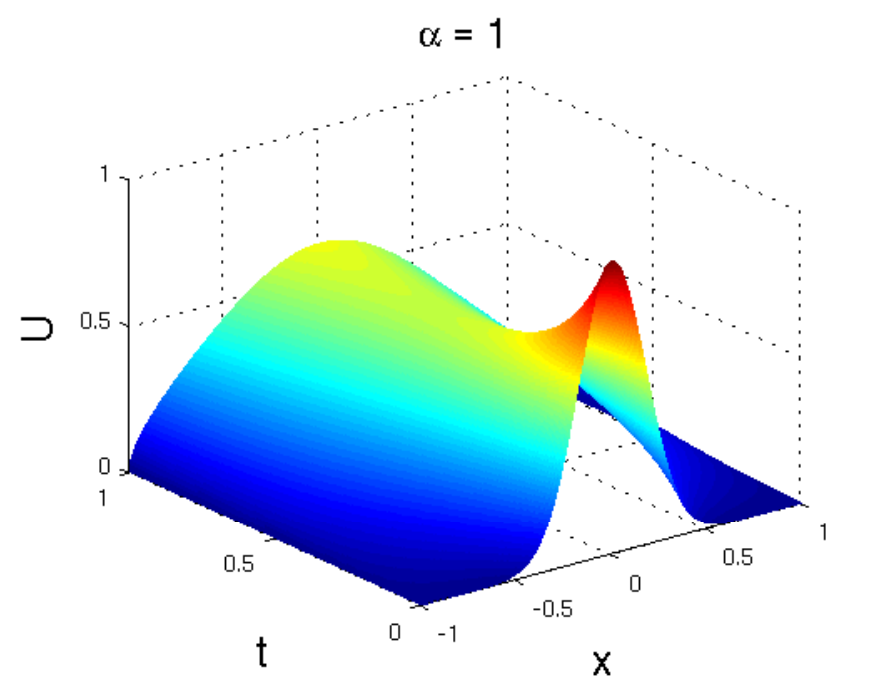} \hspace{-2mm}
\includegraphics[height=3.960cm,width=5.260cm]{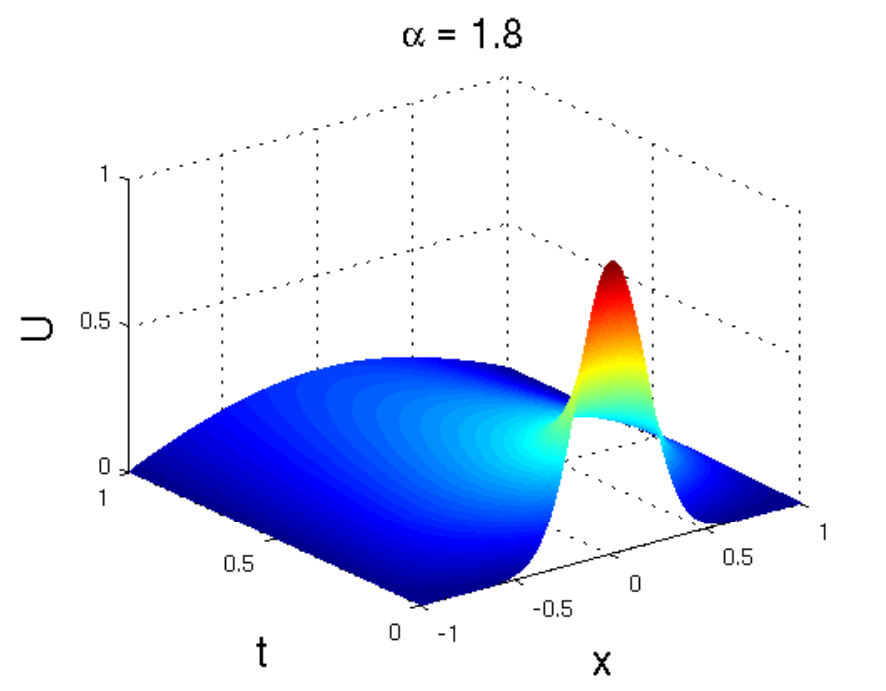}}
\caption{Time evolution of the solution $u(x, t)$ to the nonlocal diffusion-reaction equation (\ref{fDiffusion}) with the peridynamic operator ${\mathcal L}_p$,  where the horizon size $\dt = 0.1$ (top) or $\dt = 0.5$ (bottom).}\label{fig14}
\end{figure}
The results demonstrate the difference between these nonlocal operators, especially when $\ap$ is small. 
The solution from the fractional Laplacian ${\mathcal L}_h$ and the spectral fractional Laplacian ${\mathcal L}_s$ continuously decay over time for all the $\ap$ considered here. 
In contrast, the solution from the regional fractional Laplacian ${\mathcal L}_r$ or the peridynamic operator ${\mathcal L}_p$ may increase over time, depending on the competition between the diffusion and reaction terms. 
For example,  for $\ap = 0.3$, the solutions increases constantly, but as $\ap$ increases (e.g., $\ap = 1.8$), the diffusion becomes dominant, and the solution decays over time; see Figure \ref{fig16}.
\begin{figure}[htb!]
\centerline{
\includegraphics[height=4.060cm,width=5.960cm]{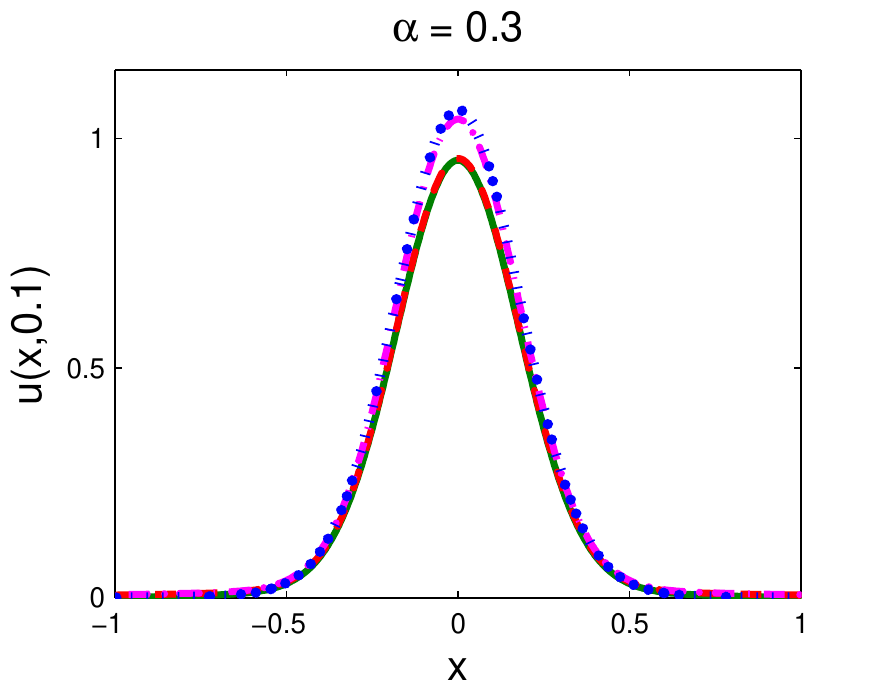}\hspace{-5mm}
\includegraphics[height=4.060cm,width=5.96cm]{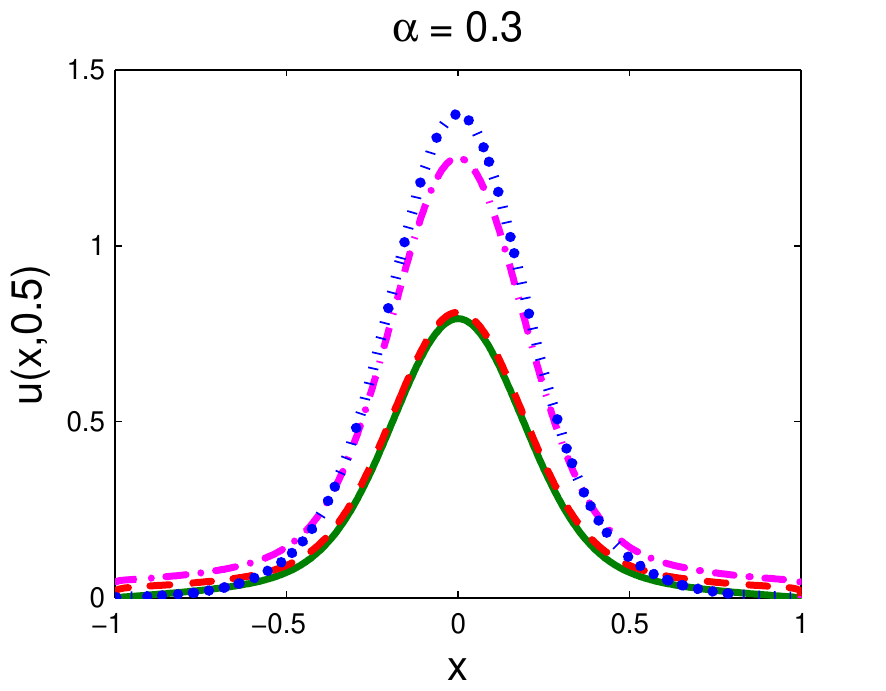}\hspace{-5mm}
\includegraphics[height=4.060cm,width=5.960cm]{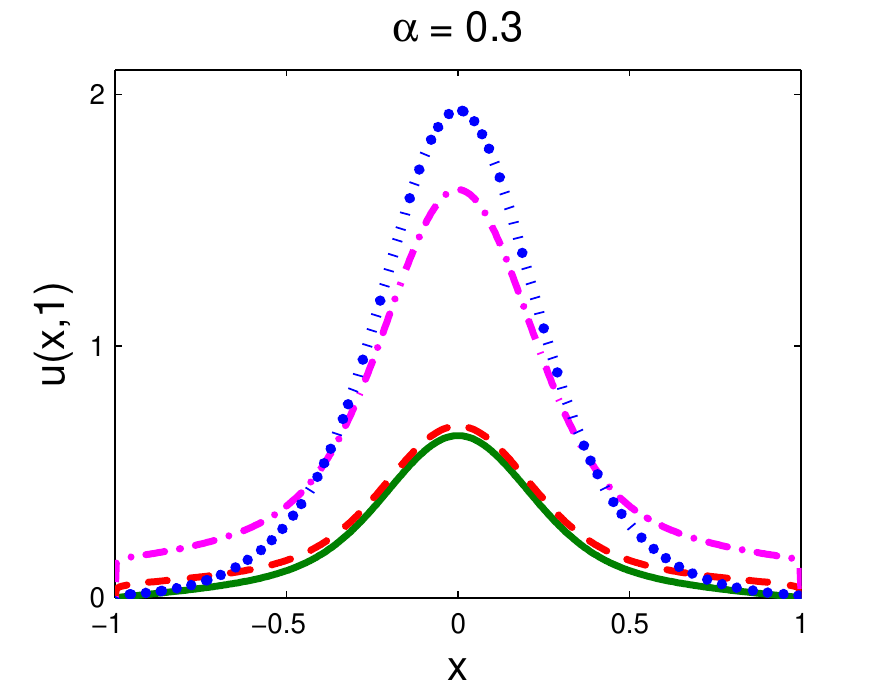}}
\centerline{
\includegraphics[height=4.060cm,width=5.960cm]{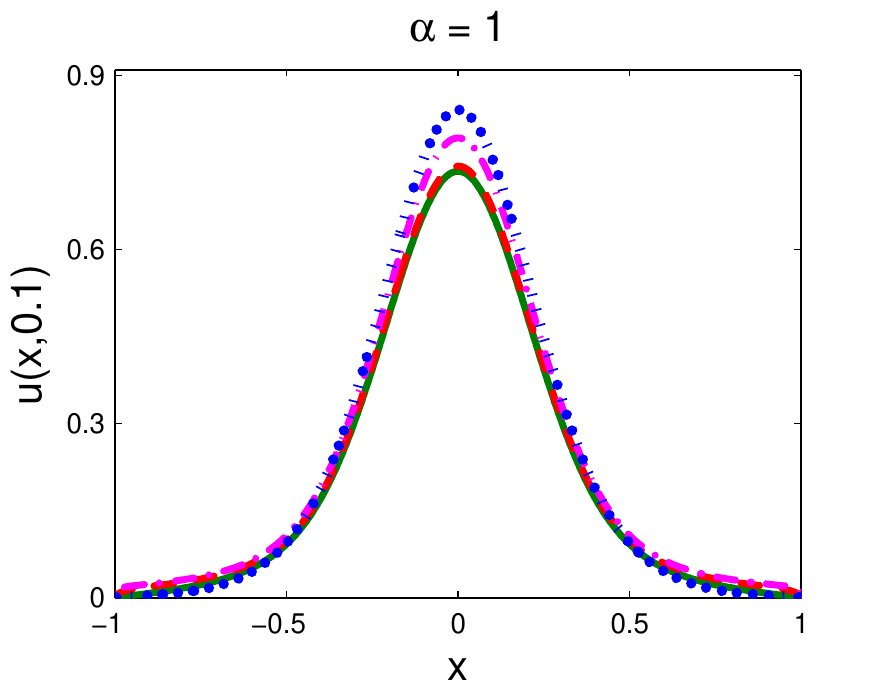}\hspace{-5mm}
\includegraphics[height=4.060cm,width=5.960cm]{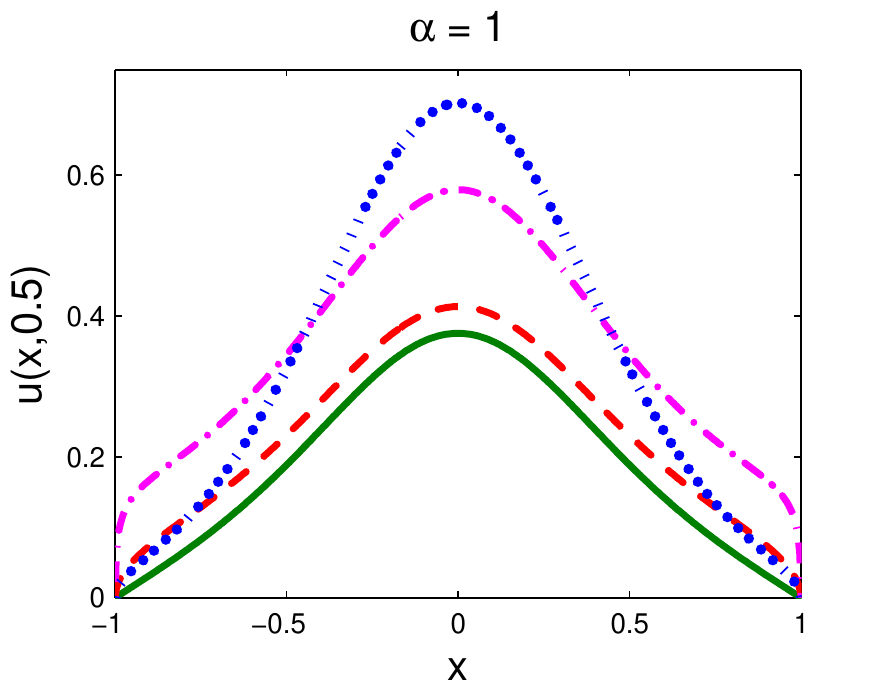} \hspace{-5mm}
\includegraphics[height=4.060cm,width=5.960cm]{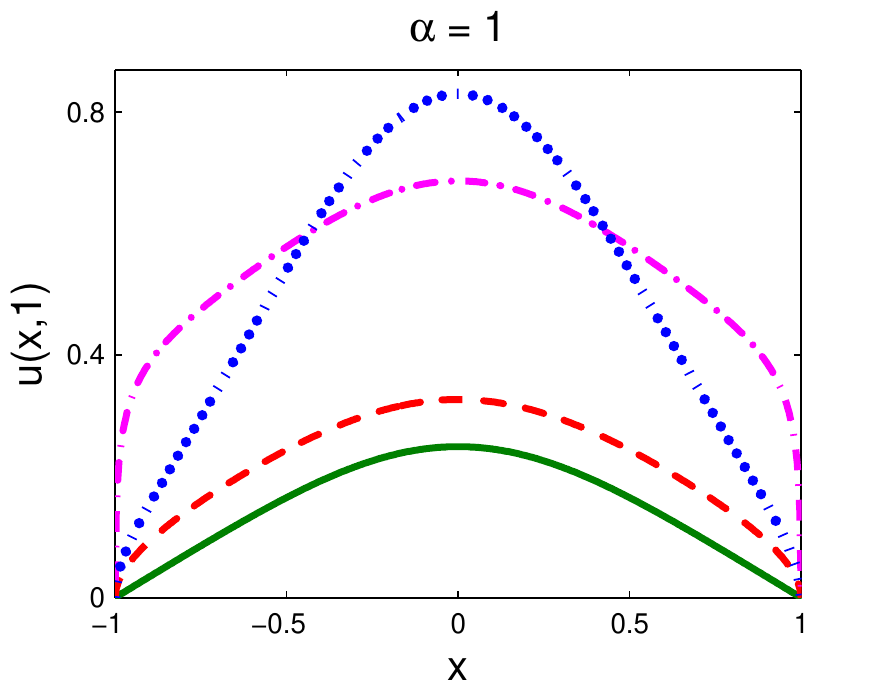}}
\centerline{
\includegraphics[height=4.060cm,width=5.960cm]{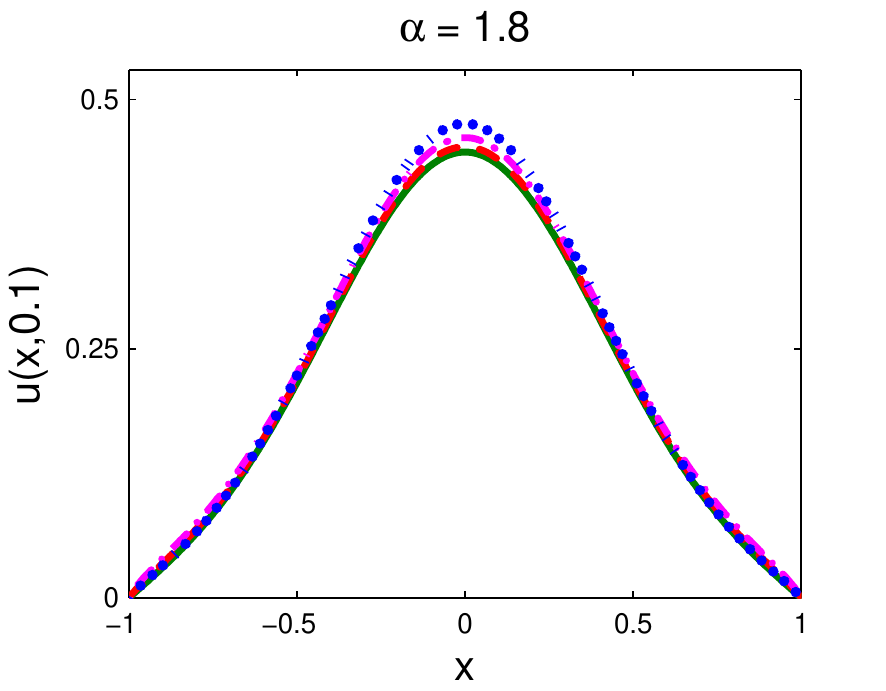}\hspace{-5mm}
\includegraphics[height=4.060cm,width=5.960cm]{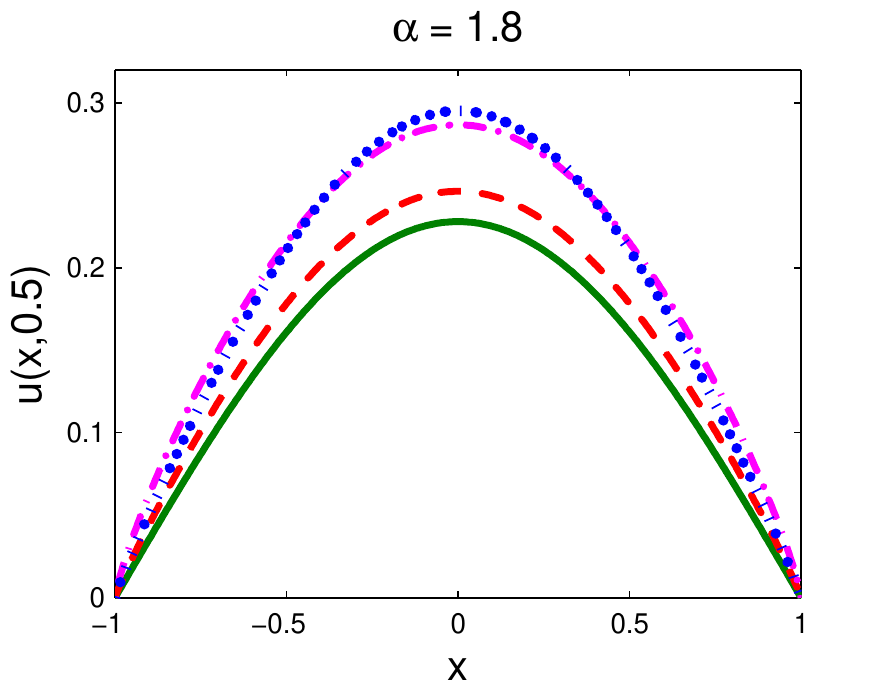} \hspace{-5mm}
\includegraphics[height=4.060cm,width=5.960cm]{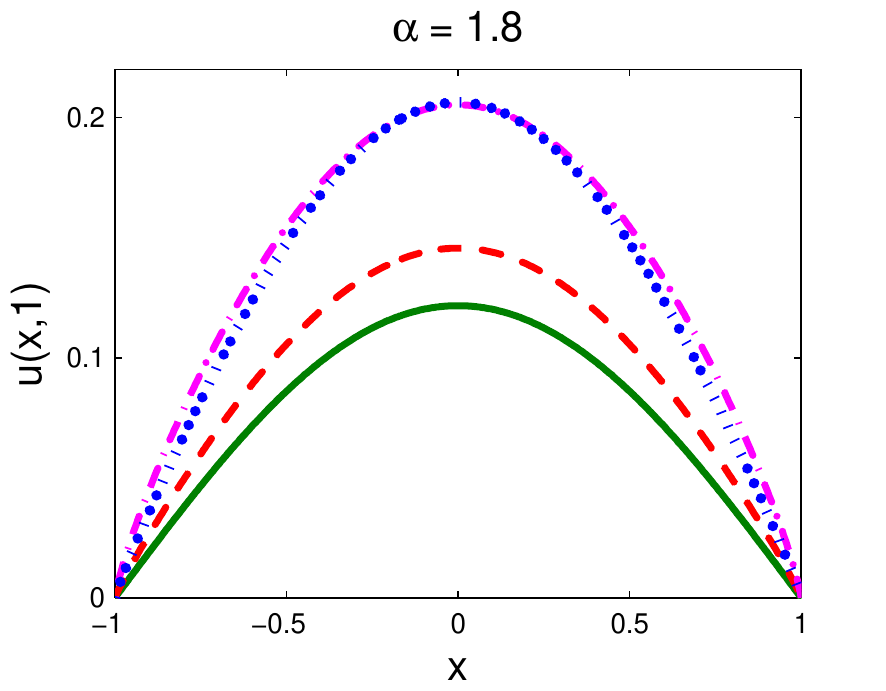}}
\centerline{
\includegraphics[height=4.060cm,width=5.960cm]{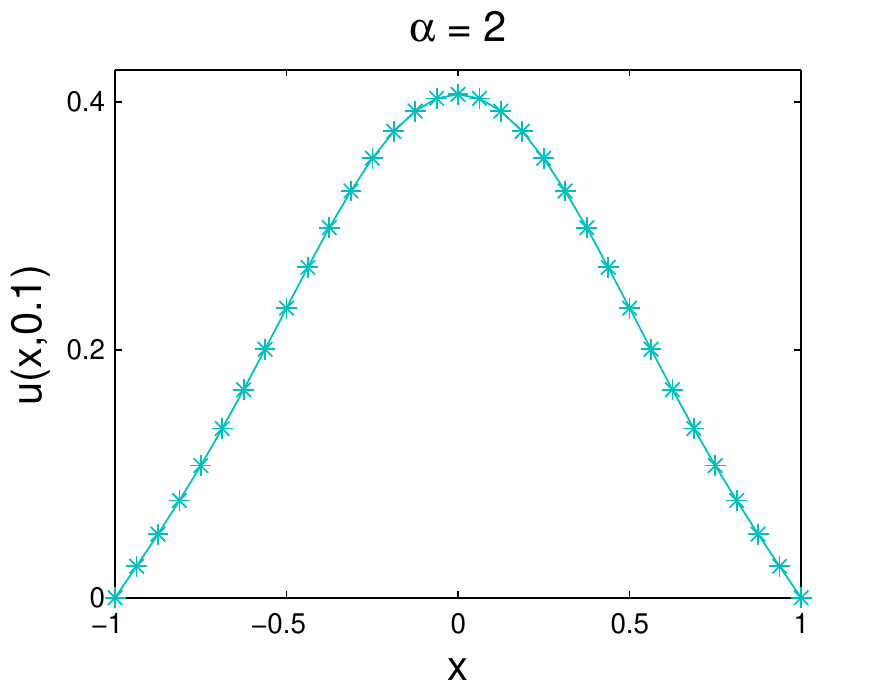}\hspace{-5mm}
\includegraphics[height=4.060cm,width=5.960cm]{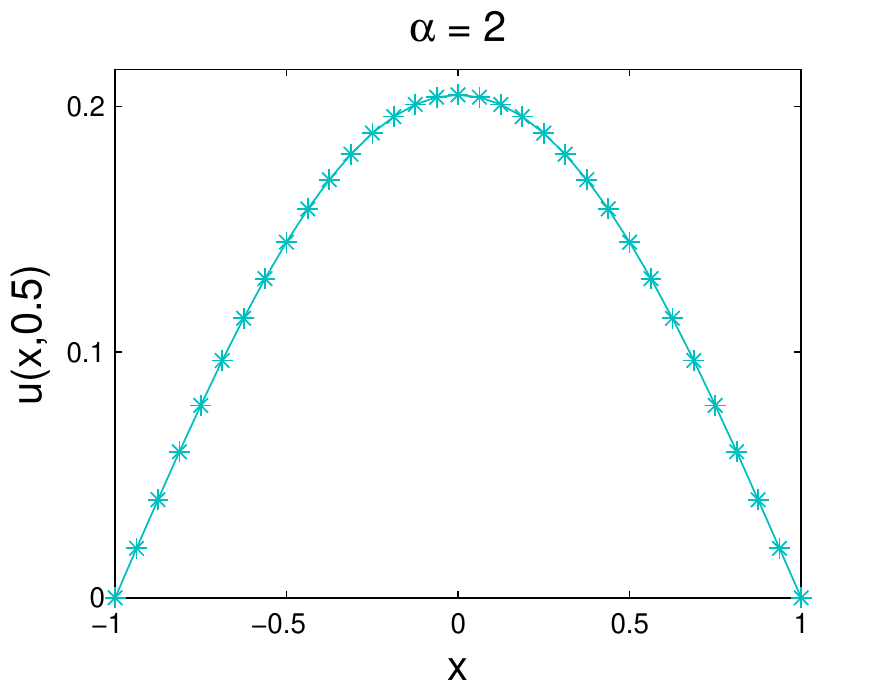} \hspace{-5mm}
\includegraphics[height=4.060cm,width=5.960cm]{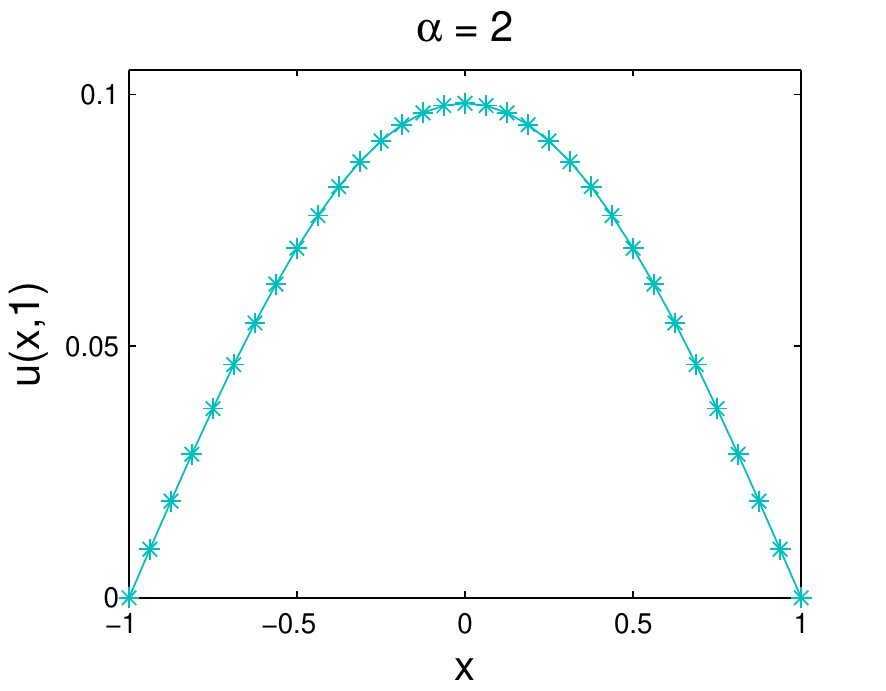}}
\caption{Solutions of the nonlocal diffusion-reaction equation (\ref{fDiffusion1}) at time $t = 0.1, 0.5, 1$, where the operator is chosen as ${\mathcal L}_s$ (solid line),  ${\mathcal L}_h$ (dashed line), ${\mathcal L}_p$ (dotted line), or ${\mathcal L}_r$ (dash-dot line). For easy comparison, we include the solution of the classical diffusion equation (i.e., ${\mathcal L}_i = -\p_{xx}$ in (\ref{fDiffusion1})) in the last row.}\label{fig16}
\end{figure}

\section{Concluding remarks}
\label{section4}

We studied four nonlocal diffusion operators, including the fractional Laplacian, spectral fractional Laplacian, regional fractional Laplacian,  and  peridynamic operator. 
These four operators are equivalent on the entire space ${\mathbb R}^n$, but their differences on a bounded domain are significant.
On a bounded domain $\Og \in {\mathbb R}^n$, they represent the generators of different stochastic processes. 
The Dirichlet fractional Laplacian represents the infinitesimal generator of a symmetric $\ap$-stable L\'evy process that particles are killed upon leaving the domain;  the regional  fractional Laplacian is the generator of a censored $\ap$-stable process that is obtained from a symmetric $\alpha$-stable L\'evy process by restricting its measure to $\Og$, while the spectral fractional Laplacian is the generator of  a subordinate killed Brownian motion (i.e. the process that first kills Brownian motion in a bounded domain $\Og$ and then subordinates it via a $\ap$-stable subordinator).
Our studies clarify the confusion existing in some literature, where the Dirichlet fractional Laplacian, spectral fractional Laplacian, and  regional fractional Laplacian are freely interchanged. 

We carried out extensive numerical investigations to understand and compare the nonlocal effects of these operators on a  bounded domain. 
Our numerical results suggest that: \vspace{-5pt}
\begin{enumerate}\itemsep -2pt
\item[i)] These four operators collapse to the classical Laplace operator as $\ap \to 2$.
\item[ii)] The eigenvalues and eigenfunctions of these four operators are different, although they all converge to those of the classical Laplace operator as $\ap \to 2$. 
For  each $k \in {\mathbb N}$, the eigenvalues of the spectral fractional Laplacian are always larger than those of the fractional Laplacian and regional fractional Laplacian, which numerically extends the conclusion in the literature \cite{Servadei2014}. 

\item[iii)]  For any $\ap \in (0, 2)$, the peridynamic operator can provide a good approximate to the fractional Laplacian, if the horizon size $\dt$ is sufficiently large.
We  found that the solution of the peridynamic model converges to that of the fractional Laplacian model at a rate of ${\mathcal O}(\dt^{-\ap})$.
Although the regional fractional Laplacian can be used to approximate the fractional Laplacian as $\ap \to 2$, it generally provides inconsistent results from the fractional Laplacian for $\ap \ll 2$.
\end{enumerate}
Moreover, we provided some conjectures from our numerical results, which might contribute to the mathematics analysis on these operators. 
Due to their nonlocality, numerical simulations of problems involving these four operators are considerably  challenging. 
We refer the readers to \cite{DuoJuZhang0016, DuoWykZhang, Delia2013, Du2012, Guan2017} for more discussions on numerical methods for these nonlocal models.

\vskip 15pt
\noindent{\bf Acknowledgements }
The second author was funded by the OSD/ARO MURI Grant W911NF-15-1-0562, and by the National Science Foundation under Grants DMS 1216923 and DMS-1620194.
The  third author was supported by the National Science Foundation under Grant DMS-1217000 and DMS-1620465.  Many thanks to the anonymous reviewers.

\bibliographystyle{plain}

\begin{landscape}
\begin{table}
\begin{center}
\begin{tabular}{|c|*{20}c|c|}
\hline
\diagbox[width=3em]{$k$}{$\ap$}  & 0.2 & & 0.5 & & 0.7 && 0.9&& 1 && 1.2 && 1.5 && 1. 8 && 1.95  && 1.999 && \\
\hline
\multirow{3}{*}{1} 
&1.0945 && 1.2533 && 1.3718 && 1.5014  && 1.5708 && 1.7193 && 1.9687 && 2.2543 && 2.4123 && 2.4663 && \multirow{3}{*}{2.4674} \\
& 0.9575 && 0.9702 && 1.0203 && 1.1032 && 1.1578 && 1.2971 && 1.5976 && 2.0488 && 2.3520 && 2.4650 && \\
& 0.0003 && 0.0038 && 0.0170 && 0.0640 && 0.1135 && 0.2939 && 0.8088  && 1.6602 && 2.2444 && 2.4628 && \\
\hline
\multirow{3}{*}{2}
& 1.2573 && 1.7725  && 2.2285 && 2.8018 && 3.1416 && 3.9498 && 5.5683  &&7.8500  && 9.3206 && 9.8583 &&\multirow{3}{*}{9.8696} \\
& 1.1966 && 1.6016 && 1.9733 && 2.4583 && 2.7549 &&  3.4870 && 5.0600 && 7.5033 && 9.2082 && 9.8559 &&\\
& 0.1878 && 0.4593 && 0.6729 && 0.9799 && 1.2026 && 1.8719 && 3.6509 && 6.7378 && 8.9854  && 9.8512 && \\
\hline
\multirow{3}{*}{3} 
&1.3635  && 2.1708 && 2.9598 &&  4.0357 && 4.7124 && 6.4252 && 10.230  && 16.287 && 20.550 && 22.172 &&\multirow{3}{*}{22.207} \\
& 1.3191  && 2.0289 && 2.7294 && 3.6987 && 4.3171 && 5.9121 && 9.5948 && 15.800 && 20.384 && 22.169 && \\
& 0.3085  && 0.8626 && 1.3646 && 2.0823 && 2.5760 && 3.9902 && 7.7500 &&14.701 && 20.049 && 22.161 &&\\
\hline
\multirow{3}{*}{4}
& 1.4442 && 2.5066 && 3.6201 && 5.2283 && 6.2832 && 9.0744 && 15.750 && 27.335 && 36.012 && 39.406 && \multirow{3}{*}{29.478} \\
& 1.4106 && 2.3873 && 3.4131 &&  4.9055 && 5.8925 && 8.5350 && 15.020 && 26.725 && 35.794 && 39.401 && \\
& 0.3981 && 1.2091 && 2.0140 && 3.2054 && 4.0292 && 6.3902 && 12.811 && 25.313  && 35.349 && 39.391 &&\\
\hline
\multirow{3}{*}{5} 
& 1.5101 && 2.8025 && 4.2322 && 6.3912 && 7.8540 && 11.861 && 22.011 && 40.847 && 55.645 && 61.558 &&\multirow{3}{*}{61.685} \\
& 1.4817 && 2.6949 && 4.0371 && 6.0733 && 7.4607 && 11.293 && 21.191 &&  40.115 && 55.374 && 61.552 &&\\
&  0.4700 && 1.5149 && 2.6231 && 4.3230 && 5.5171 && 8.9817 && 18.670 && 38.408 && 54.820 && 61.540 && \\
\hline
\multirow{3}{*}{6}
& 1.5662 && 3.0700 && 4.8083 && 7.5309 && 9.4248 && 14.761 && 28.934 && 56.714 && 79.402 && 88.627 &&\multirow{3}{*}{88.826} \\
& 1.5422  && 2.9730 && 4.6253 && 7.2206 && 9.0334 && 14.175 && 28.037 && 55.868 && 79.080 && 88.620 &&\\
&0.5306  && 1.7911 &&  3.1993&& 5.4300  && 7.0245 && 11.722 && 25.235 && 53.876 && 78.418  && 88.605 && \\
\hline
\multirow{3}{*}{8}
& 1.6590 && 3.5449  && 5.8809 && 9.7564 && 12.566 && 20.847 && 44.547 &&  95.187 && 139.14 && 157.51 &&\multirow{3}{*}{157.91}\\
& 1.6400  && 3.4612 && 5.7133 && 9.4550 && 12.175 && 20.225 && 43.509 && 94.122 && 138.72 && 157.50 && \\
& 0.6296 && 2.2799 && 4.2751 &&  7.6101 && 10.072 && 17.552 && 40.218 && 91.591 && 137.85  && 157.49 &&\\
\hline
\multirow{3}{*}{10}
& 1.7347  && 3.9633 && 6.8752 && 11.926 && 15.708 && 27.249 && 62.256 && 142.24 && 215.00 && 246.06 &&\multirow{3}{*}{246.74} \\
& 1.7189 && 3.8886 && 6.7186 && 11.632 && 15.317 && 26.598 &&  61.096 && 140.96 && 214.48 && 246.05 &&\\
& 0.7095 && 2.7090 && 5.2735 && 9.7490 && 13.145 && 23.749 && 57.377 &&137.92 && 213.39 && 246.02 && \\
\hline
\multirow{3}{*}{20}
& 1.9926 && 5.6050  && 11.169 && 22.255 && 31.416 && 62.601 && 176.09 && 495.30 && 830.70 && 983.56 &&\multirow{3}{*}{986.96 }\\
& 1.9836 && 5.5525 && 11.042 && 21.981 && 31.025 && 61.854&& 174.45 && 493.09&& 829.69 && 983.53 &&\\
& 0.9779 && 4.3810 && 9.5850&& 19.998 && 28.657 && 58.439 && 169.09 && 487.74&& 827.58 && 983.49 &&\\
\hline
\end{tabular}\caption{Comparison of the eigenvalues for different operators, where the eigenvalues of the standard Dirichlet Laplace operator $-\Dt$ are presented in most right column. For each $k$,  upper row: $\lambda_k^s$;  middle row:  $\lambda_k^h$; lower row: $\lambda_k^r$.}\label{Tab1} \end{center}
\end{table}
\end{landscape}

\end{document}